\documentclass[10pt]{article}

\usepackage{amsmath}
\usepackage{amsfonts}
\usepackage{amssymb}
\usepackage{framed}
\usepackage{mathtools}
\usepackage{enumerate}
\usepackage{mathrsfs}
\usepackage[hidelinks]{hyperref}
\usepackage{enumitem}
\usepackage{mathrsfs}
\usepackage{scrextend}
\usepackage{amsthm}
\usepackage{bbm}
\usepackage{thmtools}
\usepackage{thm-restate}
\usepackage{mathabx}
\usepackage[margin=1.5in]{geometry}
\usepackage{rotating}

\usepackage{xcolor}
\hypersetup{
    colorlinks,
    linkcolor={red!50!black},
    citecolor={blue!50!black},
    urlcolor={blue!80!black}
}

\usepackage{tikz-cd}

\theoremstyle{plain}
\newtheorem{thm}{Theorem}[section]
\newtheorem{lem}[thm]{Lemma}
\newtheorem{prop}[thm]{Proposition}
\newtheorem{cor}[thm]{Corollary}

\theoremstyle{definition}
\newtheorem{defn}[thm]{Definition}

\theoremstyle{remark}
\newtheorem*{rem}{Remark}

\newtheorem*{notn}{Notation}

\tikzset{
  symbol/.style={
    draw=none,
    every to/.append style={
      edge node={node [sloped, allow upside down, auto=false]{$#1$}}}
  }
}

\newcommand{\Spec}{\textrm{Spec} \hspace{0.15em} }

\newcommand\restr[2]{{
	\left.\kern-\nulldelimiterspace
	#1
	\vphantom{\big|}
	\right|_{#2}
	}}
\newcommand{\an}[1]{#1^{\textrm{an}}}

\newcommand{\ep}{\varepsilon}

\newcommand{\Sch}{\textrm{Sch}}

\newcommand{\Hom}{\textrm{Hom}}

\newcommand{\ch}[1]{\widecheck{{#1}}}

\newcommand{\codim}{\textrm{codim}}

\newcommand{\GL}{\textrm{GL}}

\newcommand{\rank}{\textrm{rank}\,}

\usepackage{amsmath,calligra,mathrsfs}

\usepackage[cal=boondoxo]{mathalfa}

\DeclareMathOperator{\sheafhom}{\mathcal{H \kern -1pt o \kern -2pt m}}
\DeclareMathOperator{\sheafend}{\mathcal{E \kern -1pt n \kern -2pt d}}
\DeclareMathOperator{\sheafaut}{\mathcal{A \kern -1pt u \kern -2pt t}}

\title{Algebraic Cycle Loci at the Integral Level}
\author{David Urbanik}

\begin{document}

\maketitle

\begin{abstract}
Let $f : X \to S$ be a smooth projective family defined over $\mathcal{O}_{K}[\mathcal{S}^{-1}]$, where $K \subset \mathbb{C}$ is a number field and $\mathcal{S}$ is a finite set of primes. For each prime $\mathfrak{p} \in \mathcal{O}_{K}[\mathcal{S}^{-1}]$ with residue field $\kappa(\mathfrak{p})$, we consider the algebraic loci in $S_{\overline{\kappa(\mathfrak{p})}}$ above which cohomological cycle conjectures predict the existence of non-trivial families of algebraic cycles, generalizing the Hodge loci of the generic fibre $S_{\overline{K}}$. We develop a technique for studying all such loci, together, at the integral level. As a consequence we give a non-Zariski density criterion for the union of non-trivial ordinary algebraic cycle loci in $S$. The criterion is quite general, depending only on the level of the Hodge flag in a fixed cohomological degree $w$ and the Zariski density of the associated geometric monodromy representation.
\end{abstract}

\tableofcontents

\section{Introduction}
\label{intro}

Let $f : X \to S$ be a smooth projective family defined over $R$, where $R \subset \mathbb{C}$ is a Noetherian subring, and $S$ is smooth and quasi-projective. The purpose of this manuscript is to develop a differential-algebraic method for the systematic study of loci in $S$ defined by families of algebraic cycles associated to $f$. To explain what we mean, we will write $f^{\otimes n}$ for the map
\[ X^{\otimes n} = \underbrace{X \times_{S} \cdots \times_{S} X}_{n} \to S , \]
obtained by taking the $n$-fold fibre product of $f$ with itself. By a family of algebraic cycles over a locally closed subscheme $Z \subset S$ we then mean a locally closed subscheme $Y \subset X^{\otimes n}$ which is smooth over $Z$. We are particularly interested in the situation where cohomological realizations of the fibres of $Y \to Z$ induce non-trivial classes in either Betti, algebraic de Rham, or crystalline cohomology.

As is in many ways well known, the existence of such a $Y$ lying over $Z$ can induce differential constraints on $Z$. For instance, if one considers the same situation for the associated complex-algebraic family $f_{\mathbb{C}} : X_{\mathbb{C}} \to S_{\mathbb{C}}$, and with $Z \subset S_{\mathbb{C}}$ and $Y \subset X^{\otimes n}_{\mathbb{C}}$ smooth complex algebraic subvarieties, then a finite index subgroup of $\pi_{1}(\an{Z}, s)$ will fix the Betti cycle class of $Y_{s}$ inside the cohomology of $X^{\otimes n}_{s}$, and this subgroup is independent of the point $s \in Z(\mathbb{C})$. Using the Riemann-Hilbert correspondence one can understand this fact in terms of the existence of a global flat section for an algebraic vector bundle with flat connection defined on $Z$. These vector bundles and connections admit natural descriptions in terms of the algebraic de Rham cohomology of the family $f$, and make sense at the integral level. Moreover, as we review in \S\ref{crysflatsec}, the same sort of reasoning can be carried out at the crystalline level in the case where $Z \subset S$ and $Y \subset X^{\otimes n}$ are positive characteristic subschemes. A natural question that arises is then to understand the relationship between such ``global'' differential constraints in characteristic zero and ``local'' ones that occur at finite primes; for instance, when can one conclude that all local differential constraints of this type arise as the reduction to positive characteristic of global ones?

That one cannot expect something so strong in general is clear, as it is easy to construct families of algebraic varieties that acquire additional algebraic cycles upon reduction modulo a prime. Moreover if one thinks in terms of the differential equations themselves, it is clear that spaces of solutions can degenerate upon reduction. However one could hope that such pathological behaviour remains confined to the fibre above a proper closed subscheme of $\Spec R$, and that away from these primes such differential constraints on subschemes $Z \subset S$ are controlled by corresponding constraints on subschemes lying above the generic fibre. In this paper we show that, at least when the global monodromy on $S$ is sufficiently large and we consider cohomological objects of level at least three, this is exactly what happens.

\vspace{0.5em}

The reason answering such questions is difficult, and is not merely an application of some kind of Lefschetz principle, is that we do not impose any limitations on the loci $Z \subset S$ that we consider, and in effect study all of them at once. What this means it that one cannot appeal to the fact that the loci of algebraic cycles under consideration belong to some bounded family and then employ algebro-geometric machinery to reduce the arithmetic situation away from a proper closed subscheme of $\Spec R$ to the geometric situation at the generic fibre. Instead, what is needed is an a priori way of controlling all such ``differential constraints'' of interest simultaneously at the integral level.

A possible approach to obtaining the kind of control we need is to use the theory of period maps, as is for instance illustrated in the recent work of \cite{BKU}. The idea is as follows. One considers the torsion-free variation of Hodge structure $\mathbb{V} = R^{i} \an{f_{\mathbb{C},*}} \mathbb{Z}/\textrm{tor.}$ arising from $f$ in cohomological degree $i$ and fixes a polarization $\mathcal{Q} : \mathbb{V} \otimes \mathbb{V} \to \mathbb{Z}$ on $\mathbb{V}$. One additionally fixes a polarized integral lattice $(V, Q)$ isomorphic to one (hence any) fibre $(\mathbb{V}_{s}, \mathcal{Q}_{s})$, defines $D$ to be the complex manifold parametrising all Hodge structures on $V$ polarized by $Q$, and sets $\Gamma = \textrm{Aut}(V,Q)(\mathbb{Z})$. Then we obtain a natural analytic period map $\varphi : \an{S_{\mathbb{C}}} \to \Gamma \backslash D$ which sends a point $s \in S(\mathbb{C})$ to the isomorphism class of the polarized Hodge structure $(\mathbb{V}_{s}, \mathcal{Q}_{s})$. The loci $Z \subset S_{\mathbb{C}}$ where one has ``additional differential constraints'' can then be interpreted as irreducible components of analytic loci of the form $\varphi^{-1}(\Gamma' \backslash D')$, where $D' \subset D$ is a real orbit of a certain kind of ``Hodge-theoretic'' $\mathbb{Q}$-algebraic subgroup $\mathbf{H} \subset \textrm{Aut}(V, Q)$ and $\Gamma'$ is a subgroup of $\mathbf{H}(\mathbb{Q}) \cap \textrm{Aut}(V, Q)(\mathbb{Z})$. 

One is now interested in understanding the loci $Z$ that arise as the analytic subvarieties $\Gamma' \backslash D'$ vary, and the observation made in \cite{BKU} is that although infinitely many loci $\Gamma' \backslash D' \subset \Gamma \backslash D$ may arise, the orbits $D' \subset D$ that one is interested in come from only finitely many complex analytic families, and the complex analytic geometry of these families is ``tame'' in the sense made precise by the theory of o-minimal structures. Thus while one cannot appeal to the boundedness of some set of families containing the loci $Z \subset S_{\mathbb{C}}$ that can arise, one \emph{can} appeal to the boundedness of the families of loci $D' \subset D$ that define them, and through a careful understanding of the geometry of the period map $\varphi$ one can use this ``boundedness'' to constrain the geometry of the Hodge locus in $S_{\mathbb{C}}$.

\vspace{0.5em}

To make this strategy work at an integral level we will need to develop a theory of period maps which is capable of functioning in a purely algebro-geometric setting, even over a finite field. For this we will introduce a tool for studying period maps ``infinitesimally'', building on our earlier work in \cite{periodimages} and \cite{urbanik2021sets}. The approach we will give, although similar in spirit to the approach given by the authors in \cite{BKU}, is entirely independent, as \cite{BKU} relies on an analytic theory of period maps not available in positive characteristic.

\subsection{Non-Density Criteria}
\label{nondenintro}

The main application of our ideas will be to give a non-Zariski destiny criterion for the locus where one has ``extra non-trivial families'' of primitive algebraic cycles. We assume we have $f$, $S$, $X$, $f^{\otimes n}$ and $X^{\otimes n}$ as above. We consider the variation of Hodge structure $\mathbb{V}_{\textrm{full}} = R^{k} f_{*} \mathbb{Z} / \textrm{tor}.$ modulo torsion in degree $k$ corresponding to the family $f$, let $\mathbb{V} \subset \mathbb{V}_{\textrm{full}}$ be its primitive subsystem, and let $\mathcal{Q} : \mathbb{V} \otimes \mathbb{V} \to \mathbb{Z}$ be the natural polarization induced by the Lefschetz decomposition and cup product (see \S\ref{primcupsec}). We also let $\mathcal{H}$ be the primitive subsheaf of the algebraic de Rham cohomology sheaf $R^{k} f_{*} \Omega^{\bullet}_{X/S}$, let $F^{\bullet}$ be its Hodge filtration constructed as in \S\ref{recollsec}, and $\nabla$ the Gauss-Manin connection described in \cite{katz1968} (c.f. \S\ref{recollsec}). We note that although $R^{k} f_{*} \Omega^{\bullet}_{X/S}, \mathcal{H}$ and the constituent bundles in $F^{\bullet}$ are not vector bundles in general, they become so after replacing $\Spec R$ with an affine open subscheme by spreading out from the generic fibre, and we will typically assume that we have done so. (Our results will all allow for such a replacement, and will consequently be unaffected.) We also note that the formation of these sheaves is compatible with arbitrary base change, as a consequence of \cite[\href{https://stacks.math.columbia.edu/tag/0FM0}{Lemma 0FM0}]{stacks-project}. The period map $\varphi$ is defined as before with respect to $\mathbb{V}$. 

\subsubsection{The Case of Vectors}

To illustrate the key ideas, we first state a theorem in the setting of vectors, i.e. when $n = 1$, before turning to the general setting.

\begin{thm}
\label{babythm}
Suppose that $\varphi$ is quasi-finite, and that for some $s \in S(\mathbb{C})$ one has that:
\begin{itemize}
\item[-] the image of the monodromy representation
\[ \pi_{1}(\an{S}_{\mathbb{C}}, s) \to \textrm{Aut}(\mathbb{V}_{s}, \mathcal{Q}_{s}) , \]
is Zariski dense;
\item[-] the natural Hodge structure on the Lie algebra of $\textrm{Aut}(\mathbb{V}_{s}, \mathcal{Q}_{s})$ has level at least $6$.
\end{itemize}
Then there exists a closed subscheme $E$, properly contained in $S$, with the following property: for any diagram 
\begin{equation}
\label{cyclediag}
\begin{tikzcd}
Y \arrow[r, hook] \arrow[d, swap] & X \arrow[d, "f"] \\
Z \arrow[r, hook] & S
\end{tikzcd}
\end{equation}
along an inclusion $Z \subset S$ of a locally closed subscheme which is positive dimensional over a prime $\mathfrak{p} \in \Spec R$, and for which the fibres of $Y \to Z$ induce non-zero cohomology classes in the fibres of $\mathcal{H}_{\kappa(\mathfrak{p})}$, we have $Z \subset E$.
\end{thm}

\begin{rem}
That $Z$ be positive dimensional over $\mathfrak{p}$ means that $Z$ maps into the closure of $\mathfrak{p}$ in $\Spec R$ and that the map $\pi : Z \to \Spec R$ has relative dimension at least $1$ at every point $z \in Z$; i.e., for each such $z$ one has that $\dim \pi^{-1}(\pi(z)) \geq 1$ on the level of schemes. 
\end{rem}

\begin{rem}
We emphasize that $Y$ and $Z$ are \emph{not} required to lie over the generic fibre of $R$; in particular, if $R$ is the ring of integers of a number field, $Y$ and $Z$ might correspond to algebraic varieties over a finite field. 
\end{rem}

Let us explain the monodromy and Hodge level hypotheses. As explained in the introduction above, the monodromy over a subvariety $Z \subset S_{\mathbb{C}}$ is constrained by the presence of non-trivial families of algebraic cycles lying over $Z$. To quantify how many such families lie over $S_{\mathbb{C}}$, one can consider the group $\mathbf{H}_{S}$ defined as the identity component of the Zariski closure of the image of the map 
\[ \pi_{1}(\an{S}_{\mathbb{C}}, s) \to \textrm{Aut}(\mathbb{V}_{s}, \mathcal{Q}_{s}) . \]
The first hypothesis of \autoref{babythm} then says that $\mathbf{H}_{S} = \textrm{Aut}(\mathbb{V}_{s}, \mathcal{Q}_{s})$, which means that every family of algebraic cycles over $S_{\mathbb{C}}$ in fact comes from the invariant algebra associated to the cup product pairing; i.e., such algebraic cycles induce global sections of $\mathcal{H}^{a,b}_{\mathbb{C}} = \restr{\mathcal{H}^{\otimes a} \otimes (\mathcal{H}^{\otimes b})}{S_{\mathbb{C}}}$ for some $a, b$ that lie in the tensor algebra of sections associated to the orthogonal or sympletic pairing. (For a description of this algebra see \cite[\S1]{invtheoryrem}.) 

We note that the group $\mathbf{H}_{S}$, although defined in terms of the rational structure of the local system $\mathbb{V}$, in fact admits a natural ``realization'' inside $\textrm{Aut}(\mathcal{H}_{s})$ for any scheme-theoretic point $s \in S$. This is because $\mathbf{H}_{S,\mathbb{C}}$ can, by virtue of the Riemann-Hilbert correspondence, be viewed as the stabilizer of the spaces of global flat sections of $\mathcal{H}^{a,b}_{\mathbb{C}}$ for all $a, b \geq 0$, and these spaces of sections are in fact defined over the fraction field of $R$ by the argument in \cite[Lem 3.4]{fieldsofdef}, and thus can be defined over $R$ itself as the bundles $\mathcal{H}^{a,b} = \mathcal{H}^{\otimes a} \otimes (\mathcal{H}^{*})^{\otimes b}$ are $R$-algebraic. This means that we can always ask, even for a positive-characteristic family $Y \to Z$, whether or not the fibres of this family induce algebraic cycle classes fixed by $\mathbf{H}_{S,s} \subset \textrm{Aut}(\mathcal{H}_{s})$ for any $s \in Z$. We note that in this context, we do not have a way ``identifying'' the groups $\mathbf{H}_{S,s}$ as the point $s \in S$ varies, so we will typically make sure to specify the point in question. 


As for the Hodge level hypothesis, we begin by recalling that the Lie algebra $\mathfrak{h}_{S}$ of $\mathbf{H}_{S}$ admits a natural Hodge decomposition. First, one views the Hodge structure on the fibre $\mathbb{V}_{s}$ as a morphism $h : \mathbb{S} \to \GL(\mathbb{V}_{s})_{\mathbb{R}}$, where $\mathbb{S} = \textrm{Res}_{\mathbb{C}/\mathbb{R}} \mathbb{G}_{m,\mathbb{R}}$ is the Deligne torus. Then $h$ factors through the Mumford-Tate group $\mathbf{G}_{S} \subset \GL(\mathbb{V}_{s})$ of the variation $\mathbb{V}$ (the stabilizer of global Hodge tensors of $\bigoplus_{a, b \geq 0} \mathbb{V}^{\otimes a} \otimes (\mathbb{V}^{*})^{\otimes b}$), and so the map $\textrm{Ad}\, h$ induces a Hodge structure on the Lie algebra $\mathfrak{g}_{S}$ of $\mathbf{G}_{S}$. By a theorem of Andr\'e-Deligne \cite{Andre1992}, the Lie algebra $\mathfrak{h}_{S}$ is a semisimple summand of the reductive Lie algebra $\mathfrak{g}_{S}$, hence the map $\textrm{Ad}\, h$ preserves $\mathfrak{h}_{S}$ as well and induces the desired Hodge structure. In this context we define:
\begin{defn}
If $\mathfrak{h}_{S,\mathbb{C}} = \bigoplus_{i} \mathfrak{h}^{i}$ is the associated Hodge decomposition, we call the \emph{level} of the variation $\mathbb{V}$ the largest $i$ for which $\mathfrak{h}^{i} \neq 0$.
\end{defn}
\noindent The second hypothesis of \autoref{babythm} then says that $\mathbb{V}$ has level at least three. 

\subsubsection{The Case of Tensors}
\label{tensorcasesec}

In order to give a more general result for loci defined by tensors, we begin by introducing some more background. One feature of the Zariski non-density result appearing in \cite{BKU} is that, because the methods are Hodge-theoretic and do not deal explicitly with families of algebraic cycles, the resulting theorems constrain loci which are merely \emph{conjectured} to correspond to families of algebraic cycles. To accomplish something similar we will use the theory of crystals, which we review in \S\ref{crysflatsec}. For any positive-characteristic subscheme $Z \subset S_{\overline{\kappa}}$, where $\kappa = \kappa(\mathfrak{p})$ is a positive-characteristic residue field corresponding to a point $\mathfrak{p} \in \Spec R$, we write $\mathcal{H}_{\textrm{cris},Z} = R^{k} f_{Z,*} \mathcal{O}_{X_{Z}/W}$ for the associated crystal over $Z$, where $f_{Z} : X_{Z} \to Z$ is the base-change, $W = W(\overline{\kappa})$ is the Witt vector ring, and $\mathcal{O}_{X_{Z}/W}$ is the crystalline structure sheaf.  We note that the value of the sheaf $\mathcal{H}_{\textrm{cris},Z}$ on $(Z, Z, 0)$ is exactly $\restr{\mathcal{H}}{Z}$, and as we review in \S\ref{crysflatsec} every global section of $\mathcal{H}_{\textrm{cris},Z}$ induces a flat section of $\restr{\mathcal{H}}{Z}$.

In addition to the Hodge filtration $F^{\bullet}$, for each $\mathfrak{p}$ and $Z$ as above the bundle $\restr{\mathcal{H}}{Z}$ carries a canonical ``conjugate'' filtration $F^{\bullet}_{c}$, which for fixed $\mathfrak{p}$ is compatible with restriction along subschemes of $S_{\overline{\kappa}}$. We define:

\begin{defn}
\label{ordinarydef}
Let $Z \subset S_{\overline{\kappa}}$ be a closed subscheme, with $\kappa = \kappa(\mathfrak{p})$ and $\mathfrak{p} \in \Spec R$. We say $Z$ is ordinary if either:
\begin{itemize}
\item[(i)] the field $\kappa(\mathfrak{p})$ is of characteristic zero; or
\item[(ii)] the filtrations $F^{\bullet}$ and $F^{\bullet}_{c}$ are opposed (see \autoref{oppdef}) over the generic fibre of $Z$.
\end{itemize}
If $Z$ is a point, we say it is an ordinary point. If $Z \subset S$ lies over $\mathfrak{p}$, we also say that $Z$ is ordinary (at $\mathfrak{p}$) if $Z_{\overline{\kappa(\mathfrak{p})}}$ is ordinary. 
\end{defn}

\noindent For the majority of families $f$ that one can consider, the ``ordinary'' case is the generic one, in the sense that the locus of ordinary points is dense in $S_{\overline{\kappa(\mathfrak{p})}}$ (c.f. \autoref{compinstcor} below). We will primarily consider the ordinary case in this paper; the methods can be applied more generally, but in the non-ordinary cases the analysis has to be done separately so we consider only the ordinary case for simplicity.

For a given $Z$ as above, we let $K$ be the fraction field of $W$, and consider the category $(Z/W)_{\textrm{cris},K}$ of \emph{isocrystals} obtained from the category $(Z/W)_{\textrm{cris}}$ of crystals on the crystalline site of $Z/W$ by leaving the objects unchanged and tensoring the Hom sets with $K$. The category  $(Z/W)_{\textrm{cris},K}$ is a neutral Tannakian category, so we may fix a fibre functor $F_{Z} : (Z/W)_{\textrm{cris},K} \to \textrm{Vect}_{K}$. (For now the precise choice of $F_{Z}$ is not important, but we will give a more explicit construction later.) We may use $F_{Z}$ to talk about algebraic groups defined as the stabilizers of global sections of $\mathcal{H}^{a,b}_{\textrm{cris},Z} = \mathcal{H}^{\otimes a}_{\textrm{cris},Z} \otimes (\mathcal{H}^{*}_{\textrm{cris},Z})^{\otimes b}$, where we view such global sections as morphisms from $\mathcal{O}_{Z/W}$, and then use $F_{Z}$ to associate these morphisms to elements of $F_{Z}(\mathcal{H}^{a,b}_{\textrm{cris}, Z}) = F_{Z}(\mathcal{H}_{\textrm{cris}, Z})^{a,b}$. For each $a,b$, we write $\mathcal{A}^{a,b}_{\textrm{cris}} = \mathcal{A}^{a,b}_{\textrm{cris}}(Z)$ for the collection of all global sections of $\mathcal{H}^{a,b}_{\textrm{cris},Z}$ which are induced by diagrams of the form 
\begin{equation}
\label{cyclediag2}
\begin{tikzcd}
Y \arrow[r, hook] \arrow[d, swap] & X^{\otimes n}_{\overline{\kappa}} \arrow[d, "f^{\otimes n}_{\overline{\kappa}}"] \\
Z \arrow[r, hook] & S_{\overline{\kappa}}
\end{tikzcd} ,
\end{equation}
with $Y$ smooth over $Z$. (See \S\ref{crysflatsec} for a discussion of the crystalline cycle class map.) Note that we may also use the fibre functor $F_{Z}$ to give a crystalline realization of $\mathbf{H}_{S}$ as the stabilizer of the crystalline incarnation of the polarization pairing.

In the case that $\kappa = \kappa(\mathfrak{p})$ is a field of characteristic zero and $Z \subset S_{\overline{\kappa}}$ is a scheme lying over $\mathfrak{p}$, we we write $\mathcal{A}^{a,b} = \mathcal{A}^{a,b}(Z)$ for the collection of all global sections of $\restr{\mathcal{H}^{a,b}}{Z}$ which are induced by diagrams as in (\ref{cyclediag2}), with $Y$ smooth over $Z$ (see \S\ref{flatchernsec} for a discussion of the algebraic de Rham cycle class map). 

\begin{thm}
\label{properlocusthm}
Suppose that $\varphi$ is quasi-finite, that $\mathbf{H}_{S} = \textrm{Aut}(\mathbb{V}_{s}, \mathcal{Q}_{s})$ for some $s \in S(\mathbb{C})$, and that $\mathbb{V}$ has level at least three. Consider the following three properties for subschemes $Z \subset S$:
\begin{itemize}
\item[(1)] $Z$ is positive dimensional over the prime $\mathfrak{p} \in \Spec R$.
\item[(2)] With $\kappa = \kappa(\mathfrak{p})$, that either
\begin{itemize}
\item[(i)] if $\textrm{char}\, \kappa = 0$, then the algebraic subgroup $\mathbf{H}_{\mathcal{A},s} \subset \textrm{Aut}(\mathcal{H}_{s})$ stabilizing the elements of $\mathcal{A} = \bigcup_{a,b} \mathcal{A}^{a,b}(Z_{\overline{\kappa}})$ is properly contained in $\mathbf{H}_{S,s}$ for some (hence any) $s \in Z(\overline{\kappa})$; 
\item[(ii)] if $\textrm{char}\, \kappa > 0$, then there exists a subset $\mathcal{S} \subset \mathcal{A}_{\textrm{cris}} = \bigcup_{a,b} \mathcal{A}^{a,b}_{\textrm{cris}}(Z_{\overline{\kappa}})$ such that the algebraic subgroup $\mathbf{H}_{\mathcal{S}}$ stabilizing the elements of $\mathcal{S}$ is semisimple and properly contained in $\mathbf{H}_{S}$.
\end{itemize}
\item[(3)] $Z$ is ordinary at $\mathfrak{p}$.
\end{itemize}
Then there exists a closed subscheme $E$, properly contained in $S$, which contains all locally closed subschemes $Z \subset S$ satisfying (1), (2) and (3).
\end{thm}

\begin{thm}
\label{mereconj}
More generally, suppose one replaces $\mathcal{A}$ and $\mathcal{A}_{\textrm{cris}}$ in \autoref{properlocusthm} with those sections merely \emph{conjectured} to come from families of algebraic cycles (see the proof of \autoref{masterthm} for an explanation). Then \autoref{properlocusthm} continues to hold.
\end{thm}

We comment briefly on the semisimplicity assumption in \autoref{properlocusthm}(2)(ii), which may appear unmotivated. A crucial difference between the ``vector'' and ``tensor'' cases is that we do not impose any bound on the integer $n$ that appears in (\ref{cyclediag2}), and so have to consider loci defined by algebraic cycles coming from self-products $X^{\otimes n}$ of arbitrary high order. This presents a problem for the Zariski closedness of the locus $E$, because it could be that while the loci corresponding to families of algebraic cycles associated to $f^{\otimes n}$ all lie in such a locus for each $n$, one could fail to get a Zariski closed locus after taking the union over all $n$. In characteristic zero one can appeal to the reductivity of Hodge-theoretic Mumford-Tate groups (and the fact that they lie in finitely many complex conjugacy classes) to constrain the tensors involved, but no such constraints are available in positive characteristic and the semisimplicity requirement in \autoref{properlocusthm2}(2)(ii) is one way of imposing them directly.

We also note that one can give a version of \autoref{properlocusthm} without referring to crystalline cohomology at all (and indeed, without the semisimplicity assumption) at the cost of requiring an explicit bound on $n$; this is a consequence of our main technical result \autoref{bigexpthm}. Note that the Hodge symmetry assumption appearing in \autoref{bigexpthm} is implied by semisimplicity away from a closed locus in $\Spec R$ as shown in \autoref{semsimpimpsymlem}, so in fact the Hodge symmetry requirement in \autoref{bigexpthm} is the weaker condition. 

\vspace{0.5em}

At this point we note that the quasi-finiteness assumption on $\varphi$ which appears in \autoref{babythm} and \autoref{properlocusthm} is only there so one can equate the locus $Z$ being positive dimensional over $R$ with the condition that the Hodge flag $F^{\bullet}$ on $\mathcal{H}$  varies non-trivially over $Z$. If one drops this quasi-finiteness assumption, one can obtain a more general theorem, subject to the introduction of the technical term ``$r$-limp'' which we will introduce in \S\ref{infpermapsec}. 

\begin{thm}
\label{properlocusthm2}
More generally consider the same situation as \autoref{properlocusthm} (or \ref{mereconj}) without the quasi-finiteness assumption on $\varphi$, and where condition (1) is replaced by the property that the data $\restr{(\mathcal{H}, F^{\bullet}, \nabla)}{Z}$ admits a non-constant $r$-limp of some order $r$. Then the same conclusion holds. 
\end{thm}

\begin{rem}
An $r$-limp is a formal object which describes the infinitesimal $r$'th order variation of the Hodge flag at some point. Thus, $r$-limp condition on $Z$ can be read informally as ``for some $r$, the $r$th-order variation of the Hodge flag is non-trivial at some point in $Z$.''
\end{rem}


\vspace{0.5em}

We note that the inability in \autoref{properlocusthm} and \autoref{properlocusthm2} to constrain the non-ordinary locus on the integral level does not preclude a non-Zariski density result individually at each prime, as the following corollary shows:

\begin{cor}
\label{compinstcor}
Let $f : X \to S$ be the universal family over $\mathbb{Z}$ of smooth complete intersections of dimension $j$ in $\mathbb{P}^{j+m}$ of multidegree $(a_{1}, \hdots, a_{m})$, and suppose that the associated variation of Hodge structures $\mathbb{V}$ obtained from primitive cohomology has level at least three. Let $N$ be a positive integer. Then there exists an integer $e$ such that when $p > e$, the locus in $S_{\mathbb{F}_{p}}$ which is the union of all (possibly non-ordinary) $Z$ such that
\begin{itemize}
\item[(1)] the data $\restr{(\mathcal{H}, F^{\bullet}, \nabla)}{Z}$ admits a non-constant $r$-limp of some order $r$; and 
\item[(2)] there exists families of algebraic cycles over $Z$ satisfying (2)(ii) of \autoref{properlocusthm};
\end{itemize} 
is not Zariski dense in $S_{\mathbb{F}_{p}}$. When $m = 1$ the condition (1) can be replaced by
\begin{itemize}
\item[(1')] the locus $Z$ does not lie a $\GL_{M,\mathbb{F}_{p}}$-orbit, where $M = {{j+a_{1}+1} \choose {a_{1}}} - 1$ and we consider the natural action of $\GL_{M,\mathbb{F}_{p}}$ on $S$ given by acting on the standard coordinates of $\mathbb{P}^{j+1}$.
\end{itemize}
\end{cor}

\noindent The key point is that a result of Illusie \cite{Illusie2007} shows that in the setting of \autoref{compinstcor} the ordinary locus is open in $S_{\mathbb{F}_{p}}$ for each $p$. 

\subsection{Infinitesimal Period Maps}
\label{infpermapsec}

In what follows we fix a scheme $S$ smooth over a ring $R$, and let $(\mathcal{H}, F^{\bullet}, \nabla)$ be a triple consisting of a vector bundle $\mathcal{H}$, a decreasing filtration $F^{\bullet}$, and a flat connection $\nabla : \mathcal{H} \to \Omega^{1}_{S/R} \otimes \mathcal{H}$, all defined over $R$. We have in mind the case where $\mathcal{H} = R^{n} f_{*} \Omega^{\bullet}_{X/S}$ is the relative algebraic de Rham cohomology of a smooth $R$-algebraic morphism $f : X \to S$, $F^{\bullet}$ is the Hodge filtration, and $\nabla$ is the Gauss-Manin connection as described in \cite{katz1968}. (As we mentioned above, the sheaves $\mathcal{H}$ and $F^{\bullet}$ are not always locally free in such a setting, but become so after replacing $\Spec R$ with an affine open subscheme.) We suppose that $\mathcal{H}$ has rank $m$, and write $h^{i} = \dim F^{i} / F^{i+1}$ for all $i$. 

\begin{defn}
\label{filcompframdef}
Suppose that $U \subset S$ is a Zariski open subscheme. By a \emph{filtration-compatible frame} over $U$, we mean a basis of sections $v^{1}, \hdots, v^{m}$ of $\mathcal{H}(U)$ such that for each integer $k$ there exists $j_{k}$ such that $F^{k}(U)$ is spanned by $v^{1}, \hdots, v^{j_{k}}$. 
\end{defn}

\begin{defn}
We will say that a filtration $F^{\bullet}$ on the $A$-module $A^{m}$ is graded locally free (GLF) if $F^{i} / F^{i+1}$ is locally free as an $A$-module for all $i$.
\end{defn}

\begin{notn}
We denote by $\ch{L}$ the $\mathbb{Z}$-scheme whose $A$-points, for $A$ a ring, are given by
\[ \ch{L}(A) = \left\{ \textrm{GLF Filtrations }F^{\bullet}\textrm{ of }A^{m}\textrm{ with }\dim F^{i} / F^{i+1} = h^{i} \right\} , \]
and by $\ch{L}_{R}$ its base-change along the canonical map $\Spec R \to \Spec \mathbb{Z}$. We have a natural map $q : \GL_{m} \to \ch{L}$ which associates to a point $M \in \GL_{m}(A)$ the filtration of $A^{m}$ for which $F^{k} A^{m}$ is the span of the first $j_{k}$ columns of $M$. We regard $\ch{L}$ as a homogeneous space for $\GL_{m}$ via the natural functorial action of $\GL_{m}(A)$ on filtrations of $A^{m}$. 
\end{notn}

\vspace{0.5em}

Suppose that $U \subset S$ is a Zariski open affine subscheme over $R$ admitting an \'etale map $\delta : U \to \mathbb{A}^{n}_{R}$ over $R$, and let $z_{1}, \hdots, z_{n} \in \mathcal{O}_{S}(U)$ be the induced sections. Let $I \subset \mathcal{O}_{S}(U)$ be the ideal generated by $z_{1}, \hdots, z_{n}$, and write $\mathcal{O}^{j}_{S, \delta} = \mathcal{O}_{S}(U)/I^{j+1}$ and $\mathcal{O}_{S, \delta} = \mathcal{O}^{\infty}_{S,\delta} = \varprojlim_{j} \mathcal{O}^{j}_{S,\delta}$. The data $(\mathcal{H}, F^{\bullet}, \nabla)$ naturally induces a filtered module with connection $(\mathcal{H}^{j}_{\delta}, F^{j, \bullet}_{\delta}, \nabla^{j}_{\delta})$ over the ring $\mathcal{O}^{j}_{S,\delta}$ for all $j \in \mathbb{N} \cup \{ \infty \}$. We say a frame $b^{1}, \hdots, b^{m}$ for $\mathcal{H}^{j}_{\delta}$ is flat if $\nabla b^{i} = 0$ for all $1 \leq i \leq m$.

\begin{defn}
\label{locpermapdef}
With the above setup and some $j \in \mathbb{N} \cup \{ \infty \}$, suppose that we have a filtration compatible frame $v^{1}, \hdots, v^{m}$ for $\mathcal{H}^{j}_{\delta}$, and a flat frame $b^{1}, \hdots, b^{m}$ for $\mathcal{H}^{j}_{\delta}$. Then if $M \in \GL_{m}(\mathcal{O}^{j}_{S,\delta})$ is the change-of-basis matrix from $v^{1}, \hdots, v^{m}$ to $b^{1}, \hdots, b^{m}$ then we call the composition $\psi = q \circ M$ an \emph{\'etale local infinitesimal period map} of order $j$. For ease of terminology, we will also say $\psi$ is a ``$j$-elimp''. 
\end{defn}

\begin{defn}
\label{locpermapdef2}
In the setting of \autoref{locpermapdef}, if the natural map $\delta^{-1}(0) \to \Spec R$ is an isomorphism, we call $\psi$ a \emph{local infinitesimal period map} over order $j$, or ``$j$-limp''.
\end{defn}

In general, $j$-elimps and $j$-limps are formal objects that describe a variation of the filtration on $\mathcal{H}$ in the order $j$ formal $R$-neighbourhood of the fibre $\delta^{-1}(0)$. We note that for an elimp these neighbourhoods are in general disconnected, with components corresponding to the components of $\delta^{-1}(0)$; this justifies the use of the term ``\'etale local'' as opposed to merely ``local'' which we used in \cite{periodimages} and \cite{urbanik2021sets}. Actually, for applications in this paper working exclusively with $j$-limps instead of $j$-elimps is enough, but we find the generality of $j$-elimps more natural as for a general $R$ one cannot guarantee a choice of $\delta$ for which $\delta^{-1}(0)$ consists of a single $R$-point. In situations where the ring $R$ admits a natural topology (e.g., $R = \mathbb{C}$ or $R = \mathbb{Q}_{p}$), $\infty$-limps can be identified with germs of analytic functions provided the formal solutions of an appropriate differential system determined by $\delta$ and $\nabla$ converge. In this situation, these analytic functions we will call ``local period maps'', which was the terminology used in \cite{periodimages} and \cite{urbanik2021sets}. (In our previous work we did not concern ourselves with purely formal period maps, but for our purposes in this paper situations where $R$ is a finite field will also be important, hence the need for the more general formal and infinitesimal notions.)

\begin{defn}
In the situation of \autoref{locpermapdef} we say that $\psi$ is \emph{defined} at a point $s \in S(R)$ if $\delta(s) = 0$. 
\end{defn}

In \cite{periodimages} and \cite{urbanik2021sets} we developed machinery for studying infinitesimal images of local period maps using algebro-geometric methods; the idea, roughly, is to use the differential equations that define the maps $\psi$ to ``evaluate'' $\psi$ on small infinitesimal disks, and study algebraic constraints on the images of these disks. To model such infinitesimal disks we use the notion of (higher-dimensional) \emph{jet spaces}, which are defined as follows:

\vspace{0.5em}

\begin{notn}
For any ring $R$, we define $A^{d}_{r, R} = R[t_{1}, \hdots, t_{d}]/(t_{1}, \hdots, t_{d})^{r+1}$, and write $\mathbb{D}^{d}_{r, R}$ for $\Spec A^{d}_{r, R}$. 
\end{notn}

\begin{defn}
Suppose that $S$ is an $R$-scheme, the \emph{jet space} $J^{d}_{r} S$ is defined to be the $R$-scheme representing the functor $\Sch_{R} \to \textrm{Set}$ given by
\[ T \mapsto \Hom_{R}(T \times_{R} \mathbb{D}^{d}_{r, R}, S), \hspace{1.5em} [T \to T'] \mapsto [\Hom_{R}(T' \times_{R} \mathbb{D}^{d}_{r, R}, S) \to \Hom_{R}(T \times_{R} \mathbb{D}^{d}_{r, R}, S)] , \]
where the natural map $\Hom_{R}(T' \times_{R} \mathbb{D}^{d}_{r, R}, S) \to \Hom_{R}(T \times_{R} \mathbb{D}^{d}_{r, R}, S)$ obtained by pulling back along $T \times_{R} \mathbb{D}^{d}_{r, R} \to T' \times_{R} \mathbb{D}^{d}_{r, R}$.
\end{defn}

\vspace{0.5em}

The representability of the functor defining $J^{d}_{r} S$ reduces as in \cite[\S2.1]{periodimages} to the representability of Weil restrictions, and hence holds for instance when $S$ is quasi-projective over $R$. Note that a map $g : S \to S'$ of $R$-schemes induces a natural map $J^{d}_{r} S \to J^{d}_{r} S'$ by post-composition.

The construction central to our main results, which generalizes \cite[Theorem 1.11]{periodimages} and \cite[Theorem 3.3]{urbanik2021sets}, is the following:

\begin{thm}
\label{bigjetconstr}
Let $S$ be a quasi-projective $R$-scheme, with $R$ an integral domain over $\mathbb{Z}[1/r!]$, and let $(\mathcal{H}, F^{\bullet}, \nabla)$ be a filtered vector bundle and flat connection on $S$ defined over $R$. Then there exists a canonical map of $R$-algebraic stacks
\[ \eta^{d}_{r} : J^{d}_{r} S \to \GL_{m,R} \backslash J^{d}_{r} \ch{L}_{R} , \]
with the following properties:
\begin{itemize}
\item[(i)] the formation of $\eta^{d}_{r}$ is compatible with base-change along ring maps $R \to R'$ of integral domains over $\mathbb{Z}[1/r!]$, in the sense that if $\eta'^{d}_{r}$ is the map associated to the triple $(\mathcal{H}', F'^{\bullet}, \nabla')$ obtained by pullback to $S' = S_{R'}$, the natural diagram
\begin{center}
\begin{tikzcd}
J^{d}_{r} S' \arrow[r, "\eta'^{d}_{r}"] \arrow[d] & \GL_{m,R'} \backslash J^{d}_{r} \ch{L}_{R'} \arrow[d] \\
J^{d}_{r} S \arrow[r, "\eta^{d}_{r}"] & \GL_{m,R} \backslash J^{d}_{r} \ch{L}_{R}
\end{tikzcd} .
\end{center}
commutes;
\item[(ii)] if $g : S' \to S$ is a map of smooth quasi-projective $R$-schemes, and $\eta'^{d}_{r}$ is the map associated to the triple $(\mathcal{H}', F'^{\bullet}, \nabla')$ obtained by pullback along $g$, then $\eta'^{d}_{r} = \eta^{d}_{r} \circ J^{d}_{r} g$;
\item[(iii)] given a map $v : \mathbb{D}^{d}_{r, R} \to \mathbb{D}^{d'}_{r', R}$ over $R$ with $r' \leq r$, the diagram 
\begin{center}
\begin{tikzcd}
J^{d'}_{r'} S \arrow[r, "\eta^{d'}_{r'}"] \arrow[d, "(-) \circ v"] & \GL_{m,R} \backslash J^{d}_{r} \ch{L}_{R} \arrow[d, "(-) \circ v"] \\
J^{d}_{r} S \arrow[r, "\eta^{d}_{r}"] & \GL_{m,R} \backslash J^{d}_{r} \ch{L}_{R}
\end{tikzcd}
\end{center}
commutes (see \S\ref{keypropsec} for details);
\item[(iv)] for any $r$-elimp $\psi$ associated to the \'etale neighbourhood $\delta : U \to \mathbb{A}^{n}_{R}$, and any jet $j \in (J^{d}_{r} U)(R)$ such $\psi$ is defined at the basepoint of $j$, we have $\psi \circ j = \eta^{d}_{r}(j)$ as points in $(\GL_{m} \backslash J^{d}_{r} \ch{L})(R)$. 
\end{itemize}
\end{thm}

\vspace{0.5em}

\begin{rem}
The condition that $R$ be an integral domain over $\mathbb{Z}[1/r!]$ ultimately comes from the need to consider derivatives of the maps $\psi$ up to order $r$, which may be poorly behaved in the presence of zero divisors and if the characteristic of $R$ is not greater than $r$. 
\end{rem}

\subsection{Acknowledgements}

The author thanks Jacob Tsimerman for numerous helpful conversations related to the ideas in this manuscript. The author also thanks François Charles for his careful reading and feedback.

\section{Jets of Period Maps}

\subsection{Jet Maps and Derivatives}
\label{jetmapsanddersec}

In what follows we fix a map $g : \mathbb{A}^{k}_{R} \to \mathbb{A}^{\ell}_{R}$ of affine spaces over $R$, and give an explicit description of the map $J^{d}_{r} g : J^{d}_{r} \mathbb{A}^{k}_{R} \to J^{d}_{r} \mathbb{A}^{k}_{R}$. We work in the setting where $R$ is an integral domain over $\mathbb{Z}[1/r!]$. In what follows we denote by $\mathcal{Q}^{d}_{r}$ the set of partitions of integers $0$ through $r$ which have $d$ terms, and by $\varnothing$ the empty partition. An element $\sigma \in (J^{d}_{r} \mathbb{A}^{k})(R)$ is naturally identified with $k$ formal sums $\sum_{p \in \mathcal{Q}^{d}_{r}} a_{p,i} \overline{t}^{p}$, where each $a_{p,i} \in R$, we have $1 \leq i \leq k$, and we use multi-index notation to exponentiate the tuple $\overline{t} = (t_{1}, \hdots, t_{d})$. Viewing the map $g$ as a tuple $(g_{1}, \hdots, g_{\ell})$ we compute $g \circ \sigma$ by computing, for each $1 \leq j \leq \ell$, the compositions
\begin{equation} 
\label{jeteq}
(g \circ \sigma)_{j} = g_{j}\left( \sum_{p \in \mathcal{Q}^{d}_{r}} a_{p,1} \overline{t}^{p}, \hdots, \sum_{p \in \mathcal{Q}^{d}_{r}} a_{p,\ell} \overline{t}^{p} \right) . 
\end{equation}
Let us write $(g \circ \sigma)_{j} = \sum_{p \in \mathcal{Q}^{d}_{r}} b_{p,j} \overline{t}^{p}$, and view both sides of the equality  (\ref{jeteq}) as formal expressions. Taking $p$ to be the partition $i_{1} + \cdots + i_{q}$, we may repeatedly differentiate both sides of (\ref{jeteq}) using the multivariate chain rule and evaluate at $\overline{t} = 0$ to obtain an equality
\begin{equation}
\label{curryeq}
b_{p,j} = \frac{1}{i_{1} ! \cdots i_{q} !} \left( \textrm{polynomial in } \underbrace{\left\{ \begin{array}{cc} \partial_{p} g_{j} \left( a_{\varnothing, 1}, \hdots, a_{\varnothing, \ell} \right), &p \in \mathcal{Q}^{d}_{r} ; \\ a_{p,i}, &p \in \mathcal{Q}^{d}_{r}, 1 \leq i \leq k \end{array} \right\}}_{\mathcal{S}_{j}} \right) ,
\end{equation}
where we have used the fact that $R$ is a $\mathbb{Z}[1/r!]$-algebra to invert $i_{1}! \cdots i_{q}!$. In the above, $\mathcal{S}_{j}$ denotes the finite set of formal symbols appearing on the left-hand side of the bracketed expression with indicies in the ranges given on the right-hand side of the bracketed expression. 

Now let us denote by $B$ the polynomial algebra $R[\mathcal{S}]$, where $\mathcal{S} = \mathcal{S}_{1} \cup \cdots \cup \mathcal{S}_{\ell}$. Then for each $p$ and $j$, the polynomial appearing on the right-hand side of (\ref{curryeq}), which we denote $h_{p,j}$, is an element of $R[\mathcal{S}]$. Then the above discussion proves the following:

\begin{prop}
\label{curryprop}
Let $e : \Spec R[\mathcal{S}] \to J^{d}_{r} \mathbb{A}^{\ell}_{R}$ be the map given by $b_{p,j} \mapsto h_{p,j}$, and let $\varphi_{g} : J^{d}_{r} \mathbb{A}^{k}_{R} \to \Spec R[\mathcal{S}]$ be the map given by instantiating each formal symbol in $\mathcal{S}$ with the associated polynomial in the coordinates of $J^{d}_{r} \mathbb{A}^{k}_{R}$ determined by $g$. Then the diagram
\begin{center}
\begin{tikzcd}
J^{d}_{r} \mathbb{A}^{k}_{R} \arrow[dr, "J^{d}_{r} g"] \arrow[d, "\varphi_{g}"] &  \\
\Spec R[\mathcal{S}] \arrow[r, "e"] & J^{d}_{r} \mathbb{A}^{\ell}_{R}
\end{tikzcd} 
\end{center}
commutes. 
\end{prop}

In fact, the statement of \autoref{curryprop} easily generalizes to include the situation where we consider maps between formal $R$-algebraic neighbourhoods of $\mathbb{A}^{k}_{R}$ and $\mathbb{A}^{\ell}_{R}$. To explain what we mean, let us fix points $t \in \mathbb{A}^{k}(R)$ and $u \in \mathbb{A}^{\ell}(R)$. Using the product structure of affine space, we have naturally associated tuples $(t_{1}, \hdots, t_{k})$ and $(u_{1}, \hdots, u_{\ell})$, where the entries of the tuples are points of $\mathbb{A}^{1}(R)$. If $x_{1}, \hdots, x_{k}$ are the natural coordinates of $\mathbb{A}^{k}$ and similarly $y_{1}, \hdots, y_{\ell}$ are the natural coordinates of $\mathbb{A}^{\ell}$, formal completion at the ideals $I = (x_{1} - t_{1}, \hdots, x_{k} - t_{k})$ and $J = (y_{1} - u_{1}, \hdots, y_{\ell} - u_{\ell})$ gives rings of formal power series $\mathcal{O}_{\mathbb{A}^{k}_{R}, t}$ and $\mathcal{O}_{\mathbb{A}^{\ell}_{R}, u}$.

Suppose now that we have a map $\widehat{g} : \Spec \mathcal{O}_{\mathbb{A}^{k}_{R}, t} \to \Spec \mathcal{O}_{\mathbb{A}^{\ell}_{R}, u}$. Such a map induces a map $J^{d}_{r} \widehat{g} : J^{d}_{r} \mathbb{A}^{k}_{t} \to J^{d}_{r} \mathbb{A}^{\ell}_{u}$ between the fibres of the jet spaces above $u$ and $t$, and this induced map depends only on the map $\widehat{g}^{r} : \Spec \mathcal{O}_{\mathbb{A}^{k}_{R}, t}/I^{r+1} \to \Spec \mathcal{O}_{\mathbb{A}^{\ell}_{R}, u}/J^{r+1}$ induced by $\widehat{g}$. The map $\widehat{g}^{r}$ is also induced by a polynomial map $g : \mathbb{A}^{k}_{R} \to \mathbb{A}^{\ell}_{r}$, which may be obtained by from the map $\widehat{g}$ by truncating the power series defining $\widehat{g}$ after order $r$. Applying \autoref{curryprop} and the functorial properties of the jet space construction, we obtain the following:

\begin{prop}
\label{curryprop2}
Let $e$ be as in \autoref{curryprop}, and denote by $\varphi_{\widehat{g}}$ the restriction of $\varphi_{g}$ to $J^{d}_{r} \mathbb{A}^{k}_{t}$. Then the diagram
\begin{center}
\begin{tikzcd}
J^{d}_{r} \mathbb{A}^{k}_{t} \arrow[rd, "J^{d}_{r} \widehat{g}"] \arrow[d, "\varphi_{\widehat{g}}"] & \\
\Spec R[\mathcal{S}] \arrow[r, "e"] & J^{d}_{r} \mathbb{A}^{\ell}_{u}
\end{tikzcd} 
\end{center}
commutes. 
\end{prop}

\subsection{The Main Construction}
\label{mainjetconstrsec}

In this section we prove \autoref{bigjetconstr}, which means constructing a $\GL_{m,R}$-torsor $\gamma : \mathcal{P}^{d}_{r} \to J^{d}_{r} S$ as well as a $\GL_{m,R}$-invariant algebraic map $\alpha : \mathcal{P}^{d}_{r} \to J^{d}_{r} \ch{L}_{R}$ such that the pair $(\gamma, \alpha)$ defines the map $\eta^{d}_{r}$ using the usual description of maps to a quotient stack. We will carry out the construction locally on $S$, with the construction on a general $S$ being obtained by gluing. The process is analogous to the arguments that appear in \cite[\S3]{periodimages} and \cite[\S3]{urbanik2021sets}, which the reader can consult for additional details. Note, however, that unlike in \cite[\S3]{periodimages} and \cite[\S3]{urbanik2021sets}, we do not use the existence and unicity of formal flat frames for the Hodge bundle, as this cannot be assumed when $R$ is not a field of characteristic zero. Instead of this our construction will use an infinitesimal version of the necessary property for $r$-limps which is shown in \autoref{existslem} below using the fact that $R$ is a $\mathbb{Z}[1/r!]$-algebra.

\subsubsection{Local Construction}
\label{locconstrsec}

With this in mind, after possibly shrinking the base $S$ we assume that we have $S = \Spec A$ for an $R$-algebra $A$, a filtration-compatible trivialization $v^{1}, \hdots, v^{m}$ of the filtered vector bundle $(\mathcal{H}, F^{\bullet})$, and that $\Omega^{1}_{A/R}$ is free with basis $dz_{1}, \hdots, dz_{n}$. The relative tangent sheaf $\mathcal{T}_{S/R}$, which is dual to $\Omega^{1}_{A/R}$, is then trivialized by the differential operators $\partial_{1}, \hdots, \partial_{n}$ obtained from duality. We write $\nabla v^{i} = \sum_{j = 1}^{m} c_{ij} \otimes v^{j}$ for sections $c_{ij} \in \Omega^{1}_{A/R}$, which we further expand as $c_{ij} = \sum_{\ell = 1}^{n} c_{ij,\ell} dz_{\ell}$ for functions $c_{ij,\ell} \in A$. 

Let $f_{tu}$ be a collection of $m^2$ formal symbols with $1 \leq t, u \leq m$. We now define a family of polynomials $\xi_{q,jk} \in A[f_{tu}, 1 \leq t, u \leq m]$, indexed by integers $1 \leq j, k \leq m$ and elements $q \in \mathcal{M}([n], 1) \cup \mathcal{M}([n], 2) \cup \cdots$, where $\mathcal{M}([n], a)$ denotes the set of size $a$ multiset subsets of $[n] = \{ 1, \hdots, n \}$; that is to say, the element $q$ is a finite unordered sequence of integers in the set $[n]$. The polynomials $\xi_{q,jk}$ will \emph{later} (in subsequent sections) correspond to the derivatives with respect to the differential operators $\partial_{1}, \hdots, \partial_{n}$ of a matrix-valued map $f = [f_{tu}]$ for which $f^{-1}$ is the period matrix between the frame $v^{1}, \hdots, v^{m}$ and an (infinitesimal) flat frame. However, we emphasize that we do not assume the existence of any flat frames here and that these polynomials are purely formal objects for the time being. We define
\begin{equation}
\label{derivformula}
\xi_{\{\ell\}, jk} = -\sum_{i = 1}^{m} f_{ik} c_{ij,\ell} .
\end{equation}
and for a general $q = \{ \ell_{1}, \hdots, \ell_{o} \}$ with $\ell = \ell_{o}$, we define $\xi_{q,jk}$ by applying the operator $\partial_{\ell_{1}} \cdots \partial_{\ell_{o-1}}$ to (\ref{derivformula}) and using (\ref{derivformula}) to express the result as an element of $A[f_{tu}, 1 \leq t, u \leq m]$. 

\begin{lem}
\label{xiwelldef}
The polynomials $\xi_{q, jk}$ are well-defined (independent of the choice of ordering of $\ell_{1}, \hdots, \ell_{o}$). 
\end{lem}

\begin{proof}
Because the operators $\partial_{1}, \hdots, \partial_{n}$ all pairwise commute, it suffices to check that that 
\begin{equation}
\label{derivscomm}
\partial_{\ell_{1}} \xi_{\{ \ell_{2} \}, jk} = \partial_{\ell_{2}} \xi_{\{ \ell_{1} \}, jk}
\end{equation}
for each choice of $\ell_{1}$ and $\ell_{2}$. Indeed, suppose we verify this condition, and then use the equality $\partial_{\ell} f_{jk} = \xi_{\{\ell\}, jk}$ for $1 \leq j, k \leq m$ to \emph{define}, via the Leibniz rule, an extension of the operator $\partial_{\ell}$ to the polynomial ring $A[f_{tu}, 1 \leq t, u \leq m]$. Then the equality (\ref{derivscomm}) will simply say that the extended operators commute pairwise, which implies the commutativity of all mixed partials inductively.

The condition (\ref{derivscomm}) comes from the integrability of the connection $\nabla$. To see this, simply compute that
\begin{align*}
\partial_{\ell_{2}} \xi_{\{ \ell_{1} \}, jk} - \partial_{\ell_{1}} \xi_{\{ \ell_{2} \}, jk} &= \sum_{i = 1}^{m} \left( (\partial_{\ell_{1}} f_{ik} c_{ij,\ell_{2}} - \partial_{\ell_{2}} f_{ik} c_{ij, \ell_{1}}) + (f_{ik} \partial_{\ell_{1}}c_{ij,\ell_{2}} - f_{ik} \partial_{\ell_{2}}c_{ij,\ell_{1}} ) \right) \\
&= \sum_{i = 1}^{m} \sum_{i' = 1}^{m}  \left( f_{i'k} c_{i'i, \ell_{2}} c_{ij, \ell_{1}} -  f_{i'k} c_{i'i,\ell_{1}} c_{ij,\ell_{2}}\right) + \sum_{i' = 1}^{m} (f_{i'k} \partial_{\ell_{1}}c_{i'j,\ell_{2}} - f_{i'k} \partial_{\ell_{2}}c_{i'j,\ell_{1}} ) \\
&= \sum_{i' = 1}^{m} f_{i'k} \left( \sum_{i = 1}^{m} (c_{i'i,\ell_{2}} c_{ij, \ell_{1}} - c_{i'i, \ell_{1}} c_{ij, \ell_{2}}) + \partial_{\ell_{1}} c_{i'j, \ell_{2}} - \partial_{\ell_{2}} c_{i'j,\ell_{1}} \right)
\end{align*}
The inner term for each value of $i'$ is simply the coordinate form of the curvature tensor, hence from flatness each term in the sum vanishes, hence the result. 
\end{proof}

We are now ready to construct the pair $(\gamma, \alpha)$. We define $\mathcal{P}^{d}_{r} = J^{d}_{r} S \times \GL_{m,R}$, and take for $\gamma$ the natural projection.  We let $\delta : S \to \mathbb{A}^{n}_{R}$ be the \'etale map over $R$ induced by the functions $z_{1}, \hdots, z_{n}$. The map $\alpha : \mathcal{P}^{d}_{r} \to J^{d}_{r} \ch{L}_{R}$ will be defined to be $(J^{d}_{r} (q \circ \iota))_{R} \circ e \circ \zeta$, where:

\begin{itemize}
\item[(i)] The map $J^{d}_{r} (q \circ \iota)$ is the map $J^{d}_{r} \GL_{m} \to J^{d}_{r} \ch{L}$ induced by the map $q \circ \iota$, where $q$ is as in \S\ref{infpermapsec} and $\iota : \GL_{m} \to \GL_{m}$ is the inversion.
\item[(ii)] The map $e$ is the map of \autoref{curryprop} and \autoref{curryprop2} above where we take $\mathbb{A}^{n}_{R}$ to be the ``source'' affine space and the space $\mathbb{M}_{R}$ of all $m \times m$ matrices over $R$ to be the ``target'' affine space. In particular, the map $g$ of \autoref{curryprop} is a map $\mathbb{A}^{n}_{R} \to \mathbb{M}_{R}$, and the functions $h_{p,j}$ (in the notation of \S\ref{jetmapsanddersec}) we regard as being indexed $h_{q,jk}$ where the indices $1 \leq j, k \leq m$ correspond to the entries $[M_{tu}]$ of the matrices in $\mathbb{M}_{R}$.
\item[(iii)] The map $\zeta : J^{d}_{r} S \times_{R} \GL_{m,R} \to \Spec R[\mathcal{S}]$ is the map given by
\begin{equation}
\label{zetadef}
 (\sigma, [M_{tu}]) \mapsto (\delta \circ \sigma, \xi_{q,jk}(\pi(\sigma), [M_{tu}]))), 
\end{equation}
on points, where the entries $\xi_{q,jk}(\pi(\sigma), [M_{tu}])$ are the values of the functions $h_{q,jk}$, $\pi : J^{d}_{r} S \to S$ is the natural projection, and we range over all choices of indices. This map on points can be interpreted in the sense of functors of points, and hence this gives a complete definition. In particular, we obtain an algebraic map over $R$. 
\end{itemize}

\noindent We note that, because the matrices $[M_{tu}]$ are assuemd to be points in $\GL_{m,R}$, the composition $e \circ \zeta$ has image in $J^{d}_{r} \GL_{m,R} \subset J^{d}_{r} \mathbb{M}_{R}$, so the definition makes sense.

\subsubsection{Relation to Period Maps}

\label{relpersec}

In preparation for extending the construction globally, we first relate it to the notion of \'etale infinitesimal local period maps (elimps) discussed in the introduction. All elimps in this section will be relative to the frame $v^{1}, \hdots, v^{m}$ fixed in \S\ref{locconstrsec}, which means that if $\psi = q \circ M$ is an elimp, then $M$ gives a change-of-basis matrix between a restriction to an infinitesimal neighbourhood of $v^{1}, \hdots, v^{m}$ and a flat frame.

\begin{lem}
\label{fagreescalc}
Suppose that $\psi = q \circ M$ is an $r$-elimp associated to $\delta' : U \to \mathbb{A}^{n}_{R}$, where $\delta'$ is the restriction of $\tau \circ \delta$ to an open $R$-subscheme $U \subset S$, with $\tau : \mathbb{A}^{n}_{R} \to \mathbb{A}^{n}_{R}$ a translation by an $R$-point, and $\delta$ the \'etale map fixed in \S\ref{locconstrsec}. Let $s : \Spec R \to S$ be an $R$-point above $0$. Then the components $f_{jk}$ of the matrix-valued infinitesimal map $f = M^{-1}$ have partial derivatives with respect to the operators $\partial_{1}, \hdots, \partial_{n}$ associated to $z_{1}, \hdots, z_{n}$ given at the $R$-point $s$ by $\xi_{q,jk}(s, f(s))$, where $q$ is a partition of an integer between $0$ and $r$ as in \S\ref{locconstrsec}.
\end{lem}

\begin{proof}
The flat frame $b^{1}, \hdots, b^{m}$ associated to $v^{1}, \hdots, v^{m}$ by $M$ satisfies $b^{k} = \sum_{i = 1}^{m} f_{ik} v^{i}$, and hence letting $\nabla v^{i} = \sum_{j = 1}^{m} c_{ij} \otimes v^{j}$ as in \S\ref{locconstrsec} we find that
\begin{align*}
\nabla b^{k} &= \nabla \left( \sum_{i = 1}^{m} f_{ik} v^{i} \right) \\
&= \sum_{j = 1}^{m} df_{jk} \otimes v^{j} + \sum_{i=1}^{m} f_{ik} \left( \sum_{j = 1}^{m} c_{ij} \otimes v^{j} \right) \\
&= \sum_{j = 1}^{m} \left( df_{jk} + \sum_{i = 1}^{m} f_{ik} c_{ij} \right) \otimes v^{j} , 
\end{align*}
and hence $\nabla b^{k} = 0$ implies that $df_{jk} = - \sum_{i = 1}^{m} f_{ik} c_{ij}$. This is in agreement with the differential system which defined the polynomials $\xi_{q,jk}$ in (\ref{derivformula}) as soon as one expands $c_{ij}$ in terms of the basis $dz_{1}, \hdots, dz_{n}$ for $\Omega^{1}_{S/R}$; we note that replacing $\delta$ with $\delta'$ has no effect on the basis $dz_{1}, \hdots, dz_{n}$, and hence the induced operators $\partial_{1}, \hdots, \partial_{n}$.  
\end{proof}

\begin{lem}
\label{evallem}
Suppose that $\psi = q \circ M$ is an $r$-elimp associated to the \'etale neighbourhood $\delta' : U \to \mathbb{A}^{n}_{R}$ as in \autoref{fagreescalc}, and that $s : \Spec R \to \delta^{-1}(0) \subset U$ is section of the structure map. Then if $M_{0}$ is the value of $M$ at $s$, we have $\alpha(j, M_{0}^{-1}) = \psi \circ j$ for any $j \in (J^{d}_{r} S)(R)$ lying above $s$.
\end{lem}

\begin{proof}
The fibre $\delta^{-1}(0)$ is a closed subscheme of $U$, quasi-finite over $R$, of which is the union of $s$ and a closed $R$-subscheme $c$. We may replace $U$ with $U \setminus c$ so that the map $\psi$ is an $r$-limp. To verify that $\alpha(j, M_{0}^{-1}) = \psi \circ j$, it suffices, using the construction $\alpha = (J^{d}_{r} q \circ \iota)_{R} \circ e \circ \zeta$, to verify that $(e \circ \zeta)(j, M_{0}^{-1}) = M^{-1} \circ j$. We will write $f = M^{-1}$ in what follows. Let $x_{1}, \hdots, x_{n}$ be the natural coordinates on $\mathbb{A}^{n}_{R}$, and $z_{1}, \hdots, z_{n} \in \mathcal{O}_{S}(U)$ their images under $\delta^{\sharp} : R[x_{1}, \hdots, x_{n}] \to \mathcal{O}_{S}(U)$. Then the partial derivative operators induced by these coordinates are related by $(\partial_{i} (\delta^{\sharp}(g)))_{s} = (\partial_{i} g)_{0}$. Denote by $f' \in \GL_{m}(\mathcal{O}_{\mathbb{A}^{n}_{R}/(x_{1}, \hdots, x_{n})^{r+1}, 0})$ the formal map whose pullback under $\delta^{\sharp}$ is given by $f$. We then have that
\begin{subequations}
\begin{align}
f \circ j &= \delta^{\sharp}(f') \circ j \\
&= f' \circ (\delta \circ j) \\
&= (e \circ \varphi_{f'})(\delta \circ j) \label{calclab1} \\
&= e(\delta \circ j, (\partial_{q} f'_{jk})(0)) \label{calclab2} \\
&= e(\delta \circ j, (\partial_{q} f_{jk})(s)) \label{calclab3} \\
&= e(\delta \circ j, \xi_{q,ij}(s, f_{jk}(s))) \label{calclab4} \\
&= e(\zeta(j, M_{0}^{-1})) .
\end{align}
\end{subequations}
where on line (\ref{calclab1}) we apply \autoref{curryprop2}; on lines (\ref{calclab2}), (\ref{calclab3}) and (\ref{calclab4}) the indices $q, j$ and $k$ range over all $q \in \mathcal{Q}^{d}_{r}$ and $1 \leq j, k \leq m$; on lines (\ref{calclab2}) and (\ref{calclab3}) we apply partial differentiation with respect to the coordinates $x_{1}, \hdots, x_{n}$ and $z_{1}, \hdots, z_{n}$, respectively; and lastly on line (\ref{calclab4}) we apply \autoref{fagreescalc}. 
\end{proof}

\begin{lem}
\label{existslem}
Suppose that $R = K$ is a field. Then for each $K$-point $(j, f_{0}) \in \mathcal{P}^{d}_{r}(K) = (J^{d}_{r} S \times \GL_{m})(K)$, there exists a unique $r$-limp $\psi = q \circ M$ relative to $v^{1}, \hdots, v^{m}$ such that the value of $M$ at $\delta^{-1}(0)$ is equal to $f^{-1}_{0}$. 
\end{lem}

\begin{proof}
Let $s \in S(K)$ be the image of $j$ under the projection $J^{d}_{r} S \to S$. After replacing $\delta : S \to \mathbb{A}^{n}_{R}$ with its translate by $\delta(s)$ we may assume that $s$ lies in the fibre $\delta^{-1}(0)$, and after replacing $S$ with an open subset $U \subset S$ we may further assume that $\delta^{-1}(0) = s$. To show existence, we will construct an appropriate matrix $M \in \GL_{m}(\mathcal{O}_{S,\delta}^{r})$ such that $M$ is the change of basis matrix between the image of the frame $v^{1}, \hdots, v^{m}$ in $\mathcal{H}^{r}_{\delta}$ and a flat frame $b^{1}, \hdots, b^{m}$ of $\mathcal{H}^{r}_{\delta}$. Using the coordinates $z_{1}, \hdots, z_{n}$ induced by $\delta$ we may define $f_{jk}$ to be the formal function whoose $z^{i_{1}} \cdots z^{i_{j}}$ coefficient is $\frac{1}{r!} \xi_{q, jk}(s, f_{0}^{-1})$, where $q \in \mathcal{Q}^{d}_{r}$ is the partition $i_{1} + \cdots + i_{j}$, and $\xi_{q,jk}$ is as in \S\ref{locconstrsec}; the well-definedness comes from the well-definedness of the functions $\xi_{q,jk}$ shown in \autoref{xiwelldef}. The same computation as in \autoref{fagreescalc} then shows that $M = f^{-1}$ gives the necessary change of basis matrix. We note that $f \in \GL_{m}(\mathcal{O}_{S,\delta}^{r})$ is invertible because $f_{0}$ is. 

For uniqueness, we note that $f$ is uniquely determined by its partial derivatives up to order $r$, which are determined by the differential system induced by $\nabla$ and the pair $(s, f_{0})$ by the calculation of \autoref{fagreescalc}, and this uniquely determines $M$. 
\end{proof}

\subsubsection{Gluing}
\label{glusec}

We now use the results of \S\ref{relpersec} to glue the construction in \S\ref{locconstrsec} to a global construction. In what follows we use repeatedly the fact that maps of reduced schemes are determined by their values on scheme-theoretic points. Note in particular that $J^{d}_{r} S$ is reduced: the ring $R$ is reduced because it is an integral domain, any smooth scheme over a reduced ring is reduced (see \cite[\href{https://stacks.math.columbia.edu/tag/033B}{Lemma 033B}]{stacks-project}), and the scheme $J^{d}_{r} S$ obtained by Weil restriction is smooth over $R$ if $S$ is. By observing that the construction of the map $\alpha$ in $\S\ref{relpersec}$ is preserved under base-change by maps $R \to R'$, this will allow us to reduce gluing arguments to applications of \autoref{evallem} and \autoref{existslem}.

\begin{lem}
The map $\alpha$ constructed in \S\ref{relpersec} is independent of the choice of trivializing sections $dz_{1}, \hdots, dz_{n}$. 
\end{lem}

\begin{proof}
We let $dz'_{1}, \hdots, dz'_{n}$ be another choice, and let $\alpha'$ be the associated map. To show that $\alpha$ and $\alpha'$ agree, it suffices by our reasoning above to show that the maps agree on scheme-theoretic points. After choosing such a point $\mathfrak{q} \in \Spec R$ and base-changing both constructions to $K = \kappa(\mathfrak{q})$, it suffices to show the maps agree on points over $\mathfrak{q}$. Base-changing again to the algebraic closure $\overline{K}$ we may reduce to checking the maps agree on $\overline{K}$-points. The result then follows from \autoref{evallem} and \autoref{existslem} above, since the notion of $r$-limp is independent of the sections trivializing $\Omega^{1}_{S/R}$. In particular, for any $(j, f_{0}) \in \mathcal{P}^{d}_{r}(\overline{K})$ we have that 
\[ \alpha(j, f_{0}) = \psi_{0} \circ j = \alpha'(j, f_{0}) , \] 
for a unique $\psi_{0}$ determined by $f_{0}$ and the the image $s$ in $S$ of the jet $j$. 
\end{proof}

\begin{lem}
\label{biggluinglem}
Suppose that $v'^{1}, \hdots, v'^{m}$ is a different filtration-compatible frame on $S$. Then there exists a natural isomorphism of the associated torsors $i_{vv'} : \mathcal{P}^{d}_{r} \xrightarrow{\sim} \mathcal{P}'^{d}_{r}$, such that $\alpha$ and $\alpha' \circ i_{vv'}$. If $v''^{1}, \hdots, v''^{m}$ is a third such frame, one has $i_{v'v''} \circ i_{vv'} = i_{vv''}$. 
\end{lem}

\begin{proof}
Let $[a_{k\ell}] \in \GL_{m}(S)$ be the matrix defined by the property that $v^{i} = \sum_{j = 1}^{m} a_{ji} v'^{i}$. Both $\mathcal{P}^{d}_{r}$ and $\mathcal{P}'^{d}_{r}$ are modelled as $J^{d}_{r} S \times \GL_{m}$, so we may take $i_{vv'}$ to be given by $\textrm{id} \times m_{A}$, where we denote by $m_{A}$ multiplication from the right by $A = [a_{k\ell}]$. To check the equality $\alpha' \circ i_{vv'}$, it suffices once again to reduce to the case where $R = \overline{K}$ is an algebraically closed field. Letting $(j, f_{0}) \in \mathcal{P}^{d}_{r}(\overline{K})$ be a point with $j$ lying above the point $s \in S(\overline{K})$, we let $\psi_{f_{0}}$ be the $r$-limp relative to $v^{1}, \hdots, v^{m}$ at $s$ determined by $f_{0}$ as in \autoref{existslem}, and similarly we let $\psi_{f'_{0}}$ be the $r$-limp relative to $v'^{1}, \hdots, v'^{m}$ at $s$ determined $f'_{0} =  f_{0} \cdot A(s)$.

It will suffice to check that $\psi_{f_{0}} = \psi_{f'_{0}}$. Recalling that $\psi_{f_{0}} = q \circ M$ and $\psi_{f'_{0}} = q \circ M'$ for change-of-frame matrices $M$ from $v^{1}, \hdots, v^{k}$ to $b^{1}, \hdots, b^{k}$ and $M'$ from $v'^{1}, \hdots, v'^{m}$ to $b'^{1}, \hdots, b'^{k}$, it suffices to check that $C M = M'$ for a matrix $C \in \GL_{m}(\mathcal{O}^{r}_{S, s})$ that preserves the fibres of $q$. Recalling the definition of $q$ in \S\ref{infpermapsec}, such a matrix must be block triangular with blocks of size $\dim F^{i} / F^{i+1} = h^{i}$. Letting $f = M^{-1}$ and $f' = M'^{-1}$ we have that 
\[ b^{k} = \sum_{i = 1}^{m} f_{ik} v^{i} = \sum_{i = 1}^{m} \left( \sum_{j = 1}^{m} a_{ji} f_{ik}  \right) v'^{i} . \]
It follows that the functions $a_{ji} f_{ik}$ also define a flat frame at $s$ in terms of the frame $v^{1}, \hdots, v^{m}$, and from the equality  \[ f'(s) = f'_{0} = f_{0} A(s) = f(s) A(s) , \]
and a uniqueness argument as in \autoref{existslem} we learn that $f' = f A$. Taking inverses, the result follows with $C = A^{-1}$. 

This completes the proof that $\alpha = \alpha' \circ i_{vv'}$. The claimed cocycle condition is immediate from the definition.
\end{proof}

\autoref{biggluinglem} finally lets us complete our definition of $\eta^{d}_{r}$ in general, allowing us to drop the assumption on $S$ we made at the beginning of \S\ref{locconstrsec}. We define $\eta^{d}_{r}$ locally on open neighbourhoods $U \subset S$ by the construction in \ref{locconstrsec}, and the statement of the lemma can be read as specifying a gluing datum for a torsor over $S$ with transition maps associated to different local trivializations given by the maps $i_{vv'}$. We define $\eta^{d}_{r}$ to be the map defined by the thus constructed global $\GL_{m}$-torsor. 

Lastly, let us give a coordinate-invariant description of the pair $(\gamma, \alpha)$ which will be useful later:

\begin{lem}
\label{Pinterp}
The bundle $\mathcal{P}^{d}_{r}$ is nothing other than the base-change along the map $J^{d}_{r} S \to S$ of the frame bundle associated to $\mathcal{H} \to S$. Thus a point in $\mathcal{P}^{d}_{r}(R)$ may be identified with a pair $(j, \iota)$, where $j \in (J^{d}_{r} S)(R)$ lies above $s \in S(R)$ and $\iota$ is an isomorphism $\mathcal{H}_{s} \xrightarrow{\sim} R^{m}$. Moreover, if $\psi = q \circ M$ is an $r$-elimp defined at $s$ such that $M(s)$ is the change-of-basis matrix between a filtration-compatible frame and the flat basis determined by $\iota$, then $\alpha(j,\iota) = \psi \circ j$.
\end{lem}

\begin{proof}
Immediate from the construction and the arguments of \autoref{evallem} and \autoref{biggluinglem} above.
\end{proof}

\begin{defn}
\label{framedlimpdef}
A pair $(\psi, \iota)$, where $\psi = q \circ M$ is an elimp defined at $s$ and $M(s)$ is the change-of-basis matrix from a filtration compatible frame to the frame determined by $\iota : \mathcal{H}_{s} \xrightarrow{\sim} R^{m}$, we will called a ``framed elimp''. We will also use the similar terminology ``framed limp'', ``framed $r$-limp'', ``framed local period map'', etc. 
\end{defn}

\subsubsection{Key Properties}
\label{keypropsec}

We now check that the listed properties in \autoref{bigjetconstr} hold. Most of this is essentially immediate from the construction, though we provide some details. 

Beginning with (i), suppose that $(\gamma, \alpha)$ defines $\eta^{d}_{r}$ and $(\gamma', \alpha')$ defines $\eta'^{d}_{r}$. Let $g : S' \to S$ be the base-change map. We start by observing the existence of a natural map $\xi : \mathcal{P}'^{d}_{r} \to \mathcal{P}^{d}_{r}$. Given a filtration-compatible frame $v^{1}, \hdots, v^{m}$ over a Zariski open subset $U \subset S$, the pullback $g^{*} v^{1}, \hdots, g^{*} v^{m}$ is a filtration-compatible frame for $(\mathcal{H}', F'^{\bullet}, \nabla')$. Letting $U' = g^{-1}(U)$ and using the local models $\restr{\mathcal{P}^{d}_{r}}{U'} \simeq J^{d}_{r} U' \times \GL_{m, R'}$ and $\restr{\mathcal{P}^{d}_{r}}{U} \simeq J^{d}_{r} U \times \GL_{m, R}$ associated to these frames, we define $\restr{\xi}{U'} = \restr{(J^{d}_{r} g)}{U'} \times i_{R'}$, where $i_{R'} : \GL_{m,R'} \to \GL_{m,R}$ is the base-change map. That $\xi$ is well-defined is immediate from the constructions of $\mathcal{P}'^{d}_{r}$ and $\mathcal{P}^{d}_{r}$. One easily checks that
\begin{equation}
\label{jetbcdiag}
\begin{tikzcd}
\mathcal{P}'^{d}_{r} \arrow[r, "\xi"] \arrow[d, "\delta'"] & \mathcal{P}^{d}_{r} \arrow[d, "\delta"] \\
J^{d}_{r} S' \arrow[r, "J^{d}_{r} g"] & J^{d}_{r} S
\end{tikzcd}
\end{equation}
is a fibre product diagram. To prove (i), it then suffices to check that $\alpha' = \alpha \circ \xi$. This can be done locally on $S$, so we can fix a map $\delta : U \to \mathbb{A}^{n}_{R}$ with base-change $\delta' : U \to \mathbb{A}^{n}_{R'}$. One then easily checks that, in the notation of \S\ref{locconstrsec} above, we have $\zeta_{R'} \circ (J^{d}_{r} g \times i_{R'}) = \zeta'$, where $\zeta$ and $\zeta'$ are the maps given by (\ref{zetadef}), hence the result. 

Next we consider (ii). The first part of the verification, which includes the construction of $\xi$ and verifying that the diagram (\ref{jetbcdiag}) is a fibre product diagram, is essentially identical to (i), the only difference being that we define $\restr{\xi}{U'} = \restr{(J^{d}_{r} g)}{U'} \times \textrm{id}$. To verify that $\alpha' = \alpha \circ \xi$ we once again reduce to the case of points over an algebraically-closed field $\overline{K}$. The statement is local, so we may work with local models $\restr{\mathcal{P}^{d}_{r}}{U} \simeq J^{d}_{r} U \times \GL_{m,\overline{K}}$ and $\restr{\mathcal{P}^{d}_{r}}{U'} \simeq J^{d}_{r} U' \times \GL_{m,\overline{K}}$ associated to the frames $v^{1}, \hdots, v^{m}$ and $g^{*} v^{1}, \hdots, g^{*} v^{m}$. Consider a point $(j', f'_{0}) \in (J^{d}_{r} U')(\overline{K}) \times \GL_{m}(\overline{K})$ and let $s'$ be the image of $j'$ in $U'(\overline{K})$. Let $s = g(s')$, and let $\psi = q \circ M$ be the $r$-limp defined at $s$ with $M(s)^{-1} = f'_{0}$ obtained by from \autoref{existslem}. Then $\psi \circ g$ is an $r$-limp the defined at $s'$ by the initial condition $f'_{0}$. Applying \autoref{evallem} we have
\[ \alpha'(j', f'_{0}) = (\psi \circ g) \circ j' = \psi \circ (g \circ j') = \alpha(g \circ j', f'_{0}) = \alpha(\xi(j', f'_{0})) . \]

Next we prove (iii); we begin by explaining the statement. The map $v : \mathbb{D}^{d}_{r, R} \to \mathbb{D}^{d'}_{r', R}$ induces, for any $R$-scheme $X$, a map of functors
\begin{equation}
\label{vdef}
\Hom_{R}((-) \times_{R} \mathbb{D}^{d'}_{r', R}, X) \to \Hom_{R}((-) \times_{R} \mathbb{D}^{d}_{r, R}, X) ,
\end{equation}
by precomposition. Taking $X = \ch{L}_{R}$ we obtain from representability a map $(-) \circ v: J^{d'}_{r'} \ch{L}_{R} \to J^{d}_{r} \ch{L}_{R}$. This map commutes with the actions of $\GL_{m,R}$ on $J^{d'}_{r'} \ch{L}_{R}$ and $J^{d}_{r} \ch{L}_{R}$, and so induces a map $\GL_{m,R} \backslash J^{d'}_{r'} \ch{L}_{R} \to \GL_{m,R} \backslash J^{d}_{r} \ch{L}_{R}$ of quotient stacks. We describe this latter map explicitly on points. A point of $\GL_{m,R} \backslash J^{d'}_{r'} \ch{L}_{R}$ over the $R$-scheme $T$ is a pair $(\ep, \beta)$ consisting of a $\GL_{m,R}$-torsor $\ep : \mathcal{T} \to T$ and a $\GL_{m,R}$-invariant map $\beta : \mathcal{T} \to J^{d'}_{r'} \ch{L}$. Its image in $\GL_{m,R} \backslash J^{d}_{r} \ch{L}_{R}$ is then the point associated to the pair $(\ep, \beta \circ v)$.

From this discussion we see that (iii) amounts to the statement that if $(\gamma', \alpha')$ is the pair corresponding to $\eta^{d'}_{r'}$ and $(\gamma, \alpha)$ is the pair corresponding to $\eta^{d}_{r}$, then the pair $(\gamma', \alpha' \circ v)$ is isomorphic to the pullback $((-) \circ v)^{*}(\gamma, \alpha)$ of the pair along the map $(-) \circ v : J^{d'}_{r'} S \to J^{d}_{r} S$ induced by (\ref{vdef}) with $X = S$. That the torsor $\gamma$ pulls back to the torsor $\gamma'$ is immediate from its local construction and we do not repeat this verification. That $\alpha$ pulls back to $\alpha' \circ v$ may be checked on points in a similar fashion to our arguments in (i) and (ii), where we reduce to considering everything over an algebraically closed field $\overline{K}$ and then evaluating the respective maps by constructing appropriate $r$-limps using \autoref{existslem}. 

Finally, (iv) is simply a coarser version of \autoref{evallem} which we have used throughout this verification.

\section{Additional Jet-Theoretic Constructions}

This is a short section describing some additional constructions associated to jets and jet schemes that we will find useful in our main arguments.

\begin{defn}
Given an $R$-scheme $X$ with associated jet-scheme $J^{d}_{r} X$ over $R$, the \emph{constant} subscheme $c(X) \subset J^{d}_{r} X$ is the image of the section $c : X \to J^{d}_{r} X$ which is induced on the level of functors by the map
\[ \Hom_{R}(T, X) \to \Hom_{R}(T \times_{R} \mathbb{D}^{d}_{r,R}, X) , \]
which associates to the map $T \to X$ the map $T \times_{R} \mathbb{D}^{d}_{r} \to T \to X$ coming from the natural projection $A^{d}_{r, R} \to R$.
\end{defn}

\begin{defn}
For later use, we define $J^{d}_{r,nc} X \subset J^{d}_{r} X$ to be the open subscheme of ``non-constant'' jets, i.e., the open subscheme corresponding to the complement $J^{d}_{r} X \setminus c(X)$. 
\end{defn}

In coordinates, the map $c$ amounts to associating to a point $x \in X(R)$ the map from $\mathbb{D}^{d}_{r,R}$ landing at $x$ whose partial derivatives are all constant, and doing so functorially in $R$. We use it to associate to any family $u : \mathcal{Y} \to \mathcal{M}$ over $R$ an associated family $u^{d}_{r} : \mathcal{Y}^{d}_{r} \to \mathcal{M}$, as follows: 

\begin{lem}
\label{jetfams}
Suppose that $u : \mathcal{Y} \to \mathcal{M}$ is a map of $R$-schemes, and let $\mathcal{Y}^{d}_{r} \subset J^{d}_{r} Y$ be the fibre above $c(\mathcal{M}) \subset J^{d}_{r} \mathcal{M}$, and denote by $u^{d}_{r} : \mathcal{Y}^{d}_{r} \to \mathcal{M}$ the induced projection. Then if $\mathfrak{q} \in \mathcal{M}$ is a scheme-theoretic point, the natural inclusion $J^{d}_{r} \mathcal{Y}_{m} \hookrightarrow \mathcal{Y}^{d}_{r,m}$ is an isomorphism, where the fibres are taken with respect to the map $m: \kappa(\mathfrak{q}) \to \mathcal{M}$.
\end{lem}

\begin{rem}
The statement of \autoref{jetfams} may be interpreted on the level of functors of points; that is, one does not have to worry whether the schemes $J^{d}_{r} \mathcal{Y}$ and $J^{d}_{r} \mathcal{M}$ exist or not.
\end{rem}

\begin{proof}
A map $T \to \mathcal{Y}^{d}_{r,m}$ of schemes over $\kappa(\mathfrak{q})$ is the data of a map $\sigma : T \times_{\kappa(\mathfrak{q})} \mathbb{D}^{d}_{r,\kappa(\mathfrak{q})} \to \mathcal{Y}$ which:
\begin{itemize}
\item the map $\sigma_{2} : T \times_{\kappa(\mathfrak{q})} \mathbb{D}^{d}_{r,\kappa(\mathfrak{q})} \to \mathcal{M}$ obtained as $\sigma_{2} = u \circ \sigma$, is constant; that is to say, the map $\sigma_{2}$ factors as $T \times_{\kappa(\mathfrak{q})} \mathbb{D}^{d}_{r,\kappa(\mathfrak{q})} \xrightarrow{\ep} T \xrightarrow{\xi} \mathcal{M}$ via the structure morphism $\mathbb{D}^{d}_{r,\kappa(\mathfrak{q})} \to \Spec \kappa(\mathfrak{q})$;
\item and for which the map $\xi$, which also arises by first pulling back $\sigma$ along the inclusion $\iota : T \hookrightarrow T \times_{\kappa(\mathfrak{q})} \mathbb{D}^{d}_{r,\kappa(\mathfrak{q})}$ and then composing with $u$, in fact factors as $\xi = m \circ \tau$, where $\tau : T \to \Spec \kappa(\mathfrak{q})$ is the structure morphism.
\end{itemize}
Combining these two points, we obtain the fibre product diagram
\begin{center}
\begin{tikzcd}
\mathbb{D}^{d}_{r, \kappa(\mathfrak{q})} \times_{\kappa(\mathfrak{q})} T
\arrow[bend left]{drr}{\tau \circ \ep}
\arrow[bend right,swap]{ddr}{\sigma}
\arrow[dashed]{dr} & & \\
& \mathcal{Y}_{m} \arrow{r} \arrow[d, hook]
& \Spec \kappa(\mathfrak{q}) \arrow{d}{m} \\
& \mathcal{Y} \arrow[swap]{r}{u}
& \mathcal{M}
\end{tikzcd} ,
\end{center}
showing that $\sigma$ corresponds to a point of $J^{d}_{r} \mathcal{Y}_{m}$. 
\end{proof}

We will also later find that the following subclass of jets will be useful:

\begin{defn}
For a $R$-scheme $X$, we define $J^{d}_{r,nd} \subset J^{d}_{r} X$ to be the open subscheme of ``non-degenerate jets'', defined as follows:
\begin{itemize}
\item for $r = 0$ we have $J^{d}_{0,nd} X = X$;
\item for $r = 1$ we have a natural projection $\nu : J^{d}_{1} X \to (T X)^{d}$ to the $d$-times self-product of the relative tangent bundle induced by the natural map $R[t_{1}, \hdots, t_{d}]/(t_{1}, \hdots, t_{d})^{2} \to R[t_{1}]/(t_{1})^2 \times \cdots \times R[t_{d}]/(t_{d})^2$, and we define $J^{d}_{1,nd} X$ to be the fibre above the locus of $d$-tuples of linearly independent vectors above a common point of $X$;
\item for $r \geq 2$ we let $J^{d}_{r,nd} X$ be the fibre above $J^{d}_{1,nd} X$. 
\end{itemize}
\end{defn}
\noindent The formation of the scheme $J^{d}_{r,nd} X$ is functorial in $R$. Its utility comes from the following two lemmas, which will let us use ``families'' of non-degenerate jets to study (possibly formal) local period maps. 

\begin{defn}
Suppose that we have a collection $\{ Y_{r} \}_{r \geq 0}$ of schemes over $R$ and maps $\pi_{r+1} : Y_{r+1} \to Y_{r}$ of $R$-schemes. Then by a \emph{compatible sequence} of points in $\{ Y_{r} \}_{r \geq 0}$ we mean a sequence $\{ \sigma_{r} \}_{r \geq 0}$ with $\sigma_{r} \in Y_{r}(R)$ and $\pi_{r+1}(\sigma_{r+1}) = \sigma_{r}$ for all $r \geq 0$. The two primary cases of interest are:
\begin{itemize}
\item $Y_{r}$ is a subscheme of $J^{d}_{r} X$ for some $R$-scheme $X$, and the projections are obtained by restricting the projections $J^{d}_{r+1} X \to J^{d}_{r} X$;
\item $Y_{r}$ is a subscheme of the torsor $\mathcal{P}^{d}_{r}$ constructed in \S\ref{glusec}, and the projections are obtained by restricting the projections $\mathcal{P}^{d}_{r+1} \to \mathcal{P}^{d}_{r}$.
\end{itemize}
\end{defn}

\begin{lem}
\label{existscompseqlem}
Suppose that we have a collection $\{ Y_{r} \}_{r \geq 0}$ of $\mathbb{C}$-schemes with maps $\pi_{r+1} : Y_{r+1} \to Y_{r}$, and define $\pi^{r'}_{r} : Y_{r'} \to Y_{r}$ for $r' \geq r$ to be the composition $\pi_{r+1} \circ \cdots \circ \pi_{r'}$. Let $\{ \mathcal{T}_{r} \}_{r \geq 0}$ be a family of constructible subsets, with $\mathcal{T}_{r} \subset Y_{r}$, and such that $\pi_{r+1}(\mathcal{T}_{r+1}) \subset \mathcal{T}_{r}$. Then if $\pi^{r}_{0}(\mathcal{T}_{r})$ is non-empty for each $r$ and $j_{0} \in \bigcap_{r} \pi^{r}_{0}(\mathcal{T}_{r})$, there exists a compatible sequence $\{ j_{r} \}_{r \geq 0}$ with $j_{r} \in \mathcal{T}_{r}$ for all $r$. 
\end{lem}

\begin{proof}
See \cite[Lemma 5.3]{periodimages}. 
\end{proof}

\begin{lem}
\label{jetfactorlem}
Suppose that $\psi : B \to Y$ is a map of analytic spaces with $B$ irreducible of dimension $d$, and that there exists a compatible sequence $\{ j_{r} \}_{r \geq 0}$ with $j_{r} \in J^{d}_{r,nd} B$ such that $\psi \circ j_{r}$ lies in $J^{d}_{r} Z$ for all $r$, where $Z \subset Y$ is a closed analytic subvariety. Then $\psi(B) \subset Z$.
\end{lem}

\begin{proof}
This is \cite[Lemma 4.5]{urbanik2021sets}.
\end{proof}

\section{Subvarieties of Flags Defined by Tensors}
\label{subflagsec}

Let $R$ be a ring. In this section we study algebraic loci in $\ch{L}_{R}$ which are defined by pairs $(F^{\bullet}, U)$, where $F^{\bullet}$ is a flag on $R^{m}$ of weight $w$, and for which $U$ is an $R$-module of ``Hodge-like'' elements --- let us explain what we mean. 

\vspace{0.5em}

\begin{notn}
Given a free $R$-module $V$ for $R$ a ring, and integers $a, b \geq 0$, we write
\[ V^{a,b} = V^{\otimes a} \otimes (V^{*})^{\otimes b} . \]
If $F^{\bullet}$ is a filtration on $V$, there is a natural induced filtration on $V^{a,b}$ (see \cite[\S1.1]{PMIHES_1971__40__5_0}), which we will also denote by $F^{\bullet}$. (We also describe the filtration on $V^{*}$ explicitly in the proof of \autoref{polisHodge} below, and the filtrations on the various tensor products are just obtained by taking tensor products of the filtered pieces.)
\end{notn}

\vspace{0.3em}

\begin{notn} ~
\begin{itemize}
\item[-] If $F^{\bullet}$ is any filtration of even weight $2j$ we write $F^{\textrm{mid}}$ for $F^{j}$; we do this in situations where the weight of $F^{\bullet}$ is unimportant to avoid introducing unnecessary notation.
\item[-] If $F^{\bullet}$ is a filtration of odd weight we set $F^{\textrm{mid}} = 0$.
\item[-] Given a free $R$-module $V$, we will also use the notation $F^{\textrm{mid}}$ to denote the sum
\[ \bigoplus_{a,b \geq 0} F^{\textrm{mid}} V^{a,b} , \]
with the intended meaning being clear from context.
\end{itemize}
\end{notn}

\vspace{0.3em}

\begin{defn}
\label{oppdef}
We say that filtrations $F^{\bullet}$ and $F'^{\bullet}$ of the same weight $w$ on $V$ are \emph{opposed} if $F^{p} \oplus F'^{w-p+1} = V$ for all $p$. 
\end{defn}

\vspace{0.5em}

The situation of interest is then as follows. We will study loci in $\ch{L}_{R}$ defined by pairs $(F^{\bullet} , U)$ for which the following property holds:

\begin{itemize}
\item[(M)] There exists a filtration $F^{\bullet}_{c}$, opposed to $F^{\bullet}$, and the submodule $U \subset \bigoplus_{a, b \geq 0} (R^{m})^{a, b}$ is spanned over $R$ by elements of $F^{\textrm{mid}} \cap F^{\textrm{mid}}_{c}$. 
\end{itemize}

\noindent There are two primary situations where we are interested in this setup:

\begin{itemize}
\item[(1)] The Hodge theoretic setting: here we have $R = \mathbb{C}$, the flag $F^{\bullet}$ is the Hodge flag for some polarizable Hodge structure of weight $w$ on $\mathbb{R}^{m}$, and we have $F^{\bullet}_{c} = \overline{F}^{\bullet}$.
\item[(2)] The positive-characteristic algebraic de Rham setting: here $R = \kappa$ is a field of positive characteristic $p$, and the filtrations $F^{\bullet}$ and $F^{\bullet}_{c}$ are obtained by transferring the \emph{Hodge} and \emph{conjugate} filtration along an identification $\kappa^{m} \cong H^{w}_{\textrm{dR}}(Y)$ for $Y$ a smooth projective algebraic variety over $\kappa$. 
\end{itemize}

\noindent Given such a conjugate filtration $F^{\bullet}_{c}$, we will write $H^{p,q}$ for the intersections $F^{p} \cap F^{q}_{c}$. That $F^{\bullet}_{c}$ is opposed to $F^{\bullet}$ means that $\sum H^{p, w-p}$ is a direct sum decomposition of $R^{m}$.

\subsection{Polarizations}

As it will be the case in the situations we consider, we are particularly intereted in the situation where we have non-degenerate bilinear form $Q : R^{m} \otimes R^{m} \to R$, either symmetric or alternating, and our flag $F^{\bullet}$ is \emph{polarized} by $Q$. This means it is a point of the $R$-subscheme $\ch{L}_{\textrm{pol}} \subset \ch{L}_{R}$ defined by the relation
\begin{equation}
\label{poleq}
 Q(F^{p}, F^{w-p+1}) = 0 , \hspace{2em} 0 \leq p \leq w+1 .
\end{equation}
\noindent For a flag $F^{\bullet} \in \ch{L}_{\textrm{pol}}(R)$, we will also be interested in the same condition on the conjugate filtration $F^{\bullet}_{c}$, i.e., we will be interested in the situation where additionally $F^{\bullet}_{c} \in \ch{L}_{\textrm{pol}}(R)$. In situation (1) described above this is typically automatic because the polarization $Q$ is defined over $\mathbb{R}$ and $F^{\bullet}_{c} = \overline{F}^{\bullet}$. In situation (2) however this must be explicitly checked, and is done so in our situation of interest in \autoref{polposchar} below. We note that these ``polarization conditions'' can in fact be subsumed by the discussion of elements lying inside $F^{\textrm{mid}} \cap F^{\textrm{mid}}_{c}$, as the following lemma shows:

\begin{lem}
\label{polisHodge}
Regard $Q$ as an element of $(R^{m})^{0,2}$. For a flag $F^{\bullet} \in \ch{L}(R)$, the condition that $F^{\bullet} \in \ch{L}_{\textrm{pol}}(R)$ is equivalent to the condition that $Q \in F^{\textrm{mid}}$. Consequently, for a flag $F^{\bullet} \in \ch{L}(R)$ with conjugate filtration $F^{\bullet}_{c}$, condition that both $F^{\bullet}$ and $F^{\bullet}_{c}$ lie in $\ch{L}_{\textrm{pol}}(R)$ is the same as saying $Q \in F^{\textrm{mid}} \cap F^{\textrm{mid}}_{c}$. 
\end{lem}

\begin{proof}
It suffices to check the first statement. The filtration on $(R^{m})^{*}$ has weight $-w$ and is defined via 
\[ F^{i} (R^{m})^{*} = \{ \varphi \in (R^{m})^{*} : \varphi(F^{-i+1}) = 0 \} . \]
Taking the self tensor-product, we find that 
\[ F^{\textrm{mid}} (R^{m})^{0,2} = \sum_{i} (F^{i} (R^{m})^{*}) \otimes (F^{-w-i} (R^{m})^{*}) , \]
hence if $Q \in F^{\textrm{mid}} (R^{m})^{0,2}$ we must have that $Q(F^{p}, F^{w-p+1}) = 0$ for $0 \leq p \leq w+1$. 
\end{proof}

\noindent With \autoref{polisHodge} in mind, we may now view the condition that our flag $F^{\bullet}$ with its conjugate flag be polarized as a special case of the condition (M) by requiring that $Q$ be an element of $U$.

\begin{rem}
In Hodge theory, saying that a Hodge structure is \emph{polarized} typically requires an additional positivity condition in addition to the condition \autoref{polisHodge} on the Hodge flag. When we refer to polarized \emph{flags} (instead of Hodge structures), we will always mean the weaker condition \autoref{polisHodge}.
\end{rem}

Suppose now that $R = K$ is a field. In the context of a fixed polarization $Q$, we will also often denote by $Q$ the induced bilinear form on the spaces $(K^{m})^{a,b}$ for all $a, b \geq 0$, with the intended meaning inferred from context. We will denote by $(-)^{\dagger}$ the adjunction with respect to $Q$; i.e., for $\xi \in \textrm{End}(K^{m})$ we define $\xi^{\dagger}$ via
\[ Q(\xi v, w) = Q(v, \xi^{\dagger} w), \hspace{2em} \textrm{for all }v, w \in K^{m} , \]
and with the analogous definition for the operator $(-)^{\dagger}$ on the higher tensor spaces $(K^{m})^{a,b}$. We may also extend $(-)^{\dagger}$ to a map $(K^{m})^{a,b} \xrightarrow{\sim} (K^{m})^{b,a}$ for each $a, b$, with $w^{\dagger} = Q(-,w)$. 





\subsection{Associated Orbits}
\label{orbitsec}

Fix a pair $(F^{\bullet}, U)$ as before, and let $\mathbf{G} \subset \GL_{m,R}$ be the algebraic subscheme stabilizing the elements of $U$. For the remainder of this section we will be interested in comparing two different subschemes of $\ch{L}$, which we denote by $O(F^{\bullet}, U)$ and $O(F^{\bullet}, \mathbf{G})$, induced by the pair $(F^{\bullet}, U)$. They are defined as follows:

\begin{defn}
We denote by $O(F^{\bullet}, U) \subset \ch{L}_{R}$ the irreducible component containing $F^{\bullet}$ of the $R$-algebraic subscheme whose $T$-points $F'^{\bullet}$ satisfy $U_{T} \subset F'^{\textrm{mid}}$, where $T$ is any scheme lying over $\Spec R$. 
\end{defn}

\begin{defn}
We define $O(F^{\bullet}, \mathbf{G})$ to be the orbit $\mathbf{G}^{\circ} \cdot F^{\bullet} \subset \ch{L}_{R}$, where $\mathbf{G}^{\circ} \subset \mathbf{G}$ is the identity component.
\end{defn}

\begin{prop}
\label{ABCProp}
Suppose that $R = K$ is an algebraically closed field of characteristic not equal to $2$, and fix a polarization $Q : K^{m} \otimes K^{m} \to K$. Suppose that the pair $(F^{\bullet}, U)$ satisfies property (M) above, that $U$ contains $Q$, and write $H^{p,q} = F^{p} \cap F^{q}_{c}$. Denote by $\beta : \mathbb{G}_{m,K} \to \GL_{m,K}$ the cocharacter defined by the property:
\begin{itemize}
\item[-] if the weight $w$ is odd, then $\beta$ acts as $z^{p}$ on $H^{p,w-p}$; and
\item[-] if the weight $w$ is even, then $\beta$ acts as $z^{p-w/2}$ on $H^{p,w-p}$.
\end{itemize}
Then we have
\begin{itemize}
\item[(i)] $O(F^{\bullet}, U) = \overline{O(F^{\bullet}, \mathbf{G})}^{\textrm{Zar}}$; 
\item[(ii)] the weight space decomposition of the adjoint cocharacter $\textrm{Ad}\, \beta$ is given by
\begin{equation}
\label{endiieq}
\textrm{End}(K^{m})^{i} = \sum_{p} \textrm{Hom}(H^{p, w-p}, H^{p-i, w-p+i}) ;
\end{equation}
\item[(iii)] the character $\beta$ factors through $\mathbf{G}$, and therefore induces a direct sum decomposition $\mathfrak{g} = \bigoplus_{i} \mathfrak{g}^{i}$ compatibly with (\ref{endiieq}); and
\item[(iv)] the natural map $\bigoplus_{i > 0} \mathfrak{g}^{i} \to T_{F^{\bullet}} O(F^{\bullet}, U)$ is an isomorphism.
\end{itemize}
\end{prop}

\begin{proof}
We begin by noting that the statement holds in the special case where $U = \textrm{span} \{ Q \}$. Because all symmetric (resp. alternating) bilinear forms over the algebraically closed field $K$ of characteristic not equal to $2$ are equivalent, we can reduce to the case where either $K = \overline{\mathbb{Q}}$ or $K = \overline{\mathbb{F}_{p}}$ with $p > 2$. The second case we can lift to characteristic zero, reducing to the first case, which is standard; see for instance \cite[\S3]{Schmid1973}. In what follows we denote by $\mathfrak{a} = \bigoplus_{i} \mathfrak{a}^{i}$ the Lie algebra decomposition associated to the Lie algebra $\mathfrak{a}$ of $\textrm{Aut}(K^{m}, Q)$.

Proceeding now to the general case, observe first that $O(F^{\bullet}, \mathbf{G}) \subset O(F^{\bullet}, U)$. As $O(F^{\bullet}, \mathbf{G})$ is an orbit of a connected algebraic group it is smooth and locally closed. Moreover, $O(F^{\bullet}, U)$ is (assumed to be) irreducible. For the first claim (i) it therefore suffices to show that the inclusion $T_{F^{\bullet}} O(F^{\bullet}, \mathbf{G}) \subset T_{F^{\bullet}} O(F^{\bullet}, U)$ of tangent spaces is an isomorphism.

We recall that the tangent space at $W \subset K^{m}$ to the Grassmannian $\textrm{Gr}_{\ell}(K^{m})$ of all $\ell$-dimensional subspaces of $K^{m}$, with $\ell = \dim W$, is naturally identified with $\textrm{Hom}(W, K^{m}/W)$, where the natural surjection $\textrm{End}(K^{m}) \to \textrm{Hom}(W, K^{m}/W)$ is the derivative at $(\textrm{id}, W)$ of the natural orbit map $\GL_{m,K} \times \textrm{Gr}_{\ell}(K^{m}) \to \textrm{Gr}_{\ell}(K^{m})$. Embedding $\ch{L}_{K}$ in the natural way into a product of Grassmannians, we may identify $T_{F^{\bullet}} \ch{L}_{K}$ with the image of $\textrm{End}(K^{m})$ inside $\bigoplus_{p} \Hom(F^{p}, K^{m}/F^{p})$. The subspace $T_{F^{\bullet}} O(F^{\bullet}, U)$ is identified with the image inside $\bigoplus_{p} \Hom(F^{p}, K^{m}/F^{p})$ of those elements $\xi \in \textrm{End}(K^{m})$ which additionally satisfy $v \in \xi(F^{\textrm{mid}})$ for every $v \in U$, and because $O(F^{\bullet}, U) \subset O(F^{\bullet}, \textrm{span} \{ Q \})$ it suffices to consider just those $\xi \in \mathfrak{a}$. For such $\xi$, one sees using adjunction that this condition is equivalent to the condition that $\xi (v) \in F^{\textrm{mid}}$.

We now give an interpretation of $T_{F^{\bullet}} O(F^{\bullet}, \mathbf{G})$. By differentiating the natural orbit map $\mathbf{G} \times O(F^{\bullet}, \mathbf{G}) \to O(F^{\bullet}, \mathbf{G})$ at $\{ \textrm{id} \} \times \{ F^{\bullet} \}$ we obtain a natural surjection $\mathfrak{g} \to T_{F^{\bullet}} O(F^{\bullet}, \mathbf{G})$ from the Lie algebra $\mathfrak{g}$ of $\mathbf{G}$, which is compatible with the map $\textrm{End}(K^{m}) \to \bigoplus_{p} \Hom(F^{p}, K^{m}/F^{p})$ and the inclusions 
\[ T_{F^{\bullet}} O(F^{\bullet}, \mathbf{G}) \hookrightarrow T_{F^{\bullet}} O(F^{\bullet}, U) \hookrightarrow \bigoplus_{p} \Hom(F^{p}, K^{m}/F^{p}) . \]
We observe that for $i \leq 0$ the map $\textrm{End}(K^{m}) \to \bigoplus_{p} \Hom(F^{p}, K^{m}/F^{p})$ sends $\textrm{End}(K^{m})^{i}$ to zero, from which it follows that we may identify $T_{F^{\bullet}} \ch{L}_{K}$ with the image of $\bigoplus_{i > 0} \textrm{End}(K^{m})^{i}$ in $\bigoplus_{p} \Hom(F^{p}, K^{m}/F^{p})$. Note that if we have $v \in (F^{\textrm{mid}} (K^{m})^{a,b}) \cap (F^{\textrm{mid}}_{c} (K^{m})^{a,b})$ and $\xi \in \bigoplus_{i > 0} \textrm{End}(K^{m})^{i}$ such that $\xi(v) \in F^{\textrm{mid}}$, then because $\xi(v)$ lies inside $F^{\textrm{mid}} \cap \xi(F^{\textrm{mid}}_{c}) = 0$ one in fact has that $\xi(v) = 0$. We thus obtain that
\begin{align*}
T_{F^{\bullet}} O(F^{\bullet}, U) &= \textrm{im}\left( \left( \bigoplus_{i} \mathfrak{a}^{i} \right) \cap  \left\{ \xi : \xi(v) \in F^{\textrm{mid}} \textrm{ for all } v \in U \right\} \right) \\
&= \textrm{im}\left( \left( \bigoplus_{i > 0} \mathfrak{a}^{i} \right) \cap \left\{ \xi : \xi(v) = 0 \textrm{ for all } v \in U \right\} \right) \\
&\subset \textrm{im}(\mathfrak{g}) \\
&= T_{F^{\bullet}} O(F^{\bullet}, \mathbf{G}) ,
\end{align*}
which proves the first claim.

The second claim (ii) is immediate from the definition. To show (iii), begin by noting that the polarization $Q$ induces a $\mathbf{G}$-equivariant isomorphism $K^{m} \xrightarrow{\sim} (K^{m})^{*}$, and so it follows that $\mathbf{G}$ can be defined entirely by invariants inside 
\begin{itemize}
\item[-] $\bigoplus_{j \geq 1} (K^{m})^{j,j}$, if $w$ is odd; or
\item[-] $\bigoplus_{j \geq 1} (K^{m})^{\otimes j}$, if $w$ is even. 
\end{itemize}
But $\beta$ fixes elements in these spaces which additionally lie in $F^{\textrm{mid}} \cap F^{\textrm{mid}}_{c}$, meaning that $\beta$ factors through $\mathbf{G}$. Finally, (iv) was shown in the course of the proof above since $T_{F^{\bullet}} O(F^{\bullet}, \mathbf{G}) = \textrm{im}(\bigoplus_{i > 0} \mathfrak{g}^{i})$.
\end{proof}




\begin{defn}
\label{symmhodgedef}
Given a triple $(F^{\bullet}, F^{\bullet}_{c}, U)$ as in \autoref{ABCProp}, we call the integers $h^{i} = \dim (\mathfrak{g} \cap \textrm{End}(K^{m})^{i})$ the \emph{adjoint Hodge numbers}.
\end{defn}

In the setting of \autoref{ABCProp}, we also have the following two additional facts concerning the symmetry of the adjoint Hodge numbers:

\begin{lem}
\label{nondegenimpsymlem}
Suppose that the Killing form $B$ on $\mathfrak{g} = \bigoplus_{i} \mathfrak{g}^{i}$ is non-degenerate. Then the adjoint Hodge numbers are symmetric.
\end{lem}

\begin{proof}
Consider a non-zero element $x \in \mathfrak{g}^{i}$. Because $B$ is non-degenerate, there is some $x' \in \mathfrak{g}$ such that $\textrm{tr}(\textrm{ad}\, x' \circ \textrm{ad}\, x) \neq 0$. Write $x' = \sum_{j} x'_{j}$ compatibly with the grading. Then unless $j = -i$, the map $\textrm{ad}\, x'_{j} \circ \textrm{ad}\, x$ is nilpotent, hence has zero trace. It follows that $B(x'_{-i}, x) \neq 0$. This implies that the map $\mathfrak{g}^{i} \to (\mathfrak{g}^{-i})^{*}$ induced by $B$ is injective, hence $\dim \mathfrak{g}^{i} \leq \dim \mathfrak{g}^{-i}$. By a symmetric argument one obtains $\dim \mathfrak{g}^{-i} \leq \dim \mathfrak{g}^{i}$, hence the desired equality.
\end{proof}

\begin{cor}
\label{semsimpimpsymlem}
Suppose that $\mathfrak{g}$ is semisimple, and that the field $K$ over which $\mathfrak{g}$ is defined is either of characteristic zero or of positive characteristic $p \geq m+2$. Then the adjoint Hodge numbers of $\mathfrak{g}$ are symmetric.
\end{cor}

\begin{proof}
In characteristic zero, it is well-known that the semisimplicity condition is equivalent to the non-degeneracy of the Killing form. In positive characteristic, this is shown in \cite{257785} for Lie subalgebras $\mathfrak{g} \subset \mathfrak{gl}_{m}$ provided that $p \geq m+2$. In both cases we thus reduce to \autoref{nondegenimpsymlem} above.
\end{proof}



\subsection{Orbits from Global Sections}

Consider a triple $(\mathcal{H}, F^{\bullet}, \nabla)$ as in the statement of \autoref{bigjetconstr}. We now associate to each point $s \in S(R)$ a pair $(F^{\bullet}_{s}, \mathcal{U}_{s})$, and relate the $O(F^{\bullet}_{s}, \mathcal{U}_{s})$ to the images of $r$-limps. We define $\mathcal{H}^{a,b} = \mathcal{H}^{\otimes a} \otimes (\mathcal{H}^{*})^{\otimes b}$ for each pair $(a, b)$ of non-negative integers. The filtration and connection on $\mathcal{H}$ naturally induces the same data on $\mathcal{H}^{a,b}$. We denote by $\mathcal{F}^{a,b} \subset \mathcal{H}^{a,b}$ the subbundle of $\nabla^{a,b}$-flat sections, where $\nabla^{a,b}$ is the natural connection on $\mathcal{H}^{a,b}$ induced by $\nabla$. 

\begin{defn}
\label{Udef}
Denote by $\mathcal{U}$ the submodule of $\bigoplus_{a,b \geq 0} \mathcal{H}^{a,b}$ whose $(a,b)$ component is given by $\mathcal{F}^{a,b} \cap F^{\textrm{mid}}$. Note that $\mathcal{U} \cap \mathcal{H}^{a,b}$ is a coherent sheaf for all $a, b \geq 0$. 
\end{defn}

\begin{lem}
\label{etaincl}
If $s \in S(K)$ is a point with $K$ a field of characteristic not dividing $r!$, and $\iota : \mathcal{H}_{s} \xrightarrow{\sim} K^{m}$ is any isomorphism, then $\eta^{d}_{r}$ induces a natural inclusion 
\[ \eta^{d}_{r}((J^{d}_{r} S)_{s}(K)) \hspace{0.5em} \subset \hspace{0.5em} \GL_{m}(K) \big\backslash (J^{d}_{r} O(\iota(F^{\bullet}_{s}), \iota(\mathcal{U}_{s})))(K) . \]
\end{lem}

\begin{proof}
Applying \autoref{existslem}, we may find an $r$-limp $\psi = q \circ M$ which is defined at $s$ and which is framed with respect to $\iota$: we recall that this means $M$ is a change-of-basis matrix between a filtration-compatible frame $v^{1}, \hdots, v^{m}$ and a flat frame $b^{1}, \hdots, b^{m}$ for the module $\mathcal{M}$ of formal order-$r$ sections of $\mathcal{H}$ over the ring $\mathcal{O}^{r}_{S, s}$ of order-$r$ germs of functions at $s$, and that moreover, the map $\iota$ is given by the basis $b^{1}_{s}, \hdots, b^{m}_{s}$. 

We observe that $q \circ M$ factors through $O(\iota(F^{\bullet}_{s}), \iota(\mathcal{U}_{s})) \subset \ch{L}$: when the definitions are unravelled, this is simply the statement that when restricted to $\mathcal{M}$, the flag $F^{\textrm{mid}}$ contains $\mathcal{U}$. Thus $(q \circ M) \circ j$ lies inside $(J^{d}_{r} O(\iota(F^{\bullet}_{s}), \iota(\mathcal{U}_{s})))(K)$. Passing to $\GL_{m}(K)$-quotients and using \autoref{bigjetconstr}(iv) we obtain the result. 
\end{proof}

\begin{defn}
Given a variation of Hodge structure $\mathbb{V}$ on a complex algebraic variety $Z$, we denote by $\mathbf{H}_{Z}$ the algebraic monodromy group, which is the identity component of the Zariski closure of the monodromy representation associated to $\restr{\mathbb{V}}{Z^{\textrm{nor}}}$, where $Z^{\textrm{nor}} \to Z$ is the normalization.
\end{defn}

\begin{lem}
\label{orbiteqlem}
Suppose that $R = \mathbb{C}$, and the triple $(\mathcal{H}, F^{\bullet}, \nabla)$ underlies an integral variation of Hodge structure $\mathbb{V}$ with polarization $\mathcal{Q} : \mathbb{V} \otimes \mathbb{V} \to \mathbb{Z}$. Fix a point $s \in S(\mathbb{C})$ and an isomorphism $\iota : \mathcal{H}_{s} \xrightarrow{\sim} \mathbb{C}^{m}$. Then we have \[ O(\iota(F^{\bullet}_{s}), \iota(\mathcal{U}_{s})) = \iota( \mathbf{H}_{S}(\mathbb{C}) \cdot F^{\bullet}_{s}) ,\]
where we regard $\mathbf{H}_{S}$ as a subgroup of $\textrm{Aut}(\mathbb{V}_{s})$. Moreover if we let $F^{\bullet}_{c}$ denote the complex conjugate flag to $F^{\bullet}_{s}$, with complex conjugation with respect to $\mathbb{V}_{\mathbb{R}, s} \subset \mathcal{H}_{s}$, then there exists an integer $N$ depending only on $m$ and a subspace
\[ \mathcal{U}^{\textrm{mid}}_{s} \subset \mathcal{U}_{s} \cap F^{\textrm{mid}}_{s} \cap F^{\textrm{mid}}_{c} \cap \bigoplus_{a, b \leq N} (\mathbb{V}_{\mathbb{R},s})^{a,b} \] 
which satisfies the following properties
\begin{itemize}
\item[(a)] there exists an element $g \in \GL(\mathbb{V}_{s}, \mathcal{Q}_{s})(\mathbb{R})$ which normalizes $\mathbf{H}_{S}$ such that $g \mathcal{U}^{\textrm{mid}}_{s}$ is defined over $\mathbb{Q}$;
\item[(b)] the group $\mathbf{M}$ which stabilizes $g \mathcal{U}^{\textrm{mid}}_{s}$ is a Mumford-Tate group for a Hodge structure on $\mathbb{V}_{s}$ polarized by $\mathcal{Q}_{s}$, and $\mathbf{H}_{S}$ is its derived subgroup;
\item[(c)] we have an equality of orbits
\[ O(\iota(F^{\bullet}_{s}), \iota(\mathcal{U}_{s})) = O(\iota(F^{\bullet}_{s}), \iota(\mathcal{U}^{\textrm{mid}}_{s,\mathbb{C}})) , \]
\item[(d)] we have $\mathcal{Q}_{s} \in \mathcal{U}^{\textrm{mid}}_{s}$ and the adjoint Hodge numbers for the triple $(\iota(F^{\bullet}_{s}), \iota(F^{\bullet}_{c}), \iota(\mathcal{U}^{\textrm{mid}}_{s,\mathbb{C}}))$ are symmetric.
\end{itemize}
\end{lem}

\begin{proof}
We begin by noting that the choice of isomorphism $\iota$ plays no role in the statement, but is instead only present because our notation $O(-,-)$ refers specifically to subvarieties of the variety of Hodge flags on $\mathbb{C}^{m}$. Thus it suffices to show the statement for one choice of $\iota$, and in particular we may make the identification $\mathbb{V}_{s} = \mathbb{Z}^{m}$. From $\mathcal{Q}_{s}$ we thus obtain a bilinear form $Q : \mathbb{Z}^{m} \otimes \mathbb{Z}^{m} \to \mathbb{Z}$.

Let $D$ be the space of all Hodge structures on $\mathbb{Z}^{m}$ polarized by $Q$; this is naturally a homogenous space under $\mathbf{G}(\mathbb{R})$, where $\mathbf{G} = \textrm{Aut}(\mathbb{Z}^{m}, Q)$. The space $D$ sits as an open submanifold in a flag variety $\ch{D} \subset \ch{L}$ on which $\mathbf{G}(\mathbb{C})$ acts transitively. Denote by $\mathbf{G}_{S}$ the generic Mumford-Tate group of $S$, i.e., the group stabilizing those rational tensors inside $F^{\textrm{mid}} \cap \bigoplus_{a,b \geq 0} \mathbb{V}_{\mathbb{Q}}^{a,b}(S)$. By a result of Andr\'e-Deligne, the group $\mathbf{H}_{S}$ is a semisimple normal subgroup of the reductive group $\mathbf{G}_{S}$. Let $h_{s} \in D$ be the Hodge structure induced by $F^{\bullet}_{s}$. As in \cite[III.A]{GGK}, the Mumford-Tate group $\mathbf{G}_{S}$ splits as an almost direct product $\mathbf{G}_{S} = \mathbf{H}_{S} \cdot \mathbf{H}_{2} \cdot \mathbf{T}$, where $\mathbf{H}_{2}$ is semisimple and $\mathbf{T}$ is a toral factor. The orbit $D_{S} = \mathbf{G}_{S}(\mathbb{R}) \cdot h_{s}$ splits in a corresponding fashion as a product $D_{1} \times D_{2}$, and the Mumford-Tate group $\mathbf{M}$ of a generic Hodge structure $h$ lying inside an orbit $\mathbf{H}_{S}(\mathbb{R}) \cdot h_{CM}$, where $h_{CM} \in D_{S}$ is a CM point, has $\mathbf{H}_{S}$ as its derived subgroup. From the fact that $\mathbf{H}_{S}$ is the derived subgroup of $\mathbf{M}$, one has from \cite[III.A.1]{GGK} that the real orbits $\mathbf{H}_{S}(\mathbb{R}) \cdot h$ and $\mathbf{M}(\mathbb{R}) \cdot h$ agree. As $\mathbf{M}(\mathbb{R}) \cdot h$ is open inside $\mathbf{M}(\mathbb{C}) \cdot h$ by \cite[II pg.56]{GGK}, this implies that the complex orbits $\mathbf{H}_{S}(\mathbb{C}) \cdot h$ and $\mathbf{M}(\mathbb{C}) \cdot h$ also agree.

The first claim is now demonstrated as follows. Let $\textrm{Hg}(\mathbf{M})$ denote the subspace of $\bigoplus_{a, b \geq 0} (\mathbb{Q}^{m})^{a,b}$ which consists of Hodge tensors for $h$. Choose $g \in G_{S}(\mathbb{R})$ such that $g h_{s} = h$. Note that the Theorem of the Fixed Part \cite{Steenbrink1985} implies that $\mathcal{U}_{s}$ is exactly the subset of $\bigoplus_{a, b \geq 0} (\mathbb{C}^{m})^{a,b}$ which is both stabilized by $\mathbf{H}_{S,\mathbb{C}}$ and lies inside $F^{\textrm{mid}}_{s}$. From the fact that $\mathbf{H}_{S}$ is normal in $\mathbf{G}_{S}$ it then follows that  $\textrm{Hg}(\mathbf{M})_{\mathbb{C}} \subset g \mathcal{U}_{s}$, and hence
\begin{align*}
 \mathbf{H}_{S}(\mathbb{C}) \cdot h_{s} &\subset O(F^{\bullet}_{s}, \mathcal{U}_{s}) \\ 
 &= g^{-1} O(g F^{\bullet}_{s}, g \mathcal{U}_{s}) \\
 &\subset g^{-1} O(g F^{\bullet}_{s}, \textrm{Hg}(\mathbf{M})_{\mathbb{C}}) \\
 &= g^{-1} \mathbf{M}(\mathbb{C}) \cdot h &(*) \\
 &= g^{-1} \mathbf{H}_{S}(\mathbb{C}) \cdot h \\
 &= g^{-1} \mathbf{H}_{S}(\mathbb{C}) g \cdot h_{s} \\
 &= \mathbf{H}_{S}(\mathbb{C}) \cdot h_{s} .
\end{align*}
where on the line labelled $(*)$ we have applied \cite[II.C.1]{GGK}. Note that it need not be true that $\mathbf{M}$ is normal in $\mathbf{G}_{S}$; we only use that $\mathbf{H}_{S}$ is. 

For the second claim, we define $\mathcal{U}^{\textrm{mid}}_{s} = g^{-1} \textrm{Hg}(\mathbf{M})_{\mathbb{R}} \cap \bigoplus_{a,b \leq N} (\mathbb{R}^{m})^{a,b}$, with $N$ chosen large enough so that $g^{-1} \mathbf{M} g$ is the stabilizer of $g^{-1} \textrm{Hg}(\mathbf{M})_{\mathbb{R}} \cap \bigoplus_{a,b \leq N} (\mathbb{R}^{m})^{a,b}$, and observe that $g^{-1} \textrm{Hg}(\mathbf{M})$ is a subset of $\mathcal{U}_{s}$ that lies inside $F^{\textrm{mid}} \cap F^{\textrm{mid}}_{c}$ (this uses that $g$ is a real point!). Parts (a) and (b) are immediate. For part (c) use the above calculation with \autoref{ABCProp} to get
\begin{align*}
O(F^{\bullet}_{s}, \mathcal{U}_{s}) = \mathbf{H}_{S}(\mathbb{C}) \cdot h_{s} = (g^{-1} \mathbf{M}(\mathbb{C}) g) \cdot h_{s} = O(F^{\bullet}_{s}, g^{-1} \textrm{Hg}(\mathbf{M})_{\mathbb{C}}) . 
\end{align*}
For (d) one uses the fact that the adjoint Hodge numbers of a Mumford-Tate Lie algebra are symmetric, as one obtains by complex conjugating. Finally, to see that $N$ can be taken to only depend on $m$ one then uses the fact that $\mathbf{M}$ is a Mumford-Tate group, that $\mathbf{M}_{\mathbb{C}}$ stabilizes exactly the subspace $\textrm{Hg}(\mathbf{M})_{\mathbb{C}} \subset \bigoplus_{a, b \geq 0} (\mathbb{C}^{m})^{a,b}$, and that Mumford-Tate groups corresponding to points in $D$ lie in only finitely many $\GL_{m}(\mathbb{C})$ conjugacy classes by \cite[Theorem 4.14]{voisin2010hodge}. 
\end{proof}

\section{Jets and Exceptional Intersections}

\subsection{Proper Containment of Exceptional Loci}
\label{propexplocisec}

In this section we consider triples $(\mathcal{H}, F^{\bullet}, \nabla)$ as in the statement of \autoref{bigjetconstr}. We make the additional assumption that $R$ is a subring of $\mathbb{C}$, and that the data $(\mathcal{H}_{\mathbb{C}}, F^{\bullet}_{\mathbb{C}}, \nabla_{\mathbb{C}})$ obtained by base-change underlies a polarizable integral variation of Hodge structure $\mathbb{V}$.

\begin{notn}
Given a ring $R$ and an integer $\ell$, we will write $R_{\ell} = R[1/\ell!]$. 
\end{notn}

\begin{defn}
\label{framedlocdef}
An analytic map $\psi : B \to \an{\ch{L}_{\mathbb{C}}}$ with $B$ an open ball in $\an{S_{\mathbb{C}}}$ is called a \emph{local period map} if it is of the form $\psi = \an{q_{\mathbb{C}}} \circ M$, where $M : B \to \an{\GL_{\mathbb{C}}}$ is an analytic function giving the change-of-basis matrix between an analytic filtration-compatible frame on $B$ and a flat frame on $B$. Given a point $s \in B$, we define a frame $\iota : \mathbb{V}_{\mathbb{C}, s} \xrightarrow{\sim} \mathbb{C}^{m}$ at $s$ for $\psi$ analogously to \autoref{framedlimpdef}. 
\end{defn}

\begin{defn}
Given a point $s \in S(\mathbb{C})$ and map $\iota : \mathbb{V}_{\mathbb{C}, s} \xrightarrow{\sim} \mathbb{C}^{m}$, we will denote by $\ch{D}_{\iota}$ the flag subvariety on $\ch{L}_{\mathbb{C}}$ obtained by transferring the orbit $\mathbf{H}_{S,s}(\mathbb{C}) \cdot F^{\bullet}_{s}$ along $\iota$.
\end{defn}

\begin{lem}
\label{Diindeplem}
The dimension of $\ch{D}_{\iota}$ is independent of the point $s$ and the choice of isomorphism $\iota$.
\end{lem}

\begin{proof}
The independence of $\iota$ is obvious. To show that the dimension of $\mathbf{H}_{S, s}(\mathbb{C}) \cdot F^{\bullet}_{s}$ is independent of $S$ it suffices to produce, for each $s, s' \in S(\mathbb{C})$, an isomorphism $i : \mathbb{V}_{\mathbb{C}, s} \xrightarrow{\sim} \mathbb{V}_{\mathbb{C}, s'}$ taking $\mathbf{H}_{S,s}(\mathbb{C})$ to $\mathbf{H}_{S,s'}(\mathbb{C})$ and $F^{\bullet}_{s}$ to $F^{\bullet}_{s'}$. By parallel translation from $s$ to $s'$ we obtain $i : \mathbb{V}_{s} \xrightarrow{\sim} \mathbb{V}_{s'}$ which conjugates $\mathbf{H}_{S,s}$ and $\mathbf{H}_{S,s'}$. It then suffices to show that $i(F^{\bullet}_{s})$ and $F^{\bullet}_{s'}$ both lie in the same $\mathbf{H}_{S,s'}(\mathbb{C})$-orbit. This follows from \cite[III.A]{GGK}, where a global period map is constructed on the universal cover $\widetilde{S} \to S$ and it is argued that the resulting Hodge flags map to a single algebraic monodromy orbit in the period domain. 
\end{proof}

\begin{defn}
\label{perdimdef}
We will denote the quantity $\dim \ch{D}_{\iota}$ of \autoref{Diindeplem} by $\overline{P}$.
\end{defn}

We also write $\mathcal{P}^{d}_{r}$ for the torsor constructed in \S\ref{glusec} which, together with the maps $(\gamma, \alpha)$ constructed in \S\ref{mainjetconstrsec}, defines the map $\eta^{d}_{r}$ of \autoref{bigjetconstr}. Although the pair $(\gamma, \alpha)$ also depends on both $d$ and $r$, in light of the compatibility \autoref{bigjetconstr}(iii) described in \S\ref{keypropsec} we will omit this reference as it will cause no ambiguity. We write $\Phi^{d}_{r} \subset \mathcal{P}^{d}_{r} \times J^{d}_{r} \ch{L}_{R}$ for the graph of $\alpha$.

\begin{defn}
We say an irreducible complex subvariety $Z \subset S_{\mathbb{C}}$ is \emph{weakly special} if it is maximal for its algebraic monodromy group. More precisely, given any irreducible complex subvariety $Y \subset S_{\mathbb{C}}$ which contains $Z$ the natural inclusion $\mathbf{H}_{Z} \subset \mathbf{H}_{Y}$ is an equality only if $Y = Z$.
\end{defn}

\begin{defn}
We say a weakly special subvariety $Z \subset S_{\mathbb{C}}$ is \emph{rigid} if there does not exist a weakly special subvariety $Y \subset S_{\mathbb{C}}$ properly containing $Z$ such that $\mathbf{H}_{Z}$ is normal in $\mathbf{H}_{Y}$. 
\end{defn}

In the proposition that follows, we will require the \emph{Ax-Schanuel Theorem} for variations of Hodge structures due to Bakker-Tsimerman. To state the result, suppose that $\psi : B \to \ch{L}_{\mathbb{C}}$ is a local period map with frame $\iota$, and whose image is contained in $\ch{D}_{\iota} \subset \ch{L}_{\mathbb{C}}$. Let $T \subset S_{\mathbb{C}} \times \ch{D}_{\iota}$ be a closed algebraic subvariety, and let $\Gamma_{\psi} \subset B \times \ch{L}_{\mathbb{C}}$ be the graph of $\psi$. Then we have

\begin{thm}[\cite{AXSCHAN}]
\label{axschanthm}
Suppose that $U \subset \Gamma_{\psi} \cap T$ is an analytic component such that
\begin{equation}
\label{axshanineq}
\textrm{codim}_{S_{\mathbb{C}} \times \ch{D}_{\iota}} U < \textrm{codim}_{S_{\mathbb{C}} \times \ch{D}_{\iota}} T + \textrm{codim}_{S_{\mathbb{C}} \times \ch{D}_{\iota}} \Gamma_{\psi} .
\end{equation}
Then the projection of $U$ to $B$ lies in a weakly special subvariety properly contained in $S_{\mathbb{C}}$.
\end{thm}

\begin{prop}
\label{expintprop}
Suppose that the algebraic monodromy group $\mathbf{H}_{S}$ associated to the variation $\mathbb{V}$ underlying $(\mathcal{H}_{\mathbb{C}}, F^{\bullet}_{\mathbb{C}}, \nabla_{\mathbb{C}})$ is $\mathbb{Q}$-simple. Suppose that $g : \mathcal{Y} \to \mathcal{M}$ is an $R$-algebraic family of subschemes of $\ch{L}_{R}$ over the Noetherian base $\mathcal{M}$, with projection $\pi : \mathcal{Y} \to \ch{L}_{R}$ restricting to an embedding on fibres, and with fibres of dimension $e$. Suppose that for some integer $d$ we have
\[ \overline{P} - e > \dim_{R} S - d , \]
where $\dim_{R}$ denotes the relative dimension over $R$. Let 
\[ \mathcal{E}_{r} \subset \Phi^{d}_{r} \cap (\gamma^{-1}(J^{d}_{r,nd} S) \times (\pi^{r}_{1})^{-1}(J^{d}_{1,nc} \ch{L}_{R})) \subset \mathcal{P}^d_r \times J^d_r \ch{L}_{R} , \]
with $\pi^r_1 : J^d_r \ch{L}_{R} \to J^d_1 \ch{L}_{R}$ the natural projection, be the locus of points whose image inside $J^{d}_{r} \ch{L}_{R}$ lie inside $J^{d}_{r} \mathcal{Y}_{m}$ for some $m$, where the map $m : \kappa(\mathfrak{q}) \to \mathcal{M}$ is induced by a scheme-theoretic point $\mathfrak{q} \in \mathcal{M}$. Then there exists $r$ such that, after replacing $R$ with $R_{r}$ and applying base-change, the image in $S$ of $\mathcal{E}_{r}$ lies in a closed subscheme $E$ which is properly contained in $S$.
\end{prop}

\begin{proof}
Let us denote by $E_{r}$ the image of $\mathcal{E}_{r}$ in $S$, which we obtain after base-changing to $R_{r}$ if necessary. The set $E_{r}$ is a constructible subset of $S$: it is the image in $S$ of the (constructible, by Chevalley's Theorem) intersection 
\[ \mathcal{E}_{r} := \Phi^{d}_{r} \cap (\gamma^{-1}(J^{d}_{r,nd} S) \times (\pi^{r}_{1})^{-1}(J^{d}_{1,nc} \ch{L}_{R})) \cap \textrm{im}[\mathcal{P}^{d}_{r} \times \mathcal{Y}^{d}_{r} \xrightarrow{\textrm{id} \times (J^{d}_{r} \pi)} \mathcal{P}^{d}_{r} \times J^{d}_{r} \ch{L}_{R}] , \]
where the subscheme $\mathcal{Y}^{d}_{r} \subset J^{d}_{r} \mathcal{Y}$ is as in \autoref{jetfams}. From the construction one has $E_{r+1} \subset E_{r}$ for all $r$, where the inclusion takes place inside $S_{R_{r+1}}$. We note that the formation of $\mathcal{E}_{r}$ is compatible with base-change along maps $R \to R'$, which follows from the compatibility of the maps of \autoref{bigjetconstr} with base-change. 

We first show that such an $r$ exists in the case where $R = \mathbb{C}$. We assume for the purpose of contradiction that each $E_{r}$ does not lie in a proper closed subvariety of $S$, hence that $E_{\infty} := \bigcap_{r} E_{r}$ contains the complement of a countable union $\bigcup_{r} D_{r}$ of proper closed complex subvarieties $D_{r} \subset S$. We will obtain the desired contradiction by showing that $E_{\infty}$ lies in a countable union $\bigcup_{i} Z_{i}$ of proper closed subvarieties of $S$, proving that the set $S(\mathbb{C})$ with its analytic topology is a union of countably many meagre subsets --- an absurdity. 

The countable collection $\{ Z_{i} \}_{i \in I}$ of proper closed subvarieties of $S$ we refer to here will be the collection of \emph{rigid weakly special} subvarieties for the variation $\mathbb{V}$ which are not equal to $S$. The set of rigid weakly special subvarieties is a countable set (see for instance \cite[\S2]{fieldsofdef}, where these subvarieties are called ``weakly non-factor''). The hypothesis that $\mathbf{H}_{S}$ be $\mathbb{Q}$-simple in fact implies that any weakly special $Z$ for which $\mathbf{H}_{Z}$ is neither $\mathbf{H}_{S}$ nor the trivial group lies in a rigid weakly special; see the proof of \cite[Cor. 1.13]{fieldsofdef}.

Let $s \in E_{\infty}$ be a point. To apply the Ax-Schanuel theorem to show that $s$ lies in one of the $Z_{i}$, we proceed as follows. First, let us denote by $\mathcal{E}_{r,s}$ the fibre of $\mathcal{E}_{r}$ above $s$. The sets $\mathcal{E}_{r,s}$ form a family compatible with projection as in \autoref{existscompseqlem}, so we obtain a compatible sequence $\{ (t_{r}, \sigma_{r}) \}_{r \geq 0}$ where $(t_{r}, \sigma_{r}) \in \mathcal{E}_{r,s}$ for all $r$. Using the interpretation of $\mathcal{P}^{d}_{r}$ in \autoref{Pinterp} and the natural identification $\mathcal{H}_{s} \simeq \mathbb{V}_{\mathbb{C}, s}$ the point $t_{0}$ defines a frame $\iota : \mathbb{V}_{\mathbb{C}, s} \xrightarrow{\sim} \mathbb{C}^{m}$, and if $\psi$ is a local period map for which $\iota$ is a frame then applying \autoref{Pinterp} we have that $\sigma_{r} = \alpha(t_{r}) = \psi \circ (\gamma(t_{r}))$ for all $r$. Moreover, we have by construction that $j_{r} = \gamma(t_{r})$ is a compatible sequence of $d$-dimensional non-degenerate jets, and we denote by $B_{d} \subset B$ the irreducible $d$-dimensional closed analytic locus through which all the $j_{r}$ factor.

Now consider the compatible sequence $\{ \sigma_{r} \}_{r \geq 0}$, and denote by $\mathcal{M}_{r} \subset \mathcal{M}$ the constructible locus of points $m$ such that $\sigma_{r} \in (J^{d}_{r} \pi) (\mathcal{Y}^{d}_{r,m})$. Then by construction one has that $\mathcal{M}_{r+1} \subset \mathcal{M}_{r}$ for all $r$ and each $\mathcal{M}_{r}$ is non-empty; from constructibility and Noetherianity of $\mathcal{M}$ there exists a point $m \in \bigcap_{r \geq 0} \mathcal{M}_{r}$, and we have $\sigma_{r} \in J^{d}_{r} ( \pi(\mathcal{Y}_{m}))$ for all $r$. Applying \autoref{jetfactorlem} to $(j_{r}, \restr{\psi}{B_{d}})$ we note that $\psi(B_{d}) \subset \pi(\mathcal{Y}_{m})$. 

Let $\Gamma_{\psi} \subset B \times \ch{D}_{\iota}$ be the graph of $\psi$. Setting $T = S \times \pi(\mathcal{Y}_{m})$ and letting $U \subset \Gamma_{\psi} \cap T$ the component containing the graph of $\restr{\psi}{B_{d}}$ we have that
\begin{align*}
\codim_{S \times \ch{D}_{\iota}} U \leq (\dim S - d) + \dim \ch{D}_{\iota} &< (\dim \ch{D}_{\iota} - \pi(\mathcal{Y}_{m})) + \dim \ch{D}_{\iota} \\ 
&= \textrm{codim}_{S \times \ch{D}_{\iota}} T + \textrm{codim}_{S \times \ch{D}_{\iota}} \Gamma_{\psi} ,
\end{align*}
where we have used the hypothesis of the proposition to obtain the inequality. It then follows from the Ax-Schanuel Theorem that $B_{d}$ lies inside a weakly special subvariety $Z$ which is properly contained in $S$. From \autoref{nottrivgroup} below it follows that $\mathbf{H}_{Z}$ is not the trivial group, hence by our setup $Z$ lies in one of the $Z_{i}$. Thus $s$, which lies in $Z$, also lies in one of the $Z_{i}$.

\vspace{0.5em}

Having shown the statement in the case $R = \mathbb{C}$, we return now the case of a general $R$. The above argument has shown that for large enough $r$, the projection of $\mathcal{E}_{r,\mathbb{C}}$ to $S_{\mathbb{C}}$ lies in a proper closed subscheme. Now replace $R$ with $R_{r}$, and $S$ with its base-change to $R_{r}$. Denoting by $E_{r}$ the projection of $\mathcal{E}_{r}$ to $S$, we can write $E_{r}$ as a union $\bigcup_{i = 1}^{n} U_{i} \cap V_{i}$ with $\{ U_{i} \}_{i = 1}^{n}$ (resp. $\{ V_{i} \}_{i = 1}^{n}$) a sequence of open (resp. closed) subschemes of $S$. It suffices to show that none of the $V_{i}$ is equal to $S$. Such a $V_{i}$ must lie over the generic fibre of $\Spec R$, hence $V_{i,\mathbb{C}} = S_{\mathbb{C}}$. But then $U_{i, \mathbb{C}} \cap V_{i,\mathbb{C}}$ lies in $E_{r,\mathbb{C}}$, which we showed was not Zariski dense in $S_{\mathbb{C}}$. 
\end{proof}

\begin{lem}
\label{nottrivgroup}
Suppose that $Z$ is a variety with a variation of Hodge structure $\mathbb{V}$, that $\psi : B \to \ch{L}$ is a local period map with $B$ lying in the smooth locus of $Z$, and that for some $d, r$ there exists $j \in J^{d}_{r} B$ such that $\psi \circ j$ is not constant. Then $\mathbf{H}_{Z}$ is not the trivial group.
\end{lem}

\begin{proof}
We suppose that $\mathbf{H}_{Z}$ is trivial. We will show that at each point $s$ in the smooth locus of $Z$ there exists a local period map $\psi : B \to \ch{L}$ with $s \in B$ such that $\psi$ is constant. This will show the result, since for each $s \in S(\mathbb{C})$ the germs of local period maps $\psi$ defined at $s$ lie in a single $\GL_{m}(\mathbb{C})$ orbit, i.e., if $\psi = q \circ M$ is a local period map, then so is $q \circ (C \circ M)$ for $C \in \GL_{m}(\mathbb{C})$, and all local period maps defined at $s$ are obtained this way.

The statement in question is unchanged after replacing $Z$ with its smooth locus and passing to a finite \'etale covering of $Z' \to Z$. Since $\mathbf{H}_{Z}$ is the identity component of the Zariski closure of the monodromy representation $\pi_{1}(Z, s) \to \textrm{Aut}(\mathbb{V}_{s})$, we may therefore assume after passing to such a covering that $\mathbb{V}$ is constant. The Theorem of the Fixed Part \cite{Steenbrink1985} then tells us that the space of global sections $\mathbb{V}(Z)$ admits a Hodge structure inducing the Hodge structure on each fibre of $\mathbb{V}$. As a consequence, we can choose at any point $s \in S(\mathbb{C})$ a flat frame which is \emph{also} a filtration-compatible frame. Constructing the local period map $\psi = q \circ M$ with respect to this frame the matrix $M$ is the identity, hence $\psi$ is constant. 
\end{proof}

\subsection{A Criterion for Exceptionality}

In this section we assume that $S$ is a smooth complex algebraic variety with a polarized variation of Hodge structure $\mathbb{V}$. 

\begin{notn}
If $\iota : \mathbb{V}_{\mathbb{C}, s} \xrightarrow{\sim} \mathbb{C}^{m}$ is a frame and $W$ is a subspace of $\mathbb{V}^{a,b}_{\mathbb{C},s}$ for some $a, b \geq 0$, we write $\iota(W)$ for the subspace of $(\mathbb{C}^{m})^{a,b}$ induced by $\iota$. 
\end{notn}

\begin{defn}
\label{varleveldef}
For $k \geq 0$ an integer, we say that the variation of Hodge structure $\mathbb{V}$ on $S$ has level at least $k$ if there exists a point $s \in S(\mathbb{C})$ such that the adjoint Hodge structure on the Lie algebra of the algebraic monodromy group $\mathbf{H}_{S,s}$ at $s$ has level at least $k$.
\end{defn}

The remainder of the section is devoted to the proof of \autoref{noatyplem}, which gives a criterion for the inverse image of an algebraic group orbit under a period map to have ``exceptional'' dimension. The arguments we use roughly resemble those of \cite[\S6]{BKU}, but have been modified in some respects to take into account the fact that our conjugate filtrations need not come from ``complex conjugation'' with respect to an underlying real structure.

\begin{prop}
\label{noatyplem}
Suppose that $\mathbb{V}$ has level at least three, and that $\mathbf{H}_{S}$ is absolutely simple (i.e., $\mathbf{H}_{S,\mathbb{C}}$ is simple). Suppose that $\psi : B \to \an{\ch{L}}_{\mathbb{C}}$ is a local period map with $s \in B$ and frame $\iota : \mathbb{V}_{\mathbb{C}, s} \xrightarrow{\sim} \mathbb{C}^{m}$, and consider the flag $\iota(F^{\bullet}_{s})$ transferred to $\mathbb{C}^{m}$ using $\iota$. Let $F^{\bullet}_{c}$ be a filtration opposed to $\iota(F^{\bullet}_{s})$ on $\mathbb{C}^{m}$, and let $\mathbf{U} \subset \iota(\mathbf{H}_{S,\mathbb{C}})$ be a proper complex algebraic subgroup with Lie algebra $\mathbf{u}$, such that:
\begin{itemize}
\item[-] the summands $H^{p,q} = \iota(F^{p}_{s}) \cap F^{q}_{c}$ give a direct sum decomposition of $\mathbb{C}^{m}$; 
\item[-] the induced direct sum decomposition on $\textrm{End}(\mathbb{C}^{m})$ given by (\ref{endiieq}) restricts to direct sum decompositions $\iota(\mathfrak{h}_{S}) = \bigoplus_{i} \mathfrak{h}^{i}_{S}$ and $\mathfrak{u} = \bigoplus_{i} \mathfrak{u}^{i}$, and both decompositions have symmetric Hodge numbers.
\end{itemize}
Then
\begin{equation}
\label{expcodimconc}
\dim (\mathbf{H}_{S}(\mathbb{C}) \cdot \iota(F^{\bullet}_{s})) - \dim (\mathbf{U}(\mathbb{C}) \cdot \iota(F^{\bullet}_{s})) > \dim S - \dim \psi^{-1}(\mathbf{U}(\mathbb{C}) \cdot \iota(F^{\bullet}_{s})) .
\end{equation}
\end{prop}


\vspace{0.5em}

\noindent \textbf{Reduction to generic $s$:} We begin by noting that it suffices to consider $s$ for which we have both 
\begin{align*}
\dim T_{s} \psi^{-1}(\mathbf{U}(\mathbb{C}) \cdot \iota(F^{\bullet}_{s})) &= \dim_{s} \psi^{-1}(\mathbf{U}(\mathbb{C}) \cdot \iota(F^{\bullet}_{s})) \\ 
T_{s} \psi^{-1}(\mathbf{U}(\mathbb{C}) \cdot \iota(F^{\bullet}_{s})) &= (d\psi)^{-1}_{s} (T_{\psi(s)} (\mathbf{U}(\mathbb{C}) \cdot \iota(F^{\bullet}_{s}))) .
\end{align*} 
The analytic locus where these two conditions hold is dense in $\psi^{-1}(\mathbf{U}(\mathbb{C}) \cdot \iota(F^{\bullet}_{s}))$, so let us choose a point $s'$ in this locus and show that the hypotheses of the theorem continue to hold after replacing $s$ with $s'$. The local period map $\psi$ admits a decomposition $\psi = q \circ M$, where $M$ is the change-of-basis matrix between a filtration-compatible frame $v^{1}, \hdots, v^{m}$ on $B$ and a flat frame $b^{1}, \hdots, b^{m}$ such that the isomorphism $\iota$ is given by $b^{1}_{s}, \hdots, b^{m}_{s}$. The map $\iota$ is induced by a unique map $\mathbb{V}_{\mathbb{C}}(B) \xrightarrow{\sim} \mathbb{C}^{m}$, which in turn induces a unique map $\iota' : \mathbb{V}_{\mathbb{C}, s'} \to \mathbb{C}^{m}$. Then by construction, $\iota'(F^{\bullet}_{s'})$ lies inside $\mathbf{U}(\mathbb{C}) \cdot \iota(F^{\bullet}_{s})$. Write $\iota'(F^{\bullet}_{s'}) = g \cdot \iota(F^{\bullet}_{s})$ for some $g \in \mathbf{U}(\mathbb{C})$ and define $F'^{\bullet}_{c} = g \cdot F^{\bullet}_{c}$. The required properties for the pair $(\iota'(F^{\bullet}_{s'}),  F'^{\bullet}_{c})$ are then easily checked from the construction. 

\vspace{0.5em}

\noindent \textbf{Constraint from transversality: }From now on we drop reference to $\iota$ and identify $\mathbb{V}_{\mathbb{C}, s}$ and $\mathbb{C}^{m}$ for ease of notation. Note that as $\psi(B) \subset \mathbf{H}_{S}(\mathbb{C}) \cdot F^{\bullet}_{s}$ (for instance, by applying $\psi(B) \subset O(F^{\bullet}_{s}, \mathcal{U}_{s})$ and \autoref{orbiteqlem}), we have on general grounds that 
\begin{equation}
\label{codimineq}
\dim (\mathbf{H}_{S}(\mathbb{C}) \cdot F^{\bullet}_{s}) - \dim (\mathbf{U}(\mathbb{C}) \cdot \iota(F^{\bullet}_{s})) \geq \dim S - \dim \psi^{-1}(\mathbf{U}(\mathbb{C}) \cdot \iota(F^{\bullet}_{s})) ,
\end{equation}
so it suffices to rule out the possibility of equality. By differentiating the natural orbit maps we obtain identifications
\begin{align*}
T_{F^{\bullet}_{s}} (\mathbf{H}_{S}(\mathbb{C}) \cdot F^{\bullet}_{s}) = \bigoplus_{i > 0} \mathfrak{h}^{i}_{S} \hspace{4em} T_{F^{\bullet}_{s}} (\mathbf{U}(\mathbb{C}) \cdot \iota(F^{\bullet}_{s})) = \bigoplus_{i > 0} \mathfrak{u}^{i} .
\end{align*}
where the summands indexed by $i$ are obtained by intersection with $\textrm{End}(\mathbb{C}^{m})^{i}$, defined as in (\ref{endiieq}). 

Using these identifications, Griffiths transversality, and our choice of $s$, we can reformulate the condition that (\ref{codimineq}) be an equality as
\begin{align*}
\dim \left( \bigoplus_{i > 0} \mathfrak{h}^{i}_{S} \right) -  \dim \left( \bigoplus_{i > 0} \mathfrak{u}^{i} \right) &= \dim T_{s} S - \dim (d\psi)_{s}^{-1}\left( \bigoplus_{i > 0} \mathfrak{u}^{i} \right) \\
&= \dim \, (d\psi)^{-1}_{s}\left( \mathfrak{h}^{1}_{S} \right) - \dim \, (d\psi)_{s}^{-1}\left( \mathfrak{u}^{1} \right) \\
&\leq \dim \mathfrak{h}^{1}_{S} - \dim \mathfrak{u}^{1} .
\end{align*}
This in particular implies that for $|i| \geq 2$ we have $\mathfrak{u}^{i} = \mathfrak{h}_{S}^{i}$. 

\vspace{0.5em}

\noindent \textbf{Properties of $\mathfrak{h}_{S}$:} Next we collect some facts about the direct sum decomposition of $\mathfrak{h}_{S}$ implied by semisimplicity. As the algebra $\mathfrak{h}_{S}$ is semisimple with a Lie bracket-compatible grading, a theorem of Zassenhaus \cite[Chap 3, Lemma 3.6.1]{Duistermaat2000} shows that there exists an element $x \in \mathfrak{h}_{S}$ such that $[x,y] = i y$ for any $y \in \mathfrak{h}^{i}_{S}$. Because the center of the semisimple Lie algebra $\mathfrak{h}_{S}$ is trivial, the element $x$ is unique, and lies in $\mathfrak{m} = \mathfrak{h}^{0}_{S}$. Since $\textrm{ad}\, x$ is diagonalized by a basis for $\mathfrak{h}_{S}$ compatible with the Hodge decomposition, $x$ is a semisimple element of $\mathfrak{h}_{S}$, hence lies in a Cartan subalgebra $\mathfrak{c}$. Since $[x,c] = 0$ for each $c \in \mathfrak{c}$, one necessarily has $\mathfrak{c} \subset \mathfrak{m}$. 

If one now considers a root space decomposition $\mathfrak{h}_{S} = \mathfrak{c} \oplus \bigoplus_{\mu \in R} \mathfrak{h}_{\mu}$ with respect to $\mathfrak{c}$, with $R$ the set of roots, then this decomposition necessarily refines the Hodge decomposition of $\mathfrak{h}_{S}$. From this one easily checks that both: 
\begin{itemize}
\item[(1)] the Killing form $B(y,z) = \textrm{tr}(\textrm{ad}y \circ \textrm{ad}z)$, which is non-degenerate on $\mathfrak{h}_{S}$, induces a perfect pairing between $\mathfrak{h}^{i}_{S}$ and $\mathfrak{h}^{-i}_{S}$ for each $i$; and
\item[(2)] the induced direct sum decomposition of $\mathfrak{h}^{0}_{S}$ is a root space decomposition, making $\mathfrak{h}^{0}_{S}$ into a reductive Lie subalgebra.
\end{itemize}

\vspace{0.5em}

\noindent \textbf{$\mathfrak{t}$-roots:} To carry out our argument we will need a slightly different root theory due to Kostant (we also roughly follow similar arguments in \cite[\S6]{BKU}); this is not to be confused with the traditional root theory we used above, even though we use the same notation. To begin with, we note that both $\mathfrak{q} = \bigoplus_{i \geq 0} \mathfrak{h}^{i}_{S}$ and $\bigoplus_{i \leq 0} \mathfrak{h}^{i}_{S}$ are parabolic subalgebras: this is immediate from the root space decomposition above. We call a non-zero element $\nu \in \mathfrak{t}^{*}$ a $\mathfrak{t}$-root if the space
\[ \mathfrak{h}_{\nu} = \{ z \in \mathfrak{h}_{S} : (\textrm{ad}\, y)(z) = \nu(y) z, \forall y \in \mathfrak{t} \} , \]
is non-zero. If $\mathfrak{r} = \bigoplus_{i \neq 0} \mathfrak{h}^{i}_{S}$ is the Killing form orthocomplement of $\mathfrak{m}$ in $\mathfrak{h}_{S}$, one has that
\begin{equation}
\label{rdecomp}
\mathfrak{r} = \bigoplus_{\nu \in R} \mathfrak{h}_{\nu} ,
\end{equation}
where $R \subset \mathfrak{t}^{*}$ is the set of $\mathfrak{t}$-roots. In \cite{Kostant2010}, Kostant proves that:

\begin{quote} For each $\nu \in R$, $\mathfrak{h}_{\nu}$ is an irreducible $\textrm{ad}\, \mathfrak{m}$ module, and any irreducible $\textrm{ad}\, \mathfrak{m}$-submodule of $\mathfrak{r}$ is of this form. Moreover, the decomposition (\ref{rdecomp}) is the unique decomposition of $\mathfrak{r}$ as a direct sum of irreducible $\textrm{ad}\, \mathfrak{m}$-modules.
\end{quote}

\noindent One thus has immediately that each space $\mathfrak{h}^{i}_{S}$ with $i \neq 0$ is uniquely a sum of spaces in $\{ \mathfrak{h}_{\nu} : \nu \in R \}$. One may additionally check that $\mathfrak{h}_{\nu} \subset \mathfrak{h}^{i}_{S}$ if and only if $\mathfrak{h}_{-\nu} \subset \mathfrak{h}^{-i}_{S}$ from the fact that $x \in \mathfrak{t}$. Finally, one has by \cite[Lemma 2.1]{Kostant2010} that $[\mathfrak{h}_{\nu}, \mathfrak{h}_{\mu}] = \mathfrak{h}_{\nu+\mu}$ for any two elements $\nu, \mu \in R$. 

We denote for future reference by $R_{\textrm{simp}} \subset R$ a subset of simple $\mathfrak{t}$-roots; that is a maximal set of $\mathfrak{t}$-roots such that each $\nu \in R_{\textrm{simp}}$ cannot be written as a sum $\nu = \nu_{1} + \nu_{2}$ for $\nu_{1}, \nu_{2} \in R$, and such that at most one element in $\{ \nu, -\nu \}$ lies in $R_{\textrm{simp}}$ for each $\nu \in R$. Without loss of generality we may assume that for each $\nu \in R_{\textrm{simp}}$ we have $\mathfrak{h}_{\nu} \subset \bigoplus_{i > 0} \mathfrak{h}^{i}_{S}$. Note that by \cite[Theorem 2.7]{Kostant2010} one has that $R_{\textrm{simp}}$ is a basis for $\mathfrak{t}^{*}$. 

\vspace{0.5em}

\noindent \textbf{Griffiths Transverse Generation:} Let us denote by $\mathcal{D} \subset T (\mathbf{H}_{S}(\mathbb{C}) \cdot F^{\bullet}_{s})$ the involutive distribution generated by the Griffiths transverse subbundle of $T (\mathbf{H}_{S}(\mathbb{C}) \cdot F^{\bullet}_{s})$ under $[-,-]$. As $\mathcal{D}$ is an involutive distribution we learn from the Frobenius theorem that in a small enough neighbourhood $\Delta \subset \mathbf{H}_{S}(\mathbb{C}) \cdot F^{\bullet}_{s}$ of $F^{\bullet}_{s}$ there exists a unique foliation of $\Delta$ by maximal integral manifolds for $\mathcal{D}$. Let show that the manifold $Y$ in this foliation passing through $F^{\bullet}_{s}$ is algebraic. To see this, consider the nilpotent Lie algebra $\mathfrak{n}$ generated by $\mathfrak{h}^{1}_{S}$. Because $\mathfrak{n}$ is Nilpotent, is the Lie algebra of a complex algebraic group $N$. As the natural map $\mathfrak{h}_{S} \to T_{F^{\bullet}_{s}} (\mathbf{H}_{S}(\mathbb{C}) \cdot F^{\bullet}_{s})$ preserves the Lie bracket and $\mathcal{D}$ is invariant under $N$, it follows that $Y = N \cdot F^{\bullet}_{s}$. From \cite[Lemma 4.10(ii)]{urbanik2021sets}, and the algebraicity of $Y$ we now conclude that $Y = \mathbf{H}_{S}(\mathbb{C}) \cdot F^{\bullet}_{s}$. It follows that $\mathfrak{h}^{1}_{S}$ generates $\bigoplus_{i > 0} \mathfrak{h}^{i}_{S}$. This in particular means that $\mathfrak{h}^{1}_{S} = \bigoplus_{\nu \in R_{\textrm{simp}}} \mathfrak{h}_{\nu}$. 

\vspace{0.5em}

\noindent \textbf{Showing $\mathfrak{u}^{1} = \mathfrak{h}^{1}_{S}$:} The essential consequence of the above analysis for us is that we may conclude that $\mathfrak{h}^{k}_{S} = [\mathfrak{h}^{k-1}_{S}, \mathfrak{h}^{1}_{S}]$, where $k$ is the level of $\mathbb{V}$. (We recall that $k \geq 3$ and $\mathfrak{h}^{k}_{S} \neq 0$ by assumption.)

Now consider the Lie subalgebra $\mathfrak{h} \subset \mathfrak{h}_{S}$ generated by $\bigoplus_{|i| \geq 2} \mathfrak{u}^{i} = \bigoplus_{|i| \geq 2} \mathfrak{h}^{i}_{S}$. This is naturally an $\textrm{ad}\, \mathfrak{m}$-submodule, and the set $R_{\mathfrak{h}} \subset R$ of $\mathfrak{t}$-roots $\nu$ for which $\mathfrak{h}_{\nu} \subset \mathfrak{h}$ is stable under the map $\nu \mapsto -\nu$. From the equality $[\mathfrak{h}^{k-1}_{S}, \mathfrak{h}^{1}_{S}] = \mathfrak{h}^{k}_{S}$ we find that $\mathfrak{h} \cap \mathfrak{h}_{S}^{1}$ must be non-empty: if we choose $\mu + \nu \in R$ such that $\mathfrak{h}_{\mu} \subset \mathfrak{h}^{1}_{S}$ and $\mathfrak{h}_{\nu} \subset \mathfrak{h}^{k-1}_{S}$, then $[\mathfrak{h}_{\mu + \nu}, \mathfrak{h}_{-\nu}] = \mathfrak{h}_{\mu} \subset \mathfrak{h}$. Let $\mathfrak{c}^{1}$ be the complementary $\textrm{ad}\, \mathfrak{m}$-submodule to $\mathfrak{h} \cap \mathfrak{h}^{1}_{S}$ inside $\mathfrak{h}^{1}_{S}$. Then for $\mathfrak{h}_{\gamma} \subset \mathfrak{c}^{1}$ and $\mathfrak{h}_{\mu} \subset \mathfrak{h} \cap \mathfrak{h}^{1}_{S}$ we must have $[\mathfrak{h}_{\mu}, \mathfrak{h}_{\gamma}] = 0$: otherwise $\mu + \gamma$ is a $\mathfrak{t}$-root and arguing as before one can prove $\mathfrak{h}_{\gamma} \subset \mathfrak{h}$. If $\mathfrak{h} \cap \mathfrak{h}^{1}_{S}$ is not all of $\mathfrak{h}^{1}_{S}$, we therefore obtain a non-trivial decomposition $R_{\textrm{simp}} = R_{1} \sqcup R_{2}$, with the $\mathfrak{t}$-roots in $R_{2}$ identified with the $\mathfrak{t}$-root spaces in $\mathfrak{h} \cap \mathfrak{h}^{1}_{S}$, and such that no sum of an element of $R_{1}$ and an element of $R_{2}$ is a $\mathfrak{t}$-root. Using \cite[(2.13)]{Kostant2010} and \cite[(2.33)]{Kostant2010} this implies that if $\nu_{1} \in R_{1}$ and $\nu_{2} \in R_{2}$ then $(\nu_{1}, \nu_{2}) = 0$, with $(-,-)$ the restriction of the Killing form as in \cite{Kostant2010}. 

We now show this contradicts \cite[Remark 2.9]{Kostant2010}, which says because $\mathfrak{h}_{S}$ is simple, the $\mathfrak{t}$-root obtained as the restriction to $\mathfrak{t}$ of the highest root of $\mathfrak{h}_{S}$ is of the form $\sum_{\beta \in R_{\textrm{simp}}} n_{\beta} \beta$ with all coefficients $n_{\beta} > 0$; we will show this is impossible if the $\mathfrak{t}$-root system is decomposable. Indeed, suppose that $\nu = \nu_{1} + \nu_{2}$ is a positive $\mathfrak{t}$-root, with $\nu_{1}$ a non-zero sum of elements of $R_{1}$ and $\nu_{2}$ a non-zero sum of elements of $R_{2}$. Using the definition of $R_{\textrm{simp}}$, we can find such an element so that either $\nu_{1} \in R_{\textrm{simp}}$ or $\nu_{2} \in R_{\textrm{simp}}$; assume the latter without loss of generality. Then $(\nu_{1}, \nu_{2}) = 0$ and $\nu = \nu_{1} + \nu_{2}$ is a $\mathfrak{t}$-root, so by \cite[(2.14)]{Kostant2010} so is $\nu_{1} - \nu_{2}$. Because $\nu_{2}$ is simple and $\nu_{1}$ is non-zero, $\nu_{1} - \nu_{2}$ must be a positive $\mathfrak{t}$-root, which is impossible because $R_{\textrm{simp}}$ forms a basis for $\mathfrak{t}^{*}$ and $\{ \nu_{2} \}$ is linearly independent from $R_{1}$.

\vspace{0.5em}

\noindent \textbf{Finishing Up:} We have thus shown that $\mathfrak{h}^{1}_{S} \subset \mathfrak{h}$, and so by symmetry, that $\mathfrak{h}^{-1}_{S} \subset \mathfrak{h}$. Since $\mathfrak{h}$ is by construction an $\textrm{ad} \, \mathfrak{m}$-module containing $\bigoplus_{|i| \geq 2} \mathfrak{h}^{i}_{S}$, this implies that $\mathfrak{h} = \mathfrak{h}_{S}$. As $\mathfrak{h} \subset \mathfrak{u}$, we have $\mathfrak{u} = \mathfrak{h}_{S}$, contradicting an assumption and completing the proof.

\section{Recollections from Cycle Theory}
\label{recollsec}

In the next section will apply the theory we have developed to the concrete situation arising from the relative algebraic de Rham cohomology of a smooth projective family $f : X \to S$ defined over a ring $R \subset \mathbb{C}$, with $X$ and $S$ both smooth over $R$, $S$ connected, and such that $X$ has geometrically connected fibres. We begin by introducing the triple $(\mathcal{H}, F^{\bullet}, \nabla)$ of interest. First, we denote by $\Omega^{\bullet}_{X/S}$ the relative de Rham complex associated to the morphism $f$. We also define a filtration of this complex by $F^{i} \Omega^{\bullet}_{X/S} = \Omega^{\bullet \geq i}_{X/S}$. We fix a cohomological degree $w$, define $\mathcal{H} = R^{w} f_{*} \Omega^{\bullet}_{X/S}$, and
\[ F^{i} \mathcal{H} = \textrm{image} \left[ R^{w} f_{*} F^{i} \Omega^{\bullet}_{X/S} \to R^{w} f_{*} \Omega^{\bullet}_{X/S} \right] . \]
The connection $\nabla$ is constructed as in \cite{katz1968}; note that \cite{katz1968} works over a field, but this is not necessary, see \cite{katznil}. We define a filtration $L^{\bullet}$ on the de Rham complex $\Omega^{\bullet}_{X}$ of $X$ via
\[ L^{i} \Omega^{\bullet}_{X} = \textrm{image} \left[ \Omega^{\bullet-i}_{X} \otimes f^{*}(\Omega^{i}_{S}) \to \Omega^{\bullet}_{X} \right] . \]
We obtain an exact sequence of complexes
\begin{equation}
\label{connexactseq}
0 \to f^{*} \Omega^{1}_{S} \otimes \Omega^{\bullet-1}_{X/S} \to \Omega^{\bullet}_{X} / L^{2} \Omega^{\bullet}_{X} \to \Omega^{\bullet}_{X/S} \to 0 , 
\end{equation}
and we define $\nabla$ to be the connecting homomorphism
\[ \mathcal{H} = R^{w} f_{*} \Omega^{\bullet}_{X/S} \to \Omega^{1}_{S} \otimes R^{w+1} f_{*} \Omega^{\bullet-1}_{X/S} = \Omega^{1}_{S} \otimes \mathcal{H} . \]

Later, in \S\ref{primcupsec} and onwards, we will replace the cohomology sheaf $\mathcal{H}$ with its primitive part; this will be unnecessary prior to considering the cup product, so we work with the full cohomology for the time being. 

\subsection{Flat Sections from Algebraic Cycles}
\label{flatchernsec}

In this section we suppose $w = 2k$ is even. Let us moreover suppose that we have closed immersion $\iota : Y \hookrightarrow X$ of $R$-subschemes which is of relative dimension $k$ over $S$. We explain how this induces a flat section $c^{\textrm{dR}}(Y)$ of $\mathcal{H}$ using the theory of relative Chern classes. To do this we begin by recalling how, to any vector bundle $\mathcal{V}$ on $X$, we may assign a Chern class $c^{\textrm{dR}}(\mathcal{V}) \in \bigoplus_{j} H^{2j}_{\textrm{dR}}(X/S)$ in the de Rham cohomology ring of the morphism $f$. First we consider the case where $\mathcal{V} = \mathcal{L}$ is a line bundle, and define 
\[ c^{\textrm{dR}}(\mathcal{L}) = 1 + c^{\textrm{dR}}_{1}(\mathcal{L}) = 1 + H^{2}(\textrm{dlog})(\mathcal{L}) , \]
where $\textrm{dlog} : \mathcal{O}^{\times}_{X}[-1] \to \Omega^{\bullet}_{X/S}$ is the map of complexes given by $g \mapsto dg / g$, and we use the identification $\textrm{Pic}(X) = H^{2}(X, \mathcal{O}^{\times}_{X}[-1])$. Using \cite[\href{https://stacks.math.columbia.edu/tag/0FI5}{Proposition 0FI5}]{stacks-project}, the map $c^{\textrm{dR}}_{1}$ determines a unique map
\[ c^{\textrm{dR}} : K_{0}(\textrm{Vect}(X)) \to \bigoplus_{i \geq 0} H^{2i}_{\textrm{dR}}(X/S) , \]
compatible with pullback along morphisms $X' \to X$ of quasi-compact and quasi-separated $S$-schemes, and where $K_{0}(\textrm{Vect}(X))$ is the free abelian group on isomorphisms classes of vector bundles with direct sum as addition and modulo the relations $[\mathcal{F}] - [\mathcal{F}'] - [\mathcal{F}'']$ for every exact sequence $0 \to \mathcal{F}' \to \mathcal{F} \to \mathcal{F}'' \to 0$. The map $c^{\textrm{dR}}$ satisfies $c^{\textrm{dR}}(0) = 1$ and $c^{\textrm{dR}}(\alpha + \beta) = c^{\textrm{dR}}(\alpha) c^{\textrm{dR}}(\beta)$. If additionally $R$ is a field, we may use the fact that $X$ is smooth to find a finite resolution $\mathcal{F}_{\bullet} \to \iota_{*} \mathcal{O}_{Y}$ by vector bundles, and finally define the cycle class $[Y]$ as the image under $c^{\textrm{dR}}$ of the alternating sum $\sum (-1)^{i} [\mathcal{F}_{i}]$. 

The maps $c^{\textrm{dR}}$ give, by composing with the natural maps $H^{2j}_{\textrm{dR}}(X/S) \to (R^{2j} f_{*} \Omega^{\bullet}_{X/S})(S)$, sections of $(\bigoplus_{j \geq 0} R^{2j} f_{*} \Omega^{\bullet}_{X/S})(S)$. We now prove:

\begin{prop}
Let $X$ be a quasi-compact and quasi-separated scheme which is smooth over $S$. Then for each $j$, the image of 
\[ c^{\textrm{dR}}_{j} : K_{0}(\textrm{Vect}(X)) \to (R^{2j} f_{*} \Omega^{\bullet}_{X/S})(S) , \]
lies in the subspace of sections flat for the Gauss-Manin connection $\nabla$.
\end{prop}

\begin{proof}
We note as in \cite[Theorem 1]{katz1968} that the Gauss-Manin connection $\nabla$ is compatible with the cup product via the formula $\nabla (\alpha \cup \beta) = \alpha \cup (\nabla \beta) + (\nabla \alpha) \cup \beta$, and is also compatible with pullback, in the sense that if $g : X' \to X$ is a map between quasi-compact and quasi-separated schemes smooth over $S$, then the natural maps $R^{j} f_{*} \Omega^{\bullet}_{X/S} \to R^{j} (f \circ g)_{*} \Omega^{\bullet}_{X'/S}$ respect the connections $\nabla$ and $\nabla'$, where $\nabla$ (resp. $\nabla'$) is the Gauss-Manin connection associated to $f$ (resp. $f' = f \circ g$). 

Consider a vector bundle $\mathcal{V}$ on $X$. If we consider the projective space $g : \mathbb{P}(\mathcal{V}) \to X$ over $X$ associated to $\mathcal{V}$, then the splitting principle tells us that $g^{*} \mathcal{V}$ splits as a direct sum $\mathcal{L} \oplus \mathcal{U}$, where $\mathcal{L}$ is a line bundle. The map $g$ is smooth, hence $X' = \mathbb{P}(\mathcal{V})$ is smooth over $S$, and the maps $R^{j} f_{*} \Omega^{\bullet}_{X/S} \to R^{j} (f \circ g)_{*} \Omega^{\bullet}_{X'/S}$ on cohomology are injective. We then find that showing $\nabla c^{\textrm{dR}}(\mathcal{V}) = 0$ reduces to showing that
\[ \nabla c^{\textrm{dR}}(g^{*} \mathcal{V}) = (\nabla c^{\textrm{dR}}(\mathcal{L})) \cup c^{\textrm{dR}}(\mathcal{U}) + c^{\textrm{dR}}(\mathcal{L}) \cup \left( \nabla c^{\textrm{dR}}(\mathcal{U}) \right) = 0 . \]
The claim therefore reduces inductively to the situation where $\mathcal{V} = \mathcal{L}$ is a line bundle.

In this case we proceed by directly computing $\nabla c^{\textrm{dR}}(\mathcal{L})$ using \v{C}ech cohomology, loosely following \cite{katz1968}. We fix an affine open cover $\mathcal{U} = \{ U_{\alpha} \}_{\alpha \in I}$ of $X$ and for each of the complexes $\mathcal{F}^{\bullet}$ making up the exact sequence (\ref{connexactseq}) we let $C^{q}(\mathcal{U}, \mathcal{F}^{p})$ be the set of alternating $q$-cochains $\beta$ whose value on $\beta(i_{0}, \hdots, i_{q})$ is an element of $\mathcal{F}^{p}(U_{i_{0}} \cap \cdots \cap U_{i_{q}})$. After refining $\mathcal{U}$ we may assume that the restrictions $\restr{\mathcal{L}}{U_{\alpha}}$ are trivial. We have natural differential $d : C^{q}(\mathcal{U}, \mathcal{F}^{\bullet}) \to C^{q}(\mathcal{U}, \mathcal{F}^{\bullet+1})$ given by $(d\beta)(i_{0}, \hdots, i_{q}) = d(\beta(i_{0}, \hdots, i_{q}))$. There is a second differential $\delta : C^{\bullet}(\mathcal{U}, \mathcal{F}^{p}) \to C^{\bullet+1}(\mathcal{U}, \mathcal{F}^{p})$ defined by $\beta(i_{0}, \hdots, i_{q+1}) = (-1)^{p} \sum_{j=0}^{q+1} (-1)^{j} \beta(i_{0}, \hdots, \widehat{i_{j}}, \hdots, i_{q+1})$. 

To each of the complexes $\mathcal{F}^{\bullet}$ appearing in (\ref{connexactseq}) we may associate a total complex $K^{\bullet}(\mathcal{F}^{\bullet})$ with $K^{n}(\mathcal{F}^{\bullet}) = \bigoplus_{p+q = n} C^{q}(\mathcal{U}, \mathcal{F}^{p})$ and with differential $d + \delta$. In particular, the exact sequence (\ref{connexactseq}) induces an exact sequence
\begin{equation}
\label{connexactseq2}
0 \to K^{\bullet}(f^{*} \Omega^{1}_{S} \otimes \Omega^{\bullet-1}_{X/S}) \to K^{\bullet}(\Omega^{\bullet}_{X} / L^{2} \Omega^{\bullet}_{X}) \to K^{\bullet}(\Omega^{\bullet}_{X/S}) \to 0 .
\end{equation}
The class $c^{\textrm{dR}}_{1}(\mathcal{L})$ appearing in $K^{1}(\Omega^{\bullet}_{X/S})$ is constructed as follows. We let $\{ g_{\alpha\alpha'} \}_{(\alpha, \alpha') \in I^2}$ be the transition functions associated to the line bundle $\mathcal{L}$ by the cover $\mathcal{U}$. These functions satisfy the cocycle condition $g_{\alpha\alpha'} g_{\alpha'\alpha''} = g_{\alpha \alpha''}$ on the triple intersections $U_{\alpha} \cap U_{\alpha'} \cap U_{\alpha''}$. Then as in \cite[\href{https://stacks.math.columbia.edu/tag/0FLE}{Section 0FLE}]{stacks-project}, we obtain elements of both $K^{1}(\Omega^{\bullet}_{X} / L^{2} \Omega^{\bullet}_{X})$ and $K^{1}(\Omega^{\bullet}_{X/S})$ by constructing a cocyle $\beta$ whose value on $U_{\alpha \alpha'}$ is given by $g^{-1}_{\alpha \alpha'} dg_{\alpha \alpha'}$. We then have that
\begin{align*}
(d\beta)(\alpha_{0}, \alpha_{1}) &= d(\beta(\alpha_{0}, \alpha_{1})) \\
&= -\frac{dg_{\alpha_{0} \alpha_{1}}}{g_{\alpha_{0} \alpha_{1}}^2} \wedge dg_{\alpha_{0} \alpha_{1}} = 0 \\
(\delta \beta)(\alpha_{0}, \alpha_{1}, \alpha_{2}) &= (-1)(g^{-1}_{\alpha_{0} \alpha_{1}} dg_{\alpha_{0} \alpha_{1}} - g^{-1}_{\alpha_{0} \alpha_{2}} dg_{\alpha_{0} \alpha_{2}} + g^{-1}_{\alpha_{1} \alpha_{2}} dg_{\alpha_{1} \alpha_{2}}) \\
&= g^{-1}_{\alpha_{0} \alpha_{2}} ( g_{\alpha_{1} \alpha_{2}} dg_{\alpha_{0} \alpha_{1}} + g_{\alpha_{0} \alpha_{1}} dg_{\alpha_{1} \alpha_{2}}  ) - g^{-1}_{\alpha_{0} \alpha_{1}} d g_{\alpha_{0} \alpha_{1}} - g^{-1}_{\alpha_{1} \alpha_{2}} d g_{\alpha_{1} \alpha_{2}} \\
&= ( g_{\alpha_{2} \alpha_{0}} g_{\alpha_{1} \alpha_{2}} - g_{\alpha_{1} \alpha_{0}} ) d g_{\alpha_{0} \alpha_{1}} + (g_{\alpha_{2} \alpha_{0}} g_{\alpha_{0} \alpha_{1}} - g_{\alpha_{2} \alpha_{1}}) d g_{\alpha_{1} \alpha_{2}} \\
&= 0 .
\end{align*}
Because the Gauss-Manin connection is defined via the connecting homomorphism associated to the complex (\ref{connexactseq2}), and because the class in $h^{1}(K^{\bullet}(\Omega^{\bullet}_{X/S}))$ which represents $c^{\textrm{dR}}_{1}(\mathcal{L})$ admits a lift to $K^{1}(\Omega^{\bullet}_{X}/L^{2} \Omega^{\bullet}_{X})$ which lies in the kernel of the differential $d + \delta$, it follows from the definition of the connecting homomorphism that $\nabla c^{\textrm{dR}}_{1}(\mathcal{L}) = 0$. 
\end{proof}

\subsection{Flat Sections from the Crystalline Site}
\label{crysflatsec}

Next, we describe a second source of flat sections associated to $(\mathcal{H}, F^{\bullet}, \nabla)$ coming from crystalline cohomology. Here we assume that $f : X \to S$ is a smooth projective family over $\kappa$, where $\kappa$ is a perfect field of characteristic $p > 0$, and with both $X$ and $S$ smooth; the triple $(\mathcal{H}, F^{\bullet}, \nabla)$ is defined as before. 

We begin by recalling that associated to any smooth $\kappa$-variety $Y$, we have a crystalline site $\textrm{Cris}(Y/W)$ which is the category of triples $(U, T, \delta)$ satisfying:
\begin{itemize}
\item[(i)] $U \subset Y$ is an open subscheme;
\item[(ii)] $U \hookrightarrow T$ is an infinitesimal thickening;
\item[(iii)] $\delta$ is a divided power structure on the ideal $I \subset \mathcal{O}_{T}$ defining $U$. (For the definition of divided power structure, we refer to \cite[\href{https://stacks.math.columbia.edu/tag/07GK}{Section 07GK}]{stacks-project}.)
\end{itemize}
The crystalline site of $Y$ has on it a natural sheaf $\mathcal{O}_{Y/W}$, called the crystalline structure sheaf, whose value on $(U, T, \delta)$ is $\Gamma(T, \mathcal{O}_{T})$. Given any sheaf $\mathcal{F}$ on $\textrm{Cris}(Y/W)$ and a map $g : (U, T, \delta) \to (U', T', \delta')$ of objects in $\textrm{Cris}(Y/W)$, we obtain a natural map $c_{g} : g^{-1} \mathcal{F}_{T'} \to \mathcal{F}_{T}$, where the subscript denotes restriction of $\mathcal{F}$. The sheaf $\mathcal{F}$ is called a \emph{crystal} if the maps $c_{g}$ are isomorphisms for any such $g$. In particular, $\mathcal{O}_{Y/W}$ is a crystal.

Another important sheaf on $\textrm{Cris}(Y/W)$ is the sheaf $\Omega^{1}_{Y/W}$ of divided power differentials. This is the quotient of the sheaf $(U, T, \delta) \mapsto \Omega^{1}_{T}$ by the relations $d \delta_{n}(x) = \delta_{n-1}(x) dx$ for every $n \geq 1$ and $x$ in the kernel of $\mathcal{O}_{T} \to \mathcal{O}_{U}$. An important fact about crystals $\mathcal{F}$ on $\textrm{Cris}(Y/W)$ which are modules over $\mathcal{O}_{Y/W}$ is that each comes with a canonical associated \emph{connection}, i.e., a morphism
\[ \nabla : \mathcal{F} \to \mathcal{F} \otimes_{\mathcal{O}_{Y/W}} \Omega^{1}_{Y/W} \]
such that $\nabla (\alpha s) = \alpha \nabla (s) + s \otimes d \alpha$ whenever $s$ is a local section of $\mathcal{F}$ and $\alpha$ is a local section of $\mathcal{O}_{Y/W}$. The morphism $f$ induces a natural pushforward map $f_{*}$ from sheaves on $\textrm{Cris}(X/W)$ to sheaves on $\textrm{Cris}(S/W)$.

\begin{prop}[Berthelot]
\label{berthref}
For each $w$, the sheaf $R^{w} f_{*} \mathcal{O}_{X/W}$ is a crystal, whose value on $(S, S, 0)$ agrees with the algebraic de Rham cohomology $R^{w} f_{*} \Omega^{\bullet}_{X/S}$. Moreover, the connection on $R^{w} f_{*} \Omega^{\bullet}_{X/S}$ induced by the crystal $R^{w} f_{*} \mathcal{O}_{X/W}$ is the Gauss-Manin connection.
\end{prop}

\begin{proof}
This is \cite[Ch V 3.6.2, 3.6.4]{Berthelot1974}
\end{proof}

Using \autoref{berthref}, the following Lemma will allow us to produce flat sections of $\mathcal{H} = R^{w} f_{*} \Omega^{\bullet}_{X/S}$:

\begin{lem}
\label{globalisflat}
Let $\mathcal{F}$ be a crystal on $\textrm{Cris}(Y/W)$, with associated connection $\nabla$. Then any global section $\nu$ of $\mathcal{F}$ is flat for $\nabla$, i.e., we have $\nabla(\nu) = 0$.
\end{lem}

\begin{proof}
Let us describe the construction of $\nabla$, following \cite{2011arXiv1110.5001B}. Let $(U, T, \delta)$ be an object in the crystalline site of $Y$. Then associated to $T$ is a scheme $T(1)$ and a map $\Delta : T \to T(1)$ which is a divided power thickening; as explained in \cite{2011arXiv1110.5001B}, the scheme $T(1)$ is the divided power envelope of $U$ inside $T \times_{W} T$. From the two projections $\textrm{pr}_{i} : T \times_{W} T \to T$ with $i \in \{ 1, 2 \}$ we obtain two isomorphisms $c_{i} : \textrm{pr}^{-1}_{i}\mathcal{F}_{T} \to \mathcal{F}_{T(1)}$. The connection $\nabla$ is then defined on sections $s \in \Gamma(T, \mathcal{F}_{T})$ by $\nabla(s) = c_{1}(s \otimes 1) - c_{2}(1 \otimes 2)$. In the case where $s$ comes from a global section on the crystalline site, this simply means that one has, for each object $(U', T', \delta')$ of $\textrm{Cris}(Y/W)$, a section $s_{(U', T', \delta')} \in \Gamma(T', \mathcal{F}_{T'})$, and this family of sections is compatible with restriction. Hence both $c_{1}(s \otimes 1)$ and $c_{2}(1 \otimes s)$ are equal to some section $\widetilde{s} \in \Gamma(T(1), \mathcal{F}_{T(1)})$, hence the result. 
\end{proof}

We note that the flat sections of $\mathcal{H}$ arising from \autoref{globalisflat} includes all the sections arising from algebraic cycles considered in the previous section: each $\kappa$-subvariety $Y \subset X$ of relative dimension $k$ over $S$ induces, via the natural map $H^{2k}_{\textrm{crys}}(X/W) \to H^{0}_{\textrm{cris}}(S, R^{2k} f_{*} \mathcal{O}_{X/W})$ arising from the edge map of the Leray spectral sequence
\[ E^{ab}_{2} = H^{a}_{\textrm{cris}}(S, R^{b} f_{*} \mathcal{O}_{X/W}) \implies H^{a+b}_{\textrm{cris}}(X/W) , \]
a global section of $R^{2k} f_{*} \mathcal{O}_{X/W}$ (c.f. \cite[Remark 1.1]{morrow:hal-02385601}). 

\subsection{Compatibility with Filtrations}
\label{filcompsec}

It will not be enough that the sections induced by algebraic cycles and global classes invariant under the geometric Frobenius be flat, and we will need additional compatibility with both the Hodge filtration and another filtration on algebraic de Rham cohomology called the \emph{conjugate} filtration. We begin by recalling this latter notion, which requires us to specialize to the situation where $S$ is a scheme of characteristic $p$ (so $p \mathcal{O}_{S} = 0$).

We begin by defining an additional filtration on the algebraic de Rham complex $\Omega^{\bullet}_{X/S}$, complementary to the Hodge filtration described above. Its terms are given by
\[ \tau_{p}(\Omega^{\bullet}_{X/S}) = 0 \to \Omega^{1}_{X/S} \to \cdots \to \Omega^{p-1}_{X/S} \to \textrm{ker}\ d \to 0 \to \cdots . \]
The filtration $\tau_{\bullet}$ induces a spectral sequence, called the ``second spectral sequence of hypercohomology'', denoted
\[ E^{p,q}_{2} = R^{p} f_{*} ( h^{q}(\Omega^{\bullet}_{X/S})) \implies R^{p+q} f_{*}(\Omega^{\bullet}_{X/S}) . \]
The resulting filtration on algebraic de Rham cohomology we denote by $F^{w-\bullet}_{c}$. Its formation is compatible with base-change along maps $S' \to S$ of characteristic $p$-schemes, see \cite{Katz1972}. Note that while much of the literature takes $F^{\bullet}_{c}$ to be an increasing filtration, we use a decreasing one, compatibility with our notation in \S\ref{subflagsec}. We do this for consistency with the conjugate Hodge filtration in characteristic zero. 

Let us now assume that $Y$ is a proper variety over a perfect field $\kappa$ of characteristic $p$, and let $W = W(\kappa)$ be the associated ring of Witt vectors. Then the Hodge filtration $F^{\bullet}$ and conjugate filtration $F^{\bullet}_{c}$ on $H^{\bullet}_{\textrm{dR}}(Y)$ satisfy the following compatibilities with the crystalline Frobenius, as shown by Mazur in \cite[Theorem 7.6]{10.2307/1970906}:

\begin{thm}[Mazur]
\label{Mazurthm}
Assume that $H^{\bullet}_{\textrm{cris}}(Y/W)$ has no $p$-torsion, and that the Fr\"olicher spectral sequence of $Y$ degenerates at $E_{1}$. Denote by $\varphi : H^{\bullet}_{\textrm{cris}}(Y/W) \to H^{\bullet}_{\textrm{cris}}(Y/W)$ the natural Frobenius morphism. Then
\begin{align*}
F^{i} &= \varphi^{-1}(p^{i} H^{w}_{\textrm{cris}}(X/W)) \textrm{ mod } p , \\
F^{w-i}_{c} &= \frac{\textrm{Im}(\varphi) \cap (p^{i} H^{w}_{\textrm{cris}}(X/W))}{p^{i}}  \textrm{ mod } p . 
\end{align*}
\end{thm}

\begin{cor}
\label{actbypicor}
Suppose we are in the situation where $f : X \to S$ is a $\kappa$-algebraic family obtained by base-change to a fibre $\Spec \kappa \to \Spec R$ from the smooth projective $R$-algebraic family $\widetilde{f} : \mathcal{X} \to \mathcal{S}$, where $R \subset \mathbb{C}$ and where $\mathcal{S}$ is smooth over $R$. Suppose that $R^{2k} \widetilde{f}_{*} \Omega^{\bullet}_{\mathcal{X}/\mathcal{S}}$ is a vector bundle. Then for each point $s \in S(\kappa)$, a vector in $H^{2k}_{\textrm{dR}}(X_{s})$ which is the reduction of an eigenvector in $H^{2k}_{\textrm{cris}}(X_{s}/W)$ for $\varphi$ with eigenvalue $p^{k}$ lies in both $F^{k}$ and $F^{k}_{c}$.
\end{cor}

\begin{proof}
The hypotheses ensure, by applying Hensel's Lemma, that $s$, and hence $X_{s}$, admits a characteristic zero lift $\widetilde{s}$. By identifying the crystalline cohomology with the de Rham cohomology of the lift we find that the degeneration hypothesis of \autoref{Mazurthm} is satisfied. For the $p$-torsion hypothesis one equates the reduction of $H^{2k}_{\textrm{cris}}(X_{s}/W)$ to $\kappa$ with the algebraic de Rham cohomology and uses the fact that $R^{2k} \widetilde{f}_{*} \Omega^{\bullet}_{\mathcal{X}/\mathcal{S}}$ is a vector bundle over the smooth scheme $S$, hence has constant fibre dimension. The result then follows by taking $Y = X_{s}$ and applying \autoref{Mazurthm}. 
\end{proof}

\noindent Note that if $\nu$ is a flat section of $R^{2k} f_{*} \Omega^{\bullet}_{X/S}$ induced by a family of algebraic cycles, then the fibres of $\nu$ satisfy the hypotheses of \autoref{actbypicor}.

\subsection{Primitive Cohomology and Cup Product}
\label{primcupsec}

So far we have taken the triple $(\mathcal{H}, F^{\bullet}, \nabla)$ to correspond to the ``full'' cohomology in degree $w$. To emphasize this, we refer to this triple in this section as $(\mathcal{H}_{\textrm{full}}, F^{\bullet}_{\textrm{full}}, \nabla_{\textrm{full}})$, and define a ``primitive'' subbundle $\mathcal{H} \subset \mathcal{H}_{\textrm{full}}$; the notation $(F^{\bullet}, \nabla)$ will the be repurposed for the restriction of $(F^{\bullet}_{\textrm{full}}, \nabla_{\textrm{full}})$ to $\mathcal{H}$.  We will also construct a symmetric non-degenerate bilinear pairing $\mathcal{H} \otimes \mathcal{H} \to \mathcal{O}_{S}$ by combing cup product with the trace pairing. 

The construction works as follows. The morphism $f : X \to S$ is projective, so we may fix a line bundle $\mathcal{L}$ on $X$ which is very ample over $S$. Using the same construction as in \S\ref{flatchernsec}, one may associate to $\mathcal{L}$ a relative Chern class $\omega = c_{1,S}(\mathcal{L})$ which is a global section of $R^{2} f_{*} \Omega^{\bullet}_{X/S}$. If $f$ has relative dimension $n$, one also has a canonical isomorphism $\textrm{tr} : R^{2n} f_{*} \Omega^{\bullet}_{X/S} \xrightarrow{\sim} \mathcal{O}_{S}$ which is described in \cite{PMIHES_1975__45__5_0}. The class $\omega$ induces, for each $w$, an operator $L : R^{w} f_{*} \Omega^{\bullet}_{X/S} \to R^{w+2} f_{*} \Omega^{\bullet}_{X/S}$ by means of cup product with $\omega$, i.e., we have $L(\beta) = \beta \wedge \omega$. Recalling that $\mathcal{H}_{\textrm{full}} = R^{w} f_{*} \Omega^{\bullet}_{X/S}$, we define $\mathcal{H}$ to be the kernel of $L^{n-w+1}$ inside $\mathcal{H}_{\textrm{full}}$. The pairing $\mathcal{Q} : \mathcal{H} \otimes \mathcal{H} \to \mathcal{O}_{S}$ we define by 
\[ \mathcal{Q}(\alpha, \beta) = \textrm{tr}(\alpha \wedge \beta \wedge \omega^{n-w}) . \]

The same construction on the level of the variation $\mathbb{V}_{\textrm{full}} = R^{w} \an{f_{\mathbb{C}}} \mathbb{Z}$ gives a primitive subsystem $\mathbb{V} \subset \mathbb{V}_{\textrm{full}}$, and this construction is compatible with the comparison between algebraic de Rham and Betti cohomology. Finally, we record the following fact regarding the compatibility of the Hodge and conjugate filtrations in positive characteristic and the pairing $\mathcal{Q}$:

\begin{prop}
\label{polposchar}
If $\mathfrak{p} \in \Spec R$ is a prime for which the associated residue field $\kappa(\mathfrak{p})$ has positive characteristic, and $F^{\bullet}_{c}$ is the conjugate filtration on $\mathcal{H}_{\kappa(\mathfrak{p})}$, then $\mathcal{Q}$ lies inside $F^{\textrm{mid}} \cap F^{\textrm{mid}}_{c}$.
\end{prop}

\begin{proof}
What needs to be checked by \autoref{polisHodge} is that both $F^{\bullet}$ and $F^{\bullet}_{c}$ satisfy (\ref{poleq}). This is immediate from \cite[2.3.5.1.4, 2.3.5.1.5]{Katz1972}. (Note that lines (2.3.5.1.4) and (2.3.5.1.5) in \cite{Katz1972} contain a misprint: $N+1-n$ should be $N+i-n$.)
\end{proof}

\section{Applications}
\label{appsec}

We now give our main results in the setting where $(\mathbb{V}, \mathcal{H}, F^{\bullet}, \nabla, \mathcal{Q})$ is the primitive cohomological data in cohomological degree $w$ associated to a smooth projective family $f : X \to S$ over $R \subset \mathbb{C}$, as introduced in \S\ref{primcupsec}.  We let $Q : \mathbb{Z}^{m} \otimes \mathbb{Z}^{m} \to \mathbb{Z}$ be a fixed polarization of the same type as $\mathcal{Q}$; i.e., so that $(\mathbb{Z}^{m}, Q)$ is isomorphic as a polarized lattice to $(\mathbb{V}_{s}, \mathcal{Q}_{s})$ for some (hence any) point $s \in S(\mathbb{C})$.

\subsection{Proof of Theorem 1.1}

We will begin with proving \autoref{babythm}, which will serve as a ``warm up''.

\begin{proof}[Proof of \autoref{babythm}:]
Let us first note that to prove the statement of the theorem we are free to replace $S$ with a non-empty open subscheme. This in particular means we can replace $R$ with a non-empty open subscheme $\Spec R' \subset \Spec R$, as we can always replace $S$ with the inverse image of $\Spec R'$ under the structure map $S \to \Spec R$. 

Denote by $e_{1} \in \mathbb{Z}^{m}$ the first standard basis vector, and let $U = \textrm{span}_{\mathbb{Q}} \{ e_{1}, Q \}$. Consider the family $g : \mathcal{Y} \to \mathcal{M}$ of subschemes of $\ch{L}_{R}$, with $\mathcal{M} = \GL_{m,R}$, whose fibre above a point $A \in \GL_{m}(T)$, with $T$ an $R$-scheme, is the set of flags $F'^{\bullet} \in \ch{L}(T)$ which are polarized by $A^{t} Q_{T} A$ and for which $A e_{1,T} \in F'^{\textrm{mid}}$. The family $g$ has a natural integral structure over $R$, and has constant fibre dimension $e$. We will construct $E$ as a union $\bigcup_{d = 1}^{n} E^{d}$ with $n = \dim S_{\mathbb{C}}$. 

Let $\Phi^{d}_{r}$ and $\mathcal{E}^{d}_{r}$ be as in the statement of \autoref{expintprop}; we write $\mathcal{E}^{d}_{r}$ instead of $\mathcal{E}_{r}$ to emphasize the dependence on $d$. For each $d \in \{ 1, \hdots, n \}$, the locus $E^{d}$ will be defined as the closure of the projection of $\mathcal{E}^{d}_{r}$ to $S$ for some $r = r(d)$ depending on $d$. In the case where $d$ satisfies the inequality $\overline{P} - e > \dim_{R} S - d$ of \autoref{expintprop}, we simply take $r = r(d)$ to be the $r$ given to us by \autoref{expintprop}. Let us write 
\[ E^{\textrm{exc}} = \bigcup_{\overline{P} - e > \dim_{R} S - d} E^{d} , \]
for the subscheme of $E$ which we have so far constructed.

Let us now consider the case where $\overline{P} - e \leq \dim_{R} S - d$. We claim that there is some $r = r(d)$ for which $\mathcal{E}^{d}_{r,\mathbb{C}}$ maps into $E^{\textrm{exc}}_{\mathbb{C}}$. There is nothing to show unless $\mathcal{E}^{d}_{r,\mathbb{C}}$ is non-empty for all $r$, so assume so. Then the sets $\mathcal{E}^{d}_{r,\mathbb{C}}$ are compatible with the natural projections $\mathcal{P}^{d}_{r+1,\mathbb{C}} \times J^{d}_{r+1} \ch{L}_{\mathbb{C}} \to \mathcal{P}^{d}_{r,\mathbb{C}} \times J^{d}_{r} \ch{L}_{\mathbb{C}}$, so we may apply \autoref{existscompseqlem} to find a compatible sequence $\{ (t_{r}, \sigma_{r}) \}_{r \geq 0}$ with $(t_{r}, \sigma_{r}) \in \mathcal{E}^{d}_{r, \mathbb{C}} \subset \mathcal{P}^{d}_{r,\mathbb{C}} \times J^{d}_{r} \ch{L}_{\mathbb{C}}$ for all $r$. Consider the constructible locus $\mathcal{M}_{r} \subset \mathcal{M}(\mathbb{C})$ of points $m$ for which $\sigma_{r} \in J^{d}_{r} \mathcal{Y}_{m}$; by construction each $\mathcal{M}_{r}$ is non-empty and $\mathcal{M}_{r+1} \subset \mathcal{M}_{r}$ for all $r$, so we can find some $m \in \bigcap_{r} \mathcal{M}_{r}$ from constructibility. Using the interpretation of the torsors $\mathcal{P}^{d}_{r}$ given in \autoref{Pinterp}, the point $t_{0}$ corresponds to a frame $\iota : \mathbb{V}_{\mathbb{C},s} \xrightarrow{\sim} \mathbb{C}^{m}$ with $s \in S(\mathbb{C})$ a point, and if $(\psi, \iota)$ is the associated framed local period map then $\psi \circ \gamma(t_{r}) = \sigma_{r}$ for all $r$. If we now apply \autoref{jetfactorlem} to $\psi_{d} = \restr{\psi}{B_{d}}$, with $B_{d} \subset S(\mathbb{C})$ an irreducible $d$-dimensional analytic locus through which the jets $\gamma(t_{r})$ factor, we find that $\psi^{-1}(\mathcal{Y}_{m})$ has dimension at least $d$.

Observe that $\mathcal{Y}_{m}$ is of the form $A \cdot O(A^{-1} \iota(F^{\bullet}_{s}), U_{\mathbb{C}}) = O(\iota(F^{\bullet}_{s}), A \cdot U_{\mathbb{C}})$ for some $A \in \GL_{m}(\mathbb{C})$. Notice that because $\mathcal{M}_{\mathbb{C}} = \GL_{m,\mathbb{C}}$, and because if $(\psi, \iota)$ is a framed local period map then so is $(A \psi, A \iota)$, we may actually choose any $A \in \GL_{m}(\mathbb{C})$ after adjusting $(\psi, \iota)$ appropriately. In particular, we do all of the following:
\begin{itemize}
\item[-] We choose $A$ so that $A \iota$ sends the polarization $\mathcal{Q}_{s}$ on $\mathbb{V}_{\mathbb{C}, s}$ to the polarization $Q$ on $\mathbb{Z}^{m}$. This means that $\psi(B)$ lands inside the complex algebraic subvariety $\ch{D} \subset \ch{L}_{\mathbb{C}}$ of flags satisfying the relation (\ref{poleq}). Note that $\ch{D}$ contains an open locus $D \subset \ch{D}$ of flags $F'^{\bullet}$ for which the subspaces $F'^{p} \cap \overline{F'^{q}}$ give a polarized Hodge structure for the polarized lattice $(\mathbb{Z}^{m}, Q)$.
\item[-] We can further adjust $A$ so that $A e_{1} = e_{1}$, hence so that $\mathcal{Y}_{m} = O(\iota(F^{\bullet}_{s}), U)$.
\item[-] The variety $O(\iota(F^{\bullet}_{s}), U)$ is an orbit of $\mathbf{H}_{U}(\mathbb{C})$, with $\mathbf{H}_{U}$ the pointwise stabilizer of $U$, and intersects $D$; this is a consequence of \cite[VI.B.1]{GGK}. Thus after further acting by some $A \in \mathbf{H}_{U}(\mathbb{C})$ we may even assume that $\iota(F^{\bullet}_{s}) \in D$.
\end{itemize}

We now apply \autoref{noatyplem} to this situation, observing that the hypotheses are satisfied with $F^{\bullet}_{c}$ being the complex conjugate filtration to $\iota(F^{\bullet}_{s})$, and with $\mathbf{U} = \iota(\mathbf{H}_{U})$. Indeed, because the vectors in $U$ are rational, one observes that $U$ is a subspace of Hodge tensors for the Hodge structure with summands $H^{p,q} = \iota(F^{p}_{s}) \cap F^{q}_{c}$. Moreover, the Hodge decomposition on $\textrm{End}(\mathbb{Q}^{m})$ restricts to give a Hodge decomposition on both $\iota(\mathfrak{h}_{S})$ and $\mathfrak{u}$ via the usual mechanism of inducing a Hodge structure on a Mumford-Tate Lie algebra, see for instance \cite[\S II.A]{GGK}. The Hodge symmetry assumption is therefore automatic. We thus obtain that 
\[ \overline{P} - e > \dim_{R} S - \underbrace{\dim \psi^{-1}(O(\iota(F^{\bullet}_{s}), U))}_{d_{0}} , \]
where we note that $O(\iota(F^{\bullet}_{s}), U) = \mathbf{U}(\mathbb{C}) \cdot \iota(F^{\bullet}_{s})$ as a consequence of \autoref{ABCProp}. In particular, we have $E^{d_{0}} \subset E^{\textrm{exc}}$, and the analytic locus $\psi^{-1}(O(\iota(F^{\bullet}_{s}), U))$ lies inside $E^{d_{0}}$ by construction. (One checks this by constructing non-degenerate jets $j$ at smooth points of $\psi^{-1}(O(\iota(F^{\bullet}_{s}), U))$ and using \autoref{bigjetconstr}(iv) together with the definition of $\mathcal{E}^{d_{0}}_{r}$.) This implies that $s \in E^{d_{0}} \subset E^{\textrm{exc}}$.	

We have shown that $\bigcap_{r} \mathcal{E}^{d}_{r,\mathbb{C}}$ maps into $E^{\textrm{exc}}$. From constructibility, this implies the same fact for $\mathcal{E}^{d}_{r,\mathbb{C}}$ for some large enough $r = r(d)$. For this $r$, we define $E^{d}$ to be the closure of the image of $\mathcal{E}^{d}_{r}$ in $S$. It follows that $E = \bigcup_{d = 1}^{n} E^{d}$ is a proper closed subscheme of $S$.

\vspace{0.5em}

Replace $R$ with an open subring so that the $r!$ is invertible in $R$ for every $r = r(d)$ used in the construction of $E$; let $r_{\textrm{max}}$ be the largest $r = r(d)$ that was used, and assume $r_{\textrm{max}} > 2$, increasing it if necessary. After replacing $S$ with an open subscheme, we can also ensure that every $r_{\textrm{max}}$-limp is injective on tangent spaces: the quasi-finiteness assumption on $\varphi$ ensures this is true over an open locus of $S_{\mathbb{C}}$, and hence over an open subscheme of $S$ by spreading out (c.f. the reduction in the proof of \autoref{bigexpthm} below). All that remains to be verified is that the loci $Z$ described in the statement of the theorem lie inside $E$. We may reduce to the case where $Z$ is smooth and has dimension $d = 1$. Let $\overline{\kappa}$ be the Zariski closure of the residue field $\kappa(\mathfrak{p})$, where $Z$ lies above $\mathfrak{p} \in \Spec R$. Let $\nu$ be the global flat section of $\restr{\mathcal{H}}{Z_{\overline{\kappa}}}$ induced by the family $Y \to Z$. Then for any $s \in Z(\overline{\kappa})$ the fibre $\nu_{s}$ lies inside the subspace $\mathcal{U}_{s}$ associated to $\restr{(\mathcal{H}, F^{\bullet}, \nabla)}{Z_{\overline{\kappa}}}$ (recall \autoref{Udef}), and applying \autoref{etaincl} the map $\eta^{1}_{r}$ with $r = r_{\textrm{max}}$ sends $(J^{1}_{r,nd} Z)_{s}(\overline{\kappa})$ into 
\[ \GL_{m}(\overline{\kappa}) \backslash (J^{1}_{r} O(\iota(F^{\bullet}_{s}), \iota(\mathcal{U}_{s})))(\overline{\kappa}) \subset \GL_{m}(\overline{\kappa}) \backslash (J^{1}_{r} O(F'^{\bullet}, U_{\overline{\kappa}})))(\overline{\kappa}) , \]
where $F'^{\bullet} \in \ch{L}_{\overline{\kappa}}$ is any flag for which $U_{\overline{\kappa}} \subset F'^{\textrm{mid}}$. (Here we use the fact that all the possible loci $O(F'^{\bullet}, \textrm{span}\{ Q', v' \})$, with $Q', v'$ an arbitrary polarization and vector over $\overline{\kappa}$, lie in a single $\GL_{m}(\overline{\kappa})$-orbit as long as $\textrm{char}\, \overline{\kappa} \neq 2$, and then take $Q' = \iota(\mathcal{Q}_{s})$ and $v' = \iota(\nu_{s})$.) Unravelling the definitions and using the base-change compatibility in \autoref{bigjetconstr}, this implies that the point $s$ lies below a $\overline{\kappa}$-valued point of $\mathcal{E}^{1}_{r}$, and therefore lies inside $E$.
\end{proof}

\subsection{The General Setting}

Next, we establish our main technical result, which will suffice to deduce the rest of the main theorems. We assume that $\mathbf{H}_{S} = \textrm{Aut}(\mathbb{V}_{s}, Q_{s})$, i.e., that monodromy is Zariski dense. As already mentioned in \S\ref{nondenintro}, the group $\mathbf{H}_{S}$ has a natural realization inside any fibre of $\mathcal{H}$ as the stabilizer of global flat sections of $\bigoplus_{a, b \geq 0} \mathcal{H}^{a,b}$. Suppose we are now given a subvariety $Z \subset S_{\overline{\kappa}}$, where $\kappa = \kappa(\mathfrak{p})$ is the residue field for a positive-characteristic prime $\mathfrak{p} \in \Spec R$. Then for each $a, b \geq 0$ we write $\mathcal{I}^{a,b}_{Z}$ for the set of global flat sections of $\restr{\mathcal{H}^{a,b}}{Z}$ which lie inside both $F^{\textrm{mid}}$ and $F^{\textrm{mid}}_{c}$, where the conjugate filtration on $\restr{\mathcal{H}^{a,b}}{Z}$ is as in \S\ref{filcompsec}. We also define $\mathcal{I}_{Z} = \bigcup_{a,b} \mathcal{I}^{a,b}_{Z}$ and $\mathcal{I}_{Z,N} = \bigcup_{a,b \leq N} \mathcal{I}^{a,b}_{Z}$. If $\mathcal{S} \subset \mathcal{I}_{Z}$ is any subset and $s \in Z$ is a point, we write $\mathbf{H}_{\mathcal{S},s} \subset \textrm{Aut}(\mathcal{H}_{s})$ for the algebraic group stabilizing each element of $\mathcal{S}$. 


\begin{lem}
\label{genericgrouplem}
Suppose $\mathcal{S} \subset \mathcal{I}_{Z}$ contains $\restr{\mathcal{Q}}{Z}$, that $\restr{\mathcal{Q}}{Z}$ is non-degenerate, and that $\textrm{char}\, \kappa(\mathfrak{p}) > 2$. Fix $N > 0$, and suppose $Z$ is irreducible. Then there exists a non-empty open locus $U \subset Z$ where 
\begin{itemize}
\item[(a)] the adjoint Hodge numbers of $\mathbf{H}_{\mathcal{S},s}$, the dimensions $\dim_{\kappa(s)} \mathbf{H}_{\mathcal{S},s}$, and the function
\[ s \mapsto \# \{\textrm{geometric components of } \mathbf{H}_{\mathcal{S},s} \} , \]
are all constant on $U$; 
\item[(b)] there exists a finite-dimensional subspace $\mathcal{T} \subset \mathcal{S}$ such that $\mathbf{H}_{\mathcal{S},s} = \mathbf{H}_{\mathcal{T},s}$ for any $s \in U$; and
\item[(c)] for any $s \in U$, any global section of $\bigoplus_{a,b \leq N} \restr{\mathcal{H}^{a,b}}{Z}$ stabilized by $\mathbf{H}_{\mathcal{S}, s}$ lies inside $F^{\textrm{mid}} \cap F^{\textrm{mid}}_{c}$.
\end{itemize}
\end{lem}

\begin{proof}
For each $s \in Z$, let $\beta_{s} : \mathbb{G}_{m,\kappa(s)} \to \mathbf{H}_{\mathcal{S},s} \subset \textrm{Aut}(\mathcal{H}_{s})$ be the cocharacter of \autoref{ABCProp}; note that the algebraically-closed assumption in \autoref{ABCProp} is not necessary to construct $\beta_{s}$. We obtain induced direct sum decompositions $\mathfrak{h}_{\mathcal{S},s} = \bigoplus_{i} \mathfrak{h}^{i}_{\mathcal{S},s}$. The maps $\mathfrak{h}^{i}_{\mathcal{S},s} \to s$ are fibres of an algebraic family over $Z$, so from irreducibility (c.f. \cite[\href{https://stacks.math.columbia.edu/tag/05F7}{Lemma 05F7}]{stacks-project}) there is a non-empty open subset $U_{i} \subset Z$ where these fibres have constant dimension. Define $U = \bigcap_{i} U_{i}$. Shrinking $U$ further, we may assume by \cite[\href{https://stacks.math.columbia.edu/tag/055H}{Lemma 055H}]{stacks-project} that all the groups $\mathbf{H}_{\mathcal{S},s}$ above geometric points of $U$ have the same number of irreducible components. 

Property (a) is immediate. To show (b), let $s_{0} \in U$ be the (scheme-theoretic) generic point and write $\mathcal{S}_{M} = \mathcal{S} \cap \mathcal{I}_{Z,M}$. Then for large enough $M$, we have that $\mathbf{H}_{\mathcal{S},s_{0}} = \mathbf{H}_{\mathcal{S}_{M}, s_{0}}$ from Noetherianity. If $U_{M}$ is constructed exactly as $U$ above but with respect to $\mathcal{S}_{M}$, then one has $\mathbf{H}_{\mathcal{S},s} \subset \mathbf{H}_{\mathcal{S}_{M},s}$ for each $s \in U \cap U_{M}$, and the inclusion is an equality as the dimensions and number of components is the same. Thus, we may take $\mathcal{T} = \mathcal{S}_{M}$, which shows (b) after shrinking $U$. 

After removing a closed subscheme from $U$, we may assume that any section of $\bigoplus_{a,b \leq N} \restr{\mathcal{H}^{a,b}}{Z}$ fixed by $\mathbf{H}_{\mathcal{S},s}$ for $s \in U$ is also fixed by $\mathbf{H}_{\mathcal{S},s_{0}}$. This reduces us to showing that global sections fixed by $\mathbf{H}_{\mathcal{S},s_{0}}$ lie inside $F^{\textrm{mid}} \cap F^{\textrm{mid}}_{c}$. Using the duality induced by $\mathcal{Q}_{s_{0}}$ we may reduce to considering just those tensors lying in either $\bigoplus_{j \geq 0} \mathcal{H}^{\otimes j}_{s_{0}}$ or $\bigoplus_{j \geq 0} \mathcal{H}^{\otimes j}_{s_{0}} \otimes (\mathcal{H}^{*}_{s_{0}})^{\otimes j}$, depending on the weight of the Hodge decomposition on $\mathcal{H}_{s_{0}}$. Using the interpretation of the Hodge decomposition in terms of weight spaces for $\beta$, we see that such tensors necessarily lie in $F^{\textrm{mid}} \cap F^{\textrm{mid}}_{c}$ because $\beta_{s_{0}}$ fixes them. 
\end{proof}

\begin{thm}
\label{bigexpthm}
Suppose that $\mathbb{V}$ has level at least three, suppose that $\mathbf{H}_{S} = \textrm{Aut}(\mathbb{V}_{s}, \mathcal{Q}_{s})$ for some $s \in S(\mathbb{C})$. Let $N > 0$ be a positive integer. Then there exists a closed subscheme $E$, properly contained in $S$, such that: 
\begin{itemize}
\item[-] every weakly special subvariety $Z \subset S_{\mathbb{C}}$ with $\mathbf{H}_{Z} \not\in \{ 1, \mathbf{H}_{S} \}$ lies inside $E_{\mathbb{C}}$; and
\item[-] for any smooth irreducible locally closed positive-characteristic subscheme $Z \subset S_{\overline{\kappa}}$, where $\kappa = \kappa(\mathfrak{p})$ for some prime $\mathfrak{p} \in \Spec R$ and $\overline{\kappa}$ is an algebraic closure, satisfying
\begin{itemize}
\item[(i)] $\restr{(\mathcal{H}, F^{\bullet}, \nabla)}{Z}$ admits a non-constant $r$-limp of some order $r$;
\item[(ii)] $Z$ admits a Zariski dense set of ordinary points; and
\item[(iii)] there exists a subset $\mathcal{S} \subset \mathcal{I}_{Z,N}$ such that the algebraic group $\mathbf{H}_{\mathcal{S},s_{0}} \subset \textrm{End}(\mathcal{H}_{s_{0}})$ stabilizing the elements of $\mathcal{S}$ defines a proper Lie subalgebra $\mathfrak{h}_{\mathcal{S},s_{0}} \subsetneq \mathfrak{h}_{S,s_{0}}$ with symmetric adjoint Hodge numbers, where $s_{0}$ is the generic point of $Z$; 
\end{itemize}
we have $Z \subset E_{\overline{\kappa}}$.
\end{itemize}
\end{thm}


\begin{proof}[Proof of \autoref{bigexpthm}]
~ \vspace{0.5em}

\noindent \textbf{Preliminary Reductions:} As it causes no harm, we can and will assume that the $N$ that occurs in the statement of \autoref{bigexpthm} is at least as large as the one that occurs in the statement of \autoref{orbiteqlem}. After replacing $\Spec R$ with an open subscheme we can and do assume that $\mathcal{Q}_{s}$ is non-degenerate for each $s \in S$.

Let $K$ be the fraction field of $R$ inside $\mathbb{C}$. Let us begin by noting that it suffices to show the statement after replacing $S$ with an open subscheme, and hence $\Spec R$ with an open subscheme. This in particular tells us that the particular integral model $(\mathcal{H}, F^{\bullet}, \nabla, \mathcal{Q})$ of the data $(\mathcal{H}_{K}, F^{\bullet}_{K}, \nabla_{K}, \mathcal{Q}_{K})$ that we work with is of little relevance: any other such integral model $(\mathcal{H}', F'^{\bullet}, \nabla', \mathcal{Q}')$ will be uniquely isomorphic to $(\mathcal{H}, F^{\bullet}, \nabla, \mathcal{Q})$ over a sufficiently small open subscheme of $\Spec R$, and so one can assume the two models are the same after passing to an open subscheme of $S$. These facts will be used throughout the proof, sometimes without comment.

With this in mind, we may use \cite[\S2]{CM_1981__43_2_253_0} to find a finite covering $Y_{K} \to S_{K}$ and a projective compactification $Y_{K} \subset \overline{Y}$ such that the divisor $\overline{Y} \setminus Y_{K}$ at infinity has normal crossings with unipotent monodromy. After replacing $\Spec R$ with an open subscheme $\Spec R' \subset \Spec R$ and spreading out the map $Y_{K} \to S_{K}$ to an integral map $g : Y \to S$ we may reduce to proving the theorem instead for the data associated to the family $X_{Y} \to Y$ obtained by base-change. Replacing $Y$ with $S$ we may now let $\varphi : \an{S_{\mathbb{C}}} \to \Gamma \backslash D$ be the analytic period map sending the point $s \in S(\mathbb{C})$ to the isomorphism class of the integral polarized Hodge structure $(\mathbb{V}_{s}, \mathcal{Q}_{s})$. As in \cite[Cor. 13.7.6]{CMS}, one can use the unipotence of the monodromy at infinity to construct a proper extension $\overline{\varphi}$ of $\varphi$, and then the main theorem of \cite{OMINGAGA} implies that we obtain a Stein factorization $\overline{\varphi} = u \circ \an{q}_{\mathbb{C}}$ where $q_{K} : S_{K} \to T_{K}$ is a $K$-algebraic map satisfying $q_{K, *} \mathcal{O}_{S_{K}} = \mathcal{O}_{T_{K}}$, and the map $u : \an{T_{\mathbb{C}}} \to \Gamma \backslash D$ is quasi-finite.

Applying the argument of \cite[Lem 3.4]{fieldsofdef}, the data $(\mathbb{V}, \mathcal{H}, F^{\bullet}, \nabla)$ descends to $T_{K}$ in the sense that there exists a quintuple $(\mathbb{V}_{T}, \mathcal{H}_{T_{K}}, F^{\bullet}_{T_{K}}, \nabla_{T_{K}})$ of the same kind on $T_{K}$ agreeing with the data on $S_{K}$ upon pullback via $q_{K}$. Spreading out the data $(\mathcal{H}_{T_{K}}, F^{\bullet}_{T_{K}}, \nabla_{T_{K}})$ we can extend it to data $(\mathcal{H}_{T}, F^{\bullet}_{T}, \nabla_{T})$ defined on an integral model $T$ for $T_{K}$, and moreover spread out $q$ to obtain a map $q : S \to T$, possibly passing to an open subscheme of $R$. After passing to a smooth open subscheme $V \subset T$ and replacing $S$ with $U = q^{-1}(V)$ we have reduced to the following situation:
\begin{quote}
The data $(\mathbb{V}, \mathcal{H}, F^{\bullet}, \nabla)$ is obtained by pullback of data $(\mathbb{V}_{T}, \mathcal{H}_{T}, F^{\bullet}_{T}, \nabla_{T})$ of the same type along a morphism $q : S \to T$ of smooth $R$-schemes, and the map $u : \an{T_{\mathbb{C}}} \to \Gamma \backslash D$, which is a period map for $\mathbb{V}_{T}$, is quasi-finite. 
\end{quote}

We now use this to give an alternative description of the condition (i) in the statement of \autoref{bigexpthm}. Let $c(\ch{L}_{R}) \subset J^{1}_{1} \ch{L}_{R}$ be the subscheme of ``constant'' tangent vectors (i.e., the image of the zero section). This is a $\GL_{m,R}$-invariant locus, and so we may consider the inverse image $(\eta^{1}_{1,T})^{-1}(c(\ch{L}_{R})) \subset J^{1}_{1} T$ where $\eta^{1}_{1,T}$ is as in \autoref{bigjetconstr} associated to the data $(\mathcal{H}_{T}, F^{\bullet}_{T}, \nabla_{T})$. We claim that after possibly replacing both $R$ and $T$ with open subschemes, one has that $(\eta^{1}_{1,T})^{-1}(c(\ch{L}_{R})) \subset c(T)$, where $c(T) \subset J^{1}_{1} T$ is the locus of constant tangent vectors. It suffices to show this in the case where $R = \mathbb{C}$, as the compatibility of $\eta^{1}_{1,T}$ with base-change will enable us to inherit this property from the generic fibre. Using \autoref{bigjetconstr}(iv), what then needs to be shown is that over some open locus in $T_{\mathbb{C}}$ local period maps $\psi : B \to \an{\ch{L}}_{\mathbb{C}}$ are injective on tangent spaces. This follows from the fact that $u$ is quasi-finite, as for every local period map $\psi$ there is some $g \in \GL_{m}(\mathbb{C})$ such that $g \psi = \pi \circ \restr{u}{B}$, where $\pi : D \to \Gamma \backslash D$ is the natural projection. 

From \autoref{bigjetconstr}(ii) and \autoref{bigjetconstr}(iii) the maps $\eta^{d}_{r,S}$ and $\eta^{d}_{r,T}$ obtained from the data, $(\mathcal{H}, F^{\bullet}, \nabla)$ and $(\mathcal{H}_{T}, F^{\bullet}_{T}, \nabla_{T})$, are compatible with both pullback along $q$ and pulling back jets $j : \mathbb{D}^{d}_{r,R} \to S$ along embeddings $\mathbb{D}^{1}_{1,R} \hookrightarrow \mathbb{D}^{d}_{r,R}$; it follows that after replacing $R$ with $R_{r}$ and for any $d$ one has that
\[ (J^{d}_{r} q)^{-1}(J^{d}_{r,nd} T) \subset (\eta^{d}_{r,S})^{-1}\left( J^{d}_{r} \ch{L}_{R} \setminus c(\ch{L}_{R}) \right) ; \]
using \autoref{bigjetconstr}(ii), the interpretation of this is that for any $r$-limp $\psi$ associated jets of the form $\psi \circ j$ are non-constant if $q \circ j$ is non-degenerate. Under these circumstances, one has that 
\begin{quote}
For a positive characteristic $Z \subset S_{\overline{\kappa}}$ as in the statement of \autoref{bigexpthm}, condition (i) implies the existence of an open subscheme $U \subset Z$ such that \emph{any} $r$-limp associated to $\restr{(\mathcal{H}, F^{\bullet}, \nabla)}{U}$ is non-constant to first order. 
\end{quote}
Indeed, to see this let $Y$ be the Zariski closure of $q(Z)$; the fact that $Z$ admits a non-constant $r$-limp and that $\eta^{d}_{r,Z}$ factors through $\restr{\eta^{d}_{r,T}}{Y}$ implies that $Y$ is not a point, and so has dimension $d > 0$. Then there exists open loci $U \subset Z$ and $V \subset Y$ such that both $V$ and the induced map $U \to V$ is smooth. In this situation it follows from the equivalence of smoothness and formal smoothness that $J^{d}_{r} U \to J^{d}_{r} V$ is surjective, hence the inverse image $q^{-1}(J^{d}_{r,nd} T)$ intersects $J^{d}_{r} Z$ above every point of $U$. It then follows from \autoref{bigjetconstr}(ii) that $r$-limps defined on $Z$ above points of $U$ are non-constant to first order.





\vspace{0.5em}

\noindent \textbf{Constructing E:} We are now ready to begin our construction of $E$. We will construct $E$ as a union $E = \bigcup_{p \in \mathcal{P}} E(p)$, where $E(p) \subset S$ is a proper closed locus and $\mathcal{P}$ is a finite set of parameters, i.e., tuples of non-negative integers of the form
\[ \mathcal{P} = \left\{ (d, j, h^{-w},\hdots, h^{w}) : \substack{1 \leq d \leq \dim S_{\mathbb{C}}, \, 0 < j \leq M \\ h^{1} + \cdots + h^{w} < \overline{P}, h_{0} \leq m^2 \\ h^{i} = h^{-i}\textrm{ for all }i, } \right\} , \]
where $\overline{P}$ is the period dimension of \autoref{perdimdef} and we have $M = \dim_{\mathbb{C}} \left( \bigoplus_{a, b \leq N} \mathcal{H}^{a,b}_{\mathbb{C}} \right)$. Fixing such a tuple $p = (d, j, h^{-w}, \hdots, h^{w})$, we now construct a family of subschemes of $\ch{L}_{R}$, depending on $p$, which are ``defined by'' tensors in $R_{N} = \bigoplus_{a,b \leq N} (R^{m})^{a,b}$, functionally in $R$. 

We denote by $\ch{X} = \textbf{Gr}(\mathbb{Z}_{N}, j)$ the Grassmanian of subspaces of $\mathbb{Z}_{N}$ of dimension $j$.  We then consider the locus 
\[ \mathcal{F} \subset \ch{L}_{R} \times \ch{L}_{R} \times \ch{X}_{R} , \]
where a triple $(L^{\bullet}_{c}, L^{\bullet}, W) \in \mathcal{F}(R')$, for $R \to R'$ a ring map, satisfies
\begin{itemize}
\item[-] the $R'$-module $W$ breaks up as a direct sum $W = \bigoplus_{a,b \leq N} W^{a,b}$ with $W^{a,b} \subset (R'^{m})^{a,b}$;
\item[-] the graded locally-free filtrations $L^{\bullet}$ and $L^{\bullet}_{c}$ are opposed;
\item[-] the $R'$-module $W$ lies inside $L^{\textrm{mid}} \cap L^{\textrm{mid}}_{c}$ and contains $Q_{R'}$;
\item[-] if $\mathfrak{s} \subset \mathfrak{gl}_{m,R'}$ is the Lie algebra corresponding to the subgroup $\mathbf{S} \subset \GL_{m,R'}$ which fixes $W$, then the adjoint direct sum decomposition on $\mathfrak{gl}_{m,R'}$ induced by $(L^{\bullet}, L^{\bullet}_{c})$ induces a direct sum decomposition $\mathfrak{s} = \bigoplus_{i} \mathfrak{s}^{i}$ with $\rank \mathfrak{s}^{i} = h^{i}$ for all $i$.
\end{itemize}
One may check that, after possibly passing to an open subscheme of $\Spec R$, $\mathcal{F}$ is a locally closed algebraic locus. We then let $\mathcal{Y}$ be the projection of $\mathcal{F}$ to the factor $\ch{L}_{R} \times \ch{X}_{R}$ with its reduced scheme structure, set $\mathcal{M} = \ch{X}_{R}$ and define $g$ to be the map $g : \mathcal{Y} \to \mathcal{M}$ obtained by projecting onto the last factor. The fibres of $g$ are subschemes of $\ch{L}_{R}$ parameterizing flags $L^{\bullet}$ for which there exists an opposing filtration $L^{\bullet}_{c}$ satisfying the conditions above. When we want to emphasize the dependence of $g$ on $p$ we will write $g = g(p)$. 

Each locus $E(p)$ will be constructed, for an appropriate $r = r(p)$, as the Zariski closure of the projection to $S$ of the locus
\[ \mathcal{E}_{r} \subset \Phi^{d}_{r} \cap (J^{d}_{r,nd} S \times J^{d}_{r,nc} \ch{L}_{R}) , \]
defined as in \autoref{expintprop} relative to the family $g(p)$. We will begin by showing that each locus $E(p)$ is properly contained in $S$. To do this, note that the fibres of the family $g$ have a common dimension $e = e(p)$, as \autoref{ABCProp} shows this dimension is equal to $e = \sum_{i \geq 1} h^{i}$ for any geometric fibre. Our analysis splits into two cases:

\vspace{0.5em}

\noindent \textbf{Exceptional Case:} Consider the loci $E(p)$ in the situation where
\begin{equation}
\label{expineq}
 \overline{P} - e > \dim_{R} S - d .
\end{equation}
Then the statement of \autoref{expintprop} gives us $r = r(p)$ such that $E(p)$ is properly contained in $S$. We denote by $\mathcal{P}^{\textrm{exc}} \subset \mathcal{P}$ the subset of elements $(d, j, h^{-w}, \hdots, h^{w})$ for which (\ref{expineq}) holds.

\vspace{0.5em}

Before continuing, let us observe that every weakly special subvariety $Z \subset S_{\mathbb{C}}$ with $\mathbf{H}_{Z} \not\in \{ 1, \mathbf{H}_{S} \}$ lies inside $E^{\textrm{exc}}_{\mathbb{C}}$, where $E^{\textrm{exc}} = \bigcup_{p \in \mathcal{P}^{\textrm{exc}}} E(p)$. To see this, let us consider a point $s$ in the smooth locus of $Z$, and let $\iota : \mathbb{V}_{s} \xrightarrow{\sim} \mathbb{Z}^{m}$ be a polarization-preserving frame. We may construct a local period map $\psi : B \to \ch{L}$ with frame $\iota$ such that $Z$ is locally defined as $\psi^{-1}(O(\iota(F^{\bullet}_{s}), \iota(\mathcal{U}_{Z,s})))$, where $\mathcal{U}_{Z,s}$ is the restriction to $s$ of global sections of $\bigoplus_{a, b \geq 0} \restr{\mathcal{H}^{a,b}}{Z}$ which are flat for $\nabla$. Note that 
\[ O(\iota(F^{\bullet}_{s}), \iota(\mathcal{U}_{Z,s})) = O(\iota(F^{\bullet}_{s}), \iota(\mathcal{U}^{\textrm{mid}}_{Z,s})) = \iota(\mathbf{H}_{Z}(\mathbb{C})) \cdot \iota(F^{\bullet}_{s}) \]
from \autoref{orbiteqlem}, where $\mathcal{U}^{\textrm{mid}}_{Z,s} \subset \bigoplus_{a,b \leq N} (\mathbb{V}_{\mathbb{C},s})^{a,b} \cap \mathcal{U}_{Z,s} \cap F^{\bullet}_{s} \cap F^{\bullet}_{c}$ and $F^{\bullet}_{c}$ is the complex conjugate of $F^{\bullet}_{s}$.


If we consider the group $\mathbf{S}$ stabilizing $W = \iota(\mathcal{U}^{\textrm{mid}}_{Z,s})$, then one has that $Q \subset W$ by construction, and one moreover knows from \autoref{orbiteqlem} that the adjoint Hodge numbers $h^{-w}, \hdots, h^{w}$ on the Lie algebra of $\mathbf{S}$ are symmetric. It follows that $O(\iota(F^{\bullet}_{s}), \mathcal{U}^{\textrm{mid}}_{Z,s})$ is a fibre of $g(p)$ with $p = (\dim Z, \dim W, h^{-w}, \hdots, h^{w})$. We thus obtain that $Z \subset E(p)$ by construction, noting that the map $\psi$ is non-constant since $\mathbf{H}_{Z}$ is non-trivial. Applying \autoref{noatyplem} with $\mathbf{U} = \iota(\mathbf{H}_{Z})$ and with respect to the conjugate flag $\iota(F^{\bullet}_{c})$ one obtains the inequality 
\[ \underbrace{\dim (\mathbf{H}_{S}(\mathbb{C}) \cdot F^{\bullet}_{s})}_{= \overline{P}} - \underbrace{\dim O(\iota(F^{\bullet}_{s}), \mathcal{U}^{\textrm{mid}}_{Z,s})}_{= e} > \dim S_{\mathbb{C}} - \dim Z . \] 
Hence $p \in \mathcal{P}^{\textrm{exc}}$.


\vspace{1em}

\noindent \textbf{Unexceptional Case:} Next we consider the loci $E(p)$ in the situation where
\begin{equation}
\label{unexpineq}
 \overline{P} - e \leq \dim_{R} S - d .
\end{equation}
In this case we show that for large enough $r$ the generic fibre of the locus $\mathcal{E}_{r}$ maps into $E^{\textrm{exc}} = \bigcup_{p \in \mathcal{P}^{\textrm{exc}}} E(p)$. Since we have just established that $E^{\textrm{exc}}$ is properly contained in $S$, this will let us choose $r = r(p)$ such that each $E(p)$ with $p$ satisfying (\ref{unexpineq}) is as well, hence so will be the union $E = \bigcup_{p \in \mathcal{P}} E(p)$.

From compatibility with base-change it suffices to show that for large enough $r$ the locus $\mathcal{E}_{r,\mathbb{C}}$ maps into $E^{\textrm{exc}}_{\mathbb{C}}$. If write
\[ \mathcal{L} = \bigcap_{r} \textrm{im}(\mathcal{E}_{r,\mathbb{C}} \to \mathcal{P}^{d}_{0,\mathbb{C}} \times \ch{L}_{\mathbb{C}}) , \]
then using constructibility this is the same as showing that $\mathcal{L}$ maps into $E^{\textrm{exc}}_{\mathbb{C}}$. 

There is nothing to show unless $\mathcal{L}$ is non-empty, so we assume so. The loci $\mathcal{E}_{r,\mathbb{C}}$ form a compatible sequence of sets as in \autoref{existscompseqlem}, so for any point $(t_{0}, \sigma_{0}) \in \mathcal{L}$ there exists a compatible sequence $(t_{r}, \sigma_{r}) \in \mathcal{E}_{r}(\mathbb{C})$ such that $(t_{r+1}, \sigma_{r+1})$ projects onto $(t_{r}, \sigma_{r})$ for all $r \geq 0$. The element $t_{0}$ corresponds by \autoref{Pinterp} to the pair $(s, \iota)$, where $\iota : \mathcal{H}_{s} \xrightarrow{\sim} \mathbb{C}^{m}$ is a choice of frame. Using the local period map version of \autoref{existslem} we can construct a framed local period map $(\psi, \iota)$ on an analytic ball $B$ containing $s$ such that $\psi \circ j_{r} = \sigma_{r}$ for all $r$, where $j_{r} = \gamma(t_{r})$. We may furthermore let $B_{d} \subset B$ be a closed analytic locus such that $j_{r} \in J^{d}_{r} B_{d}$ for all $r$.

Now by definition of $\mathcal{E}_{r,\mathbb{C}}$, for each $r$ the constructible locus $\mathcal{M}(\sigma_{r}) \subset \mathcal{M}_{\mathbb{C}}$ of points $m$ for which $\sigma_{r} \in (J^{d}_{r} \mathcal{Y}_{m})(\mathbb{C})$ is non-empty. Since $\mathcal{M}(\sigma_{r+1}) \subset \mathcal{M}(\sigma_{r})$ for all $r$ and all these sets are constructible, it follows that there exists a point $m \in \bigcap_{r} \mathcal{M}(\sigma_{r})$. We recall that $m$ corresponds to a subspace $W \subset \mathbb{C}_{N}$. We thus find by applying \autoref{jetfactorlem} that $\psi^{-1}(\mathcal{Y}_{m})$ contains the locus $B_{d}$, hence has dimension at least $d > 0$. The variety $\mathcal{Y}_{m}$ has as its component through $\iota(F^{\bullet}_{s})$ the variety $O(\iota(F^{\bullet}_{s}), W)$. Let $\mathbf{U}$ be the pointwise stabilizer of $W$. By the construction of the family $g$, there exists a filtration $F^{\bullet}_{c}$, opposed to $\iota(F^{\bullet}_{s})$, such that the pair $(\iota(F^{\bullet}_{s}), F^{\bullet}_{c})$ induces a symmetric direct sum decomposition on the Lie algebra $\mathfrak{u}$ of $\mathbf{U}$. Moreover, the same is true for the Lie algebra $\iota(\mathfrak{h}_{S})$: the fact that we obtain a direct sum decomposition comes from the fact that $Q \in W$ and \autoref{ABCProp}, and the symmetry is due to the fact that the filtration on $\mathfrak{h}_{S}$ induced by $F^{\bullet}_{s}$ has graded pieces with symmetric dimensions. Thus we are in the situation of \autoref{noatyplem}, and we conclude that (\ref{expcodimconc}) holds. Applying the Ax-Schanuel Theorem (\ref{axschanthm}) with $T = S_{\mathbb{C}} \times O(\iota(F^{\bullet}_{s}), W)$, the locus $\psi^{-1}(O(\iota(F^{\bullet}_{s}), W))$ lies inside a proper weakly special subvariety $Z$ of $S_{\mathbb{C}}$ which moreover admits a non-constant local period map (since $\sigma_{r}$ is non-constant). But due to our argument above this means that $s$ lies inside $E^{\textrm{exc}}_{\mathbb{C}}$.

\vspace{1em}

\noindent \textbf{Finishing Up:} The only thing left to show is that each positive characteristic subscheme $Z \subset S_{\overline{\kappa}}$ as in the statement of the theorem lies inside $E$. Note that here we assume that the characteristic of the field associated to $Z$ is at least $r_{\textrm{max}} = \textrm{max}_{p \in \mathcal{P}} \{ r(p) \}$, which we can do after passing to an open subscheme of $R$. Using our initial reduction, we may assume that any $r$-limp associated to $Z$ is non-constant to first order. Moreover, applying \autoref{genericgrouplem} and shrinking $Z$ we may assume that condition (iii) holds for any $s \in Z(\overline{\kappa})$. We set $d = \dim Z$.

It suffices to show that for any ordinary point $s \in Z(\overline{\kappa})$ we have that $s \in E(\overline{\kappa})$. Let $\iota : \mathcal{H}_{s} \xrightarrow{\sim} \overline{\kappa}^{m}$ be a frame sending $\mathcal{Q}_{s}$ to $Q_{\overline{\kappa}}$, and let $\psi = q \circ M$ be an associated framed $r_{\textrm{max}}$-limp for the data $\restr{(\mathcal{H}, F^{\bullet}, \nabla)}{Z}$. Write $\mathcal{S}_{s}$ for the set of restrictions of elements of $\mathcal{S}$ to the fibre above $s$. Then $\mathcal{S}_{s} \subset F^{\textrm{mid}}_{s} \cap F^{\textrm{mid}}_{c}$, where we write $F^{\bullet}_{c}$ for the conjugate filtration on algebraic de Rham cohomology. Because the elements of $\mathcal{S}_{s}$ extend to flat sections over $Z$, we have 
\begin{equation}
\label{impsiliesin}
 \psi((J^{d}_{r,nd} Z)_{s}) \subset J^{d}_{r,nc} O(\iota(F^{\bullet}_{s}), \iota(\mathcal{S}_{s})) .
\end{equation}
It therefore suffices to show that $O(\iota(F^{\bullet}_{s}), W)$, with $W = \iota(\mathcal{S}_{s})$, is a fibre of the family $g(p)$ for 
\[ p = (d, \dim W, \dim \mathfrak{h}^{-w}_{\mathcal{S}}, \hdots, \dim \mathfrak{h}^{w}_{\mathcal{S}}) . \] 
The only thing that needs to be verified is that 
\[ \dim \mathfrak{h}^{1}_{\mathcal{S}} + \cdots + \dim \mathfrak{h}^{w}_{\mathcal{S}} < \overline{P} = \dim \mathfrak{h}^{1}_{S} + \cdots + \dim \mathfrak{h}^{w}_{S} . \]
If not, then $\mathfrak{h}_{\mathcal{S}} \subsetneq \mathfrak{h}_{S}$ lies inside a proper $\textrm{ad}\, \mathfrak{h}^{0}_{S}$ module which contains $\bigoplus_{i \neq 0} \mathfrak{h}^{i}_{S}$. This is impossible in characteristic zero, hence impossible if the characteristic of $\overline{\kappa}$ is large enough by a Lefschetz principle argument. Passing to an open subscheme of $\Spec R$, this completes the proof.
\end{proof}

\subsection{Results on Algebraic Cycle Loci}


\begin{thm}
\label{masterthm}
For large enough $N > 0$, the locus $E \subset S$ specified in \autoref{bigexpthm} contains all the loci $Z$ specified in \autoref{properlocusthm2}. Moreover, this continues to hold if the sets $\mathcal{A}^{a,b}$ (resp. the sets $\mathcal{A}^{a,b}_{\textrm{cris}}$) in \autoref{properlocusthm} are instead taken to consist of the span of those global sections of $\restr{\mathcal{H}^{a,b}}{Z_{\overline{\kappa}}}$ (resp. of $\mathcal{H}^{a,b}_{\textrm{cris}, Z_{\overline{\kappa}}}$) for which a $Y$ as in (\ref{cyclediag2}) is merely expected to exist according to the either the Hodge conjecture or the (variational) crystalline Tate conjectures. 
\end{thm}

\begin{proof}[Proof of (\ref{masterthm}), (\ref{properlocusthm}), (\ref{mereconj}), (\ref{properlocusthm2})] ~ \\

The statement of \autoref{masterthm} implies the others, so we will simply focus on \autoref{masterthm}, beginning by explaining its meaning.

\vspace{0.5em}

\noindent \textbf{Meaning of ``expected to exist'':} We assume that $\kappa = \kappa(\mathfrak{p})$ is a positive characteristic field for the time being. We denote by $h^{\otimes n} : X^{\otimes n}_{Z_{\overline{\kappa}}} \to Z_{\overline{\kappa}}$ the base-change of the map $f^{\otimes n} : X^{\otimes n} \to S$ to $Z_{\overline{\kappa}} \subset S_{\overline{\kappa}}$. Consider a point $s \in Z(\overline{\kappa})$, and let $z \in \textrm{CH}^{i}(X^{\otimes n}_{s})_{\mathbb{Q}}$ be an element of the Chow group of the fibre in codimension $i$. Then it is expected that $z$ arises from a codimension $i$ element $\widetilde{z} \in \textrm{CH}^{i}(X^{\otimes n}_{Z_{\overline{\kappa}}})$ precisely when the crystalline realization of $z$ is the restriction of a global section over $Z$ of the isocrystal $R^{2i} h^{\otimes n}_{*} \mathcal{O}_{X^{\otimes n}_{Z_{\overline{\kappa}}}/K}$, with $K$ the fraction field of the Witt vector ring $W = W(\overline{\kappa})$ and $\mathcal{O}_{X^{\otimes n}_{Z_{\overline{\kappa}}}/K}$ the crystalline structure sheaf; this is the content of the crystalline variational Tate conjecture of \cite{morrow:hal-02385601}. The crystals $\mathcal{H}_{\textrm{cris}}^{a,b}$ of even weight all arise as subobjects of $R^{2i} h^{\otimes n}_{\textrm{cris}, *} \mathcal{O}_{X^{\otimes n}/W(\kappa)}$ for appropriate $i$ and $n$. Moreover, as was explained in \S\ref{crysflatsec} and \S\ref{filcompsec}, global sections of this type give elements of the set $\mathcal{I}_{Z}$ of \autoref{bigexpthm}. Thus, if we combine Morrow's variational Tate conjecture with the crystalline Tate conjecture for the fibre $X^{\otimes n}_{s}$, we obtain a conjectural descriptions of the sets $\mathcal{A}^{a,b}_{\textrm{cris}}$ as consisting of those global sections of $\mathcal{H}^{a,b}_{\textrm{cris},Z_{\overline{\kappa}}}$ whose fibres at points of $Z_{\overline{\kappa}}$ are Frobenius eigenvectors of the expected type.

\vspace{0.5em}

\noindent \textbf{Reducing to Bounded Tensors:} To begin the argument, let us first make explicit the fibre functor $F_{Z}$ alluded to in \S\ref{tensorcasesec}. Fix a point $s : \Spec \overline{\kappa} \to Z_{\overline{\kappa}}$. Then by pullback along $s$ one obtains a functor $F_{Z} : (Z_{\overline{\kappa}}/W(\overline{\kappa}))_{\textrm{cris},K} \to (\Spec \overline{\kappa}/W(\overline{\kappa}))_{\textrm{cris},K}$, and the second category is canonically isomorphic to the category of vector spaces over $K$. If we consider the crystal $\mathcal{H}_{\textrm{cris}, Z_{\overline{\kappa}}}$, then the object $F_{Z}(\mathcal{H}_{\textrm{cris}, Z_{\overline{\kappa}}})$ may be identified with $\mathcal{H}_{\textrm{cris},s} = H^{k}_{\textrm{cris}}(X_{s}/K)$, and admits a crystalline Frobenius endomorphism $\varphi : \mathcal{H}_{\textrm{cris},s} \to \mathcal{H}_{\textrm{cris},s}$. Using the cup product form to identify $\mathcal{H}_{\textrm{cris},s}$ with $\mathcal{H}^{*}_{\textrm{cris},s}$, we may reduce to the situation where either:
\begin{itemize}
\item[-] we have $\mathcal{S} \subset \bigoplus_{j \geq 1} (\mathcal{H}_{\textrm{cris}, Z_{\overline{\kappa}}})^{\otimes j}$ with the fibre at $s$ fixed by $\varphi / p^{k/2}$ (if $k$ is even); or
\item[-] we have $\mathcal{S} \subset \bigoplus_{j \geq 1} (\mathcal{H}_{\textrm{cris}, Z_{\overline{\kappa}}})^{\otimes j} \otimes (\mathcal{H}^{*}_{\textrm{cris}, Z_{\overline{\kappa}}})^{\otimes j}$ with the fibre at $s$ fixed by $\varphi$ (if $k$ is odd).
\end{itemize}

We now make a choice of $N > 0$: we take $N$ to be large enough so that all connected semisimple subgroups of $\GL_{m,F}$, where $F$ is any characteristic zero field, are defined by their invariants in $\bigoplus_{a,b \leq N} F^{a,b}$. That such an $N$ exists follows from the fact that there are only finitely many semisimple subgroups of $\GL_{m,\mathbb{C}}$ up to conjugacy. With such a choice made, we now claim that $\mathcal{S}$ can be taken to be a subset of global sections of $\bigoplus_{a,b \leq N} \mathcal{H}^{a,b}_{\textrm{cris},Z_{\overline{\kappa}}}$. Indeed, as $\mathbf{H}_{\mathcal{S}}$ is assumed semisimple, we may replace $\mathcal{S}$ with (sections corresponding to) its set of invariants in $\bigoplus_{a, b \leq N} \mathcal{H}^{a,b}_{\textrm{cris},s}$. That these invariants correspond to global sections of $\bigoplus_{a,b \leq N} \mathcal{H}^{a,b}_{\textrm{cris},Z}$ is a consequence of the Tannakian formalism, and that they are Frobenius eigenvectors of the expected type, and therefore conjecturally given by algebraic cycles, is a consequence of the fact that either $\varphi / p^{k}$ or $\varphi$ (as appropriate) acts through the group $\mathbf{H}_{\mathcal{S}}$. (Note that this latter fact is true at \emph{every} point $s \in Z_{\overline{\kappa}}$, so if we started with a set $\mathcal{S}$ given by families of algebraic cycles we will still have a collection of sections which are conjecturally algebraic in every fibre.)

\vspace{0.5em}

The above reasoning has reduced the statement of the proposition (in the case where $\kappa$ has positive characteristic), to the situation where $\mathcal{S}$ consists tensors with degree bounded by $N$. We now wish to use $\mathcal{S}$ to obtain global flat sections $\overline{\mathcal{S}}$ of $\restr{\mathcal{H}}{Z_{\overline{\kappa}}}$ which satisfy the hypotheses outlined in \autoref{bigexpthm}. The main difficulty is showing that the stabilizer subgroup $\mathbf{H}_{\overline{\mathcal{S}}}$ of the reduced sections remains semisimple, at least after passing to an open subscheme of $\Spec R$. We therefore begin with a digression on the moduli of semisimple Lie algebras.

\vspace{0.5em}

\noindent \textbf{Semisimple Moduli: } Let us construct a crude moduli space for semisimple Lie subalgebras of $\mathfrak{gl}_{m}$ as follows. We begin by considering the constructible locus $\mathcal{M}$ in the union $\bigcup_{d} \textrm{Gr}(\mathfrak{gl}_{m}, d)$ of Grassmanian schemes which parameterizes all semisimple Lie algebras. Over $\mathbb{C}$, this is a finite union of locally closed orbits of $\GL_{m,\mathbb{C}}$ under the conjugation action. Let $F$ be a number field over which the geometric components of $\mathcal{M}$ are defined, and let $B$ be its ring of integers. After passing to an open subscheme of $\Spec B$, we may give $\mathcal{M}$ the structure of a smooth reduced locally closed subscheme of $\bigcup_{d} \textrm{Gr}(\mathfrak{gl}_{R}, d)_{B}$. Letting $\mathcal{M}_{1}, \hdots, \mathcal{M}_{j}$ be the geometric components of $\mathcal{M}$, then we may choose uniformization maps $\pi_{i} : \GL_{m,B} \to \mathcal{M}_{i}$ by spreading out the orbit maps and passing to an open subscheme of $\Spec B$. Finally, each component $\mathcal{M}_{i,\mathbb{C}}$ consists of semisimple Lie algebras whose associated connected semisimple groups are the stabilizers of a tensor subspace of $\bigoplus_{a, b \leq N} (\mathbb{C}^{m})^{a,b}$ of dimension exactly $e_{i}$; after replacing $B$ with an open subscheme, we may assume this is true for all field-valued points of $\mathcal{M}_{i}$.

As we are free to pass to an open subscheme of $R$, and replace $R$ by a finite integral extension, we now assume that $B \subset R$. Let $L$ be an algebraic closure of $K$, and let $A$ be its associated valuation ring; note that the residue field of $A$ may be identified with $\overline{\kappa}$. We claim that every $L$-point $m$ of $\mathcal{M}_{i}$ is the generic fibre of an $A$-point of $\mathcal{M}_{i}$. To see this, choose an $L$-point $n \in \pi^{-1}(m)$. After replaying $n$ with $\lambda n$ for some scalar $\lambda \in L$, we may assume that $n$ is the restriction of $\widetilde{n} \in \GL_{m,B}(A)$, because the maps $\pi_{i}$ factor through the quotient of $\GL_{m,B}$ by the diagonal subgroup. Then $\pi(\widetilde{n})$ has $m$ as its generic fibre, and lies inside $\mathcal{M}_{i}(A)$.

\vspace{0.5em}

\noindent \textbf{Applying Semisimplicity:} The above discussion has shown that, at least after replacing $R$ with an open subscheme and passing to an integral extension, we may assume the following property holds:

\begin{quote} 
Given any connected semisimple Lie subgroup $H \subset \GL_{m,K}$ defined over $K$, whose set of fixed tensors is $\mathcal{T}_{H} \subset \bigoplus_{a, b \leq N} (K^{m})^{a,b}$, there exists a model $\mathbf{H}$ of $H_{L}$ over $A$ such that the special fibre of $\mathbf{H}$ is a connected semisimple group whose fixed locus in $\bigoplus_{a, b \leq N} (\overline{\kappa}^{m})^{a,b}$ is exactly the reduction to $\overline{\kappa}$ of $\mathcal{T}_{H} \cap \left(\bigoplus_{a, b \leq N} (W^{m})^{a,b} \right)$. 
\end{quote} 
Applying this to the specific case of the (identify component of) the group $\mathbf{H}_{\mathcal{S}}$, we may choose the elements of $\mathcal{S}$ such that their reductions $\overline{\mathcal{S}}$ define a semisimple group $\mathbf{H}_{\overline{\mathcal{S}}}$ in every fibre of $\restr{\mathcal{H}}{Z_{\overline{\kappa}}}$. 

Because the sections $\overline{\mathcal{S}}$ lie inside $\mathcal{I}_{Z,N}$ it suffices to check that the adjoint Hodge numbers on the Lie algebra of $\mathbf{H}_{\overline{\mathcal{S}}}$ are symmetric. This follows from \autoref{semsimpimpsymlem} after possibly passing to an open subscheme of $\Spec R$.

\vspace{0.5em}

\noindent \textbf{Characteristic Zero:} Finally, we consider the case where $\kappa(\mathfrak{p})$ has characteristic zero. For any such $Z$, the complexification $Z_{\mathbb{C}}$ is contained in the (tensorial) Hodge locus of $S_{\mathbb{C}}$. Since any component of the Hodge locus is a weakly special subvariety by \cite[Def 4.1, Cor. 4.14]{2021InMat.225..857K}, necessarily such a $Z$ lies inside $E$ by the first property of $E$ listed in \autoref{bigexpthm}.
\end{proof}

\noindent Let us also note that the fact that \autoref{properlocusthm2} implies \autoref{properlocusthm} may be argued exactly as in the beginning of the ``reduction'' step in the proof of \autoref{bigexpthm}, where we saw that the quasi-finiteness of $\varphi$ can be used to imply that $r$-limps are immersive away from a closed subscheme of $S$; as there are no new ideas here we have omitted this reduction.

Lastly, we resolve:

\begin{proof}[Proof of \autoref{compinstcor}:]
Observe that the hypotheses of \autoref{properlocusthm2} are satisfied for the family $f$: all that needs to be checked is the Zariski density of the monodromy representation, which is due to \cite{beauville1986groupe}. (The exceptional cases described in \cite{beauville1986groupe} are ruled out by our assumption on level.) Since the non-ordinary locus in $S_{\mathbb{F}_{p}}$ is not Zariski dense due to a result of Illusie \cite{Illusie2007}, this shows the first claim. 

Restricting now to the case where $m = 1$, it suffices to check that after restricting to some dense open locus inside $S_{\mathbb{F}_{p}}$ the condition that all $r$-limps on $Z$ are constant implies that $Z$ lies inside a $\GL_{M,\mathbb{F}_{p}}$ orbit. As in the proof of \autoref{bigexpthm}, this can be established at the integral level (away from finitely many primes), by observing that the canonical period map $\varphi : \an{S_{\mathbb{C}}} \to \Gamma \backslash D$ factors through a period map $\an{(S/\GL_{M})_{\mathbb{C}}} \to \Gamma \backslash D$ which is generically injective as a consequence of the generic Torelli theorem for the hypersurface variation proved in \cite{CM_1983__50_2-3_325_0}.
\end{proof}

\bibliography{hodge_theory}

\begin{thebibliography}{CPMS03}

\bibitem[BBT18]{OMINGAGA}
Benjamin {Bakker}, Yohan {Brunebarbe}, and Jacob {Tsimerman}.
\newblock {o-minimal GAGA and a conjecture of Griffiths}.
\newblock {\em arXiv e-prints}, page arXiv:1811.12230, November 2018.

\bibitem[Bd11]{2011arXiv1110.5001B}
Bhargav {Bhatt} and Aise~Johan {de Jong}.
\newblock {Crystalline cohomology and de Rham cohomology}.
\newblock {\em arXiv e-prints}, page arXiv:1110.5001, October 2011.

\bibitem[Bea86]{beauville1986groupe}
Arnaud Beauville.
\newblock Le groupe de monodromie des familles universelles d'hypersurfaces et
  d'intersections compl{\`e}tes.
\newblock In {\em Complex analysis and algebraic geometry}, pages 8--18.
  Springer, 1986.

\bibitem[Ber74]{Berthelot1974}
Pierre Berthelot.
\newblock {\em Cohomologie cristalline et cohomologie de de Rham}.
\newblock Springer Berlin Heidelberg, Berlin, Heidelberg, 1974.

\bibitem[BKU21]{BKU}
Gregorio Baldi, Bruno Klingler, and Emmanuel Ullmo.
\newblock On the distribution of the {H}odge locus.
\newblock {\em arXiv preprint arXiv:2107.08838}, 2021.

\bibitem[BT17]{AXSCHAN}
Benjamin Bakker and Jacob Tsimerman.
\newblock The {A}x-{S}chanuel conjecture for variations of {H}odge structures.
\newblock {\em Inventiones mathematicae}, 217:77--94, 2017.

\bibitem[CPMS03]{CMS}
James~A Carlson, C.~(Chris) Peters, and Stefan M\"uller-Stach.
\newblock {\em Period mappings and period domains}.
\newblock Cambridge, U.K. : Cambridge University Press, 2003.

\bibitem[Del71]{PMIHES_1971__40__5_0}
Pierre Deligne.
\newblock Th\'eorie de {Hodge} : {II}.
\newblock {\em Publications Math\'ematiques de l'IH\'ES}, 40:5--57, 1971.

\bibitem[DK00]{Duistermaat2000}
J.~J. Duistermaat and J.~A.~C. Kolk.
\newblock {\em Compact {L}ie {G}roups}.
\newblock Springer Berlin Heidelberg, Berlin, Heidelberg, 2000.

\bibitem[Don83]{CM_1983__50_2-3_325_0}
Ron Donagi.
\newblock Generic torelli for projective hypersurfaces.
\newblock {\em Compositio Mathematica}, 50(2-3):325--353, 1983.

\bibitem[GGK12]{GGK}
Mark Green, Phillip Griffiths, and Matt Kerr.
\newblock {\em Mumford-Tate Groups and Domains: Their Geometry and Arithmetic
  (AM-183)}.
\newblock Princeton University Press, 2012.

\bibitem[Har75]{PMIHES_1975__45__5_0}
Robin Hartshorne.
\newblock On the de {Rham} cohomology of algebraic varieties.
\newblock {\em Publications Math\'ematiques de l'IH\'ES}, 45:5--99, 1975.

\bibitem[How89]{invtheoryrem}
Roger Howe.
\newblock Remarks on {C}lassical {I}nvariant {T}heory.
\newblock {\em Transactions of the American Mathematical Society},
  313(2):539--570, 1989.

\bibitem[hp]{257785}
Alexander~Premet (https://mathoverflow.net/users/24386/alexander premet).
\newblock {S}emisimplicity of {L}ie algebra in positive characteristic.
\newblock MathOverflow.
\newblock URL:https://mathoverflow.net/q/257785 (version: 2016-12-31).

\bibitem[Ill07]{Illusie2007}
Luc Illusie.
\newblock Ordinarit{\'e} des intersections compl{\`e}tes g{\'e}n{\'e}rates.
\newblock In Pierre Cartier, Nicholas~M. Katz, Yuri~I. Manin, Luc Illusie,
  G{\'e}rard Laumon, and Kenneth~A. Ribet, editors, {\em The Grothendieck
  Festschrift: A Collection of Articles Written in Honor of the 60th Birthday
  of Alexander Grothendieck}, pages 375--405. Birkh{\"a}user Boston, Boston,
  MA, 2007.

\bibitem[Kat70]{katznil}
Nicholas Katz.
\newblock Nilpotent connections and the monodromy theorem : applications of a
  result of {Turrittin}.
\newblock {\em Publications Math\'ematiques de l'IH\'ES}, 39:175--232, 1970.

\bibitem[Kat72]{Katz1972}
Nicholas Katz.
\newblock Algebraic {S}olutions of {D}ifferential equations (p-curvature and
  the {H}odge {F}iltration).
\newblock {\em Inventiones mathematicae}, 18:1--118, 1972.

\bibitem[Kaw81]{CM_1981__43_2_253_0}
Yujiro Kawamata.
\newblock Characterization of abelian varieties.
\newblock {\em Compositio Mathematica}, 43(2):253--276, 1981.

\bibitem[KO68]{katz1968}
Nicholas~M. Katz and Tadao Oda.
\newblock On the differentiation of {D}e {R}ham cohomology classes with respect
  to parameters.
\newblock {\em J. Math. Kyoto Univ.}, 8(2):199--213, 1968.

\bibitem[KO21]{2021InMat.225..857K}
B.~{Klingler} and A.~{Otwinowska}.
\newblock {On the closure of the {H}odge locus of positive period dimension}.
\newblock {\em Inventiones Mathematicae}, 225(3):857--883, September 2021.

\bibitem[Kos10]{Kostant2010}
Bertram Kostant.
\newblock Root systems for {L}evi {F}actors and {B}orel--de {S}iebenthal
  {T}heory.
\newblock In H.~E.~A. Campbell, Aloysius~G. Helminck, Hanspeter Kraft, and
  David Wehlau, editors, {\em Symmetry and Spaces: In Honor of Gerry Schwarz},
  pages 129--152. Birkh{\"a}user Boston, Boston, 2010.

\bibitem[KOU20]{fieldsofdef}
Bruno {Klingler}, Anna {Otwinowska}, and David {Urbanik}.
\newblock {On the fields of definition of Hodge loci}.
\newblock {\em To appear in Ann. Sci. de l'École Nor. Sup. arXiv:2010.03359},
  October 2020.

\bibitem[Maz73]{10.2307/1970906}
B.~Mazur.
\newblock Frobenius and the {H}odge {F}iltration (estimates).
\newblock {\em Annals of Mathematics}, 98(1):58--95, 1973.

\bibitem[Mor19]{morrow:hal-02385601}
Matthew Morrow.
\newblock {A Variational Tate Conjecture in crystalline cohomology}.
\newblock {\em {Journal of the European Mathematical Society}},
  21(11):3467--3511, 2019.

\bibitem[Sch73]{Schmid1973}
Wilfried Schmid.
\newblock Variation of {H}odge {S}tructure: {T}he {S}ingularities of the
  {P}eriod {M}apping.
\newblock {\em Inventiones mathematicae}, 22:211--320, 1973.

\bibitem[{Sta}20]{stacks-project}
The {Stacks project authors}.
\newblock The stacks project.
\newblock \url{https://stacks.math.columbia.edu}, 2020.

\bibitem[SZ85]{Steenbrink1985}
Joseph Steenbrink and Steven Zucker.
\newblock Variation of mixed {H}odge structure. {I}.
\newblock {\em Inventiones mathematicae}, 80:489--542, 1985.

\bibitem[Urb21a]{periodimages}
David Urbanik.
\newblock On the {T}ranscendence of {P}eriod {I}mages.
\newblock {\em arXiv preprint arXiv:2106.09342}, 2021.

\bibitem[Urb21b]{urbanik2021sets}
David Urbanik.
\newblock Sets of {S}pecial {S}ubvarieties of {B}ounded {D}egree.
\newblock {\em arXiv preprint arXiv:2109.07663}, 2021.

\bibitem[Voi10]{voisin2010hodge}
Claire Voisin.
\newblock Hodge loci.
\newblock {\em Handbook of moduli}, 3:507--546, 2010.

\bibitem[Yve92]{Andre1992}
Andr\'e Yves.
\newblock Mumford-{T}ate groups of mixed {H}odge structures and the theorem of
  the fixed part.
\newblock {\em Compositio Mathematica}, 82(1):1--24, 1992.

\end{thebibliography}
\bibliographystyle{alpha}

\end{document}